\documentclass[11pt,leqno]{amsart}
\usepackage[usenames,dvipsnames,svgnames]{xcolor}
\date{}
\usepackage{hyperref}
\usepackage{amssymb,amsthm,amsmath,amscd,amsfonts,bbm,mathrsfs,fullpage,comment}
\usepackage[all]{xy}
\DeclareMathAlphabet{\mathpzc}{OT1}{pzc}{m}{it}
\DeclareMathAlphabet{\mathdutchcal}{U}{dutchcal}{m}{n}

\newcommand{\addappendix}{%
  \section*{\appendixname}% start the appendix
  \stepcounter{section}% reset counters related to section
  \renewcommand{\thesection}{A}% we want A
 \stepcounter{equation}
  }

 \def\noproof{{\unskip\nobreak\hfill\penalty50\hskip2em\hbox{}%
      \nobreak\hfill$\Box$\parfillskip=0pt%
     \finalhyphendemerits=0\par}}
\def\enddemo{\ifmmode\eqno\Box\else\noproof\vskip0.8truecm\fi}
\newtheorem{theo}{Theorem}[section]
\newtheorem{theorem}[theo]{Theorem}
\newtheorem{df}[theo]{Definition}

\newtheorem{lemma/def}[theo]{Lemma/Definition}

\newtheorem{prop}[theo]{Proposition}
\newtheorem{coro}[theo]{Corollary}

\newtheorem{remark}[theo]{Remark}
\newtheorem{remarks}[theo]{Remarks}

\newtheorem{lemma}[theo]{Lemma}

\newtheorem{examples}[theo]{Examples}
\newtheorem{example}[theo]{Example}

\newcommand{\lra}{\longrightarrow}
\newcommand{\lla}{\longleftarrow}
\newcommand{\hra}{\hookrightarrow}

\DeclareMathOperator{\Aff}{Aff}
\DeclareMathOperator{\pt}{pt}

\DeclareMathOperator{\Spec}{Spec}

\DeclareMathOperator{\Eis}{Eis}
\DeclareMathOperator{\bEis}{\mathbb{E}is}

\DeclareMathOperator{\Cov}{Cov}

\DeclareMathOperator{\Ind}{Ind}

\DeclareMathOperator{\aug}{aug}

\DeclareMathOperator{\inv}{inv}

\DeclareMathOperator{\res}{res}
\DeclareMathOperator{\cor}{cor}

\DeclareMathOperator{\rec}{rec}

\DeclareMathOperator{\ord}{ord}

\DeclareMathOperator{\id}{id}

\DeclareMathOperator{\ev}{ev}

\DeclareMathOperator{\Aut}{Aut}
\DeclareMathOperator{\Hom}{Hom}

\DeclareMathOperator{\univ}{univ}

\DeclareMathOperator{\uM}{\underline{M}}

\DeclareMathOperator{\Ext}{Ext}

\DeclareMathOperator{\pr}{pr}

\DeclareMathOperator{\rank}{rank}

\DeclareMathOperator{\supp}{supp}
\DeclareMathOperator{\tor}{tor}
\DeclareMathOperator{\wcdot}{\, \cdot\, }

\DeclareMathOperator{\Maps}{Maps}
\DeclareMathOperator{\Real}{Re}

\DeclareMathOperator{\opp}{opp}

\DeclareMathOperator{\sh}{sh}

\DeclareMathOperator{\Gal}{Gal}
\DeclareMathOperator{\Norm}{N}

\DeclareMathOperator{\Mod}{Mod}

\DeclareMathOperator{\spaces}{\underline{Top}}

\DeclareMathOperator{\open}{open}

\DeclareMathOperator{\GL}{GL}

\DeclareMathOperator{\Ann}{Ann}
\DeclareMathOperator{\disc}{disc}
\DeclareMathOperator{\incl}{incl}

\DeclareMathOperator{\fin}{fin}

\DeclareMathOperator{\pol}{pol}
\DeclareMathOperator{\loc}{loc}
\DeclareMathOperator{\mloc}{--loc}
\DeclareMathOperator{\lpol}{lpol}
\DeclareMathOperator{\mlpol}{--lpol}

\DeclareMathOperator{\Int}{Int}
\DeclareMathOperator{\Lat}{\mathdutchcal{Lat}}

\DeclareMathOperator{\cycl}{cycl}

\DeclareMathOperator{\prolim}{\underset{\lla}{\lim}}

\DeclareMathOperator{\prolimm}{\underset{\underset{m}{\lla}}{\lim}}
\DeclareMathOperator{\dlim}{\underset{\lra}{\lim}}

\DeclareMathOperator{\sgn}{sign}

\newcommand{\io}{{\iota}}
\newcommand{\La}{{\Lambda}}
\newcommand{\la}{{\lambda}}
\newcommand{\Ga}{{\Gamma}}
\newcommand{\ga}{{\gamma}}
\newcommand{\Si}{{\Sigma}}

\newcommand{\fa}{{\mathfrak a}}
\newcommand{\fA}{{\mathfrak A}}
\newcommand{\fb}{{\mathfrak b}}
\newcommand{\fB}{{\mathfrak B}}
\newcommand{\fc}{{\mathfrak c}}
\newcommand{\fC}{{\mathfrak C}}

\newcommand{\fD}{{\mathfrak D}}
\newcommand{\ff}{{\mathfrak f}}

\newcommand{\fm}{{\mathfrak m}}

\newcommand{\fn}{{\mathfrak n}}
\newcommand{\fp}{{\mathfrak p}}

\newcommand{\fq}{{\mathfrak q}}
\newcommand{\fr}{{\mathfrak r}}

\newcommand{\fX}{{\mathfrak X}}
\newcommand{\barfC}{{\overline{\fC}}}
\newcommand{\wfC}{{\widetilde{\fC}}}
\newcommand{\hatfa}{{\widehat{\fa}}}
\newcommand{\wfa}{{\widetilde{\fa}}}

\newcommand{\fOpen}{{\mathfrak{Open}}}

\newcommand{\bC}{{\mathbb C}}

\newcommand{\bF}{{\mathbb F}}

\newcommand{\bN}{{\mathbb N}}

\newcommand{\bQ}{{\mathbb Q}}
\newcommand{\bR}{{\mathbb R}}

\newcommand{\bZ}{{\mathbb Z}}
\newcommand{\barQ}{{\overline{\mathbb Q}}}

\newcommand{\bA}{{\mathbb A}}

\newcommand{\bhatZ}{{\widehat{\bZ}}}

\newcommand{\barbF}{{\overline{\bF}}}

\newcommand{\barK}{{\overline{K}}}

\newcommand{\barR}{{\overline{R}}}

\newcommand{\barX}{{\overline{X}}}

\DeclareMathOperator{\bga}{\overline{\gamma}}

\newcommand{\cA}{{\mathcal A}}
\newcommand{\cB}{{\mathcal B}}
\newcommand{\cC}{{\mathcal C}}
\newcommand{\cD}{{\mathcal D}}

\newcommand{\cF}{{\mathcal F}}
\newcommand{\cG}{{\mathcal G}}
\newcommand{\cH}{{\mathcal H}}
\newcommand{\cI}{{\mathcal I}}
\newcommand{\cJ}{{\mathcal J}}

\newcommand{\cM}{{\mathcal M}}
\newcommand{\cN}{{\mathcal N}}
\newcommand{\cO}{{\mathcal O}}
\newcommand{\cP}{{\mathcal P}}

\newcommand{\cR}{{\mathcal R}}
\newcommand{\cS}{{\mathcal S}}

\newcommand{\cU}{{\mathcal U}}
\newcommand{\cV}{{\mathcal V}}

\newcommand{\cX}{{\mathcal X}}
\newcommand{\cY}{{\mathcal Y}}
\newcommand{\Sh}{{\mathcal S}h}
\newcommand{\PSh}{{\mathcal PS}h}

\newcommand{\cOs}{{\cO_{\Sigma}}}
\newcommand{\cOsu}{\cO_{\Sigma}^{\,*}}

\newcommand{\cOscu}{\cO_{\Sigma, \fc}^{\,*}}

\newcommand{\sD}{{\mathscr D}}

\newcommand{\sF}{{\mathscr F}}
\newcommand{\sG}{{\mathscr G}}
\newcommand{\sH}{{\mathscr H}}
\newcommand{\sI}{{\mathscr I}}
\newcommand{\sJ}{{\mathscr J}}

\newcommand{\sL}{{\mathscr L}}

\DeclareMathOperator{\sLog}{{\sL og}}

\DeclareMathOperator{\brC}{\langle C \rangle}

\newcommand{\wGa}{\widetilde{\Gamma}}
\newcommand{\wDelta}{\widetilde{\Delta}}

\newcommand{\wM}{\widetilde{M}}

\newcommand{\wcR}{\widetilde{\cR}}

\newcommand{\wga}{{\widetilde{\ga}}}

\newcommand{\hA}{\widehat{A}}
\newcommand{\hB}{\widehat{B}}
\newcommand{\hT}{\widehat{T}}
\newcommand{\hL}{\widehat{L}}

\newcommand{\hX}{\widehat{X}}
\newcommand{\hY}{\widehat{Y}}

\newcommand{\hLa}{\widehat{\Lambda}}

\newcommand{\hf}{\hat{f}}

\newcommand{\tf}{\tilde{f}}

\newcommand{\ua}{{\underline{a}}}

\newcommand{\uk}{{\underline{k}}}
\newcommand{\ul}{{\underline{l}}}
\newcommand{\um}{{\underline{m}}}

\newcommand{\ut}{{\underline{t}}}

\newcommand{\uX}{{\underline{X}}}
\newcommand{\uY}{{\underline{Y}}}
\newcommand{\uz}{{\underline{z}}}

\newcommand{\ep}{{\epsilon}}
\newcommand{\vep}{{\varepsilon}}

\newcommand{\ula}{{\underline{\lambda}}}

\newcommand{\uxi}{{\underline{\xi}}}

\newcommand{\bu}{{\bullet}}

\newcommand{\noi}{\noindent}

\begin{document}

\title{Adelic Eisenstein classes and divisibility properties of Stickelberger elements}

\author{Alexandros Galanakis and Michael Spie{\ss}}
\thanks{Both authors acknowledge support by the Deutsche
    Forschungsgemeinschaft (DFG, German Research Foundation) via the grant
    SFB-TRR 358/1 2023 — 491392403}
\address{Universität Bielefeld, Fakultät für Mathematik, Germany,
E-Mail: \url{alexandros.galanakis@gmail.com}}
\address{Universität Bielefeld, Fakultät für Mathematik, Germany,
E-Mail: \url{mspiess@math.uni-bielefeld.de}}
\date{}
\maketitle

\begin{abstract}
Nori's Eisenstein cohomology classes and their integral refinements due to Beilinson, Kings and Levin can be used to obtain simple proofs of the rationality and integrality properties of special values of abelian $L$-functions of totally real fields. Here we introduce an adelic refinement of these constructions. This will be used to establish new divisibility properties of Stickelberger elements associated to abelian extensions of totally real fields.
\end{abstract}

\tableofcontents

\section{Introduction}

In this paper we introduce an adelic refinement of the Eisenstein cohomology classes introduced by
Beilinson, Kings and Levin in \cite{bkl} (which in turn are an integral refinement of the Eisenstein classes introduced by Nori \cite{nori}). The classes of Beilinson, Kings and Levin are elements in the cohomology of $\GL_n(\bZ)$ -- or some of its subgroups -- in degree $n\!-\!1$ with coefficients in the completed group ring of the free-abelian group $\bZ^n$. Their construction is given in terms of $\GL_n(\bZ)$-equivariant sheaf cohomology on the $n$-dimension real torus $\bR^n/\bZ^n$ with coefficients in the so called {\it logarithm sheaf} $\sLog$, the locally constant sheaf associated to the completed group ring. 

The adelic Eisenstein classes introduced in this paper are cohomology classes of the group $\GL_n(\bQ)$ in degree $n\!-\!1$  with coefficients in a module $\cD_{\lpol}$ which we call the module of {\it locally polynomial distributions on $\bQ^n$.}\footnote{More precisely for technical reasons we have to work with a proper ''large'' subgroup $\Gamma\subseteq \GL_n(\bQ)$ and distributions on a certain $\Gamma$-stable subset of $\bQ^n$; for the purpose of keeping this introduction less technical we divert to section \ref{section:eiscocycle} for details.} In principle our construction is modelled after that introduced in \cite{bkl}. We consider the $\GL_n(\bQ)$-equivariant sheaf cohomology of the adelic solenoid $(\bA/\bQ)^n$. However one of the difficulties in generalizing the approach of \cite{bkl} to the adelic setting is that the usual topology of the adeles and the adele class group is too fine to produce meaningful sheaf cohomology groups. Instead we introduce a certain coarse Grothendieck topology -- which we call {\it lattice topology} -- and develop a $\GL_n(\bQ)$-equivariant sheaf cohomology with respect to this topology. Another difficulty that we encountered is the fact that the logarithm sheaf does not seem to admit a natural adelic counterpart. We work instead with the {\it sheaf of locally polynomial distributions $\sD$}.

The interest in studying Eisenstein cohomology classes lies in their relation to special values of $L$-functions. This was first exploited by Nori (and independently by Sczech) who reproved the Theorem of Siegel and Klingen regarding the rationality of the values of partial zeta functions of totally real fields at non-positive integers. Using their integral refinement of Nori's classes Beilinson, Kings and Levin were able to reprove the integrality results of Cassou-Nogu{\`e}s and Deligne-Ribet for these special $L$-values. Our adelic variants of the Eisenstein classes allow us to refine these results further by proving certain divisibility properties for Stickelberger elements at non-positive integers. A first result in this direction had been obtained in joint work of the second author with S.\ Dasgupta \cite{dassp}. There a different and rather concrete construction of Eisenstein cohomology classes -- based on Shintani's method to study special $L$-values -- had been used. We feel that the construction using equivariant sheaf cohomology introduced here is conceptually particularly satisfying and that it adds a new tool to study properties of special $L$-values. The relation between the Eisenstein classes introduced here and in \cite{dassp} remains unclear.\footnote{It should be noted that both construction differ insofar that the two methods produce classes in different cohomology groups, i.e.\ in cohomology groups for different ''large'' subgroups $\Gamma$ of $\GL_n(\bQ)$. It turns out that one obtains also slightly different divisibility properties for Stickelberger elements; see Remark \ref{remarks:history} (b) below for details.}

We now give the reader some idea about the topological (i.e.\ topos-theoretic) constructions on which our definition of the adelic Eisenstein classes is based. For the purpose of keeping the technical details to a minimum we consider here only a somewhat simplified situation and refer to sections \ref{section:locpol}, \ref{section:latspaces} and \ref{section:eiscocycle} for the general framework. We fix a totally real number field $F\subseteq \bC$ of degree $n$ over $\bQ$ with ring of intergers $\cO_F$. Let $\cI$ denote the set of all fractional ideals of $F$. We attach to the ring of finite adeles $\bA_f$, the ring of adeles $\bA$ and the adele class group $\bA/F$ of $F$ certain sites\footnote{We use here a definition of a site that is not the standard notion familiar to researchers in algebraic or arithmetic geometry. For this reason we review  in the appendix the necessary background material regarding sites and topoi that we use throughout this paper.} denoted by $B$, $A$ and $T$. We refer to these as {\it lattice topologies} on $\bA_f$, $\bA$ and $\bA/F$ respectively. They are much coarser than the usual topologies. For example an open subset $U\subseteq B$ (i.e.\ an object of $B$) is a subset $U\subseteq \bA_f$ that is $\hatfa$-stable for some fractional ideal $\fa\in \cI$ (i.e.\ we have $\hatfa + U=U$). Here $\hatfa := \fa\otimes \bhatZ$ is the closure of $\fa$ in $\bA_f$. A covering of an open subset $U\subset B$ consists of a collection of open subsets $\{U_i\}_i$ of $B$ so that $U=\bigcup_i U_i$ and so that there exists $\fa\in \cI$ such that each $U_i$ is $\hatfa$-stable. The lattice topologies $A$ and $T$ are defined in a similar way. 

The group $\wGa:= \Aff(F)= F^*\ltimes F$ acts continuously on $B$ and $A$ whereas $\Ga := F^*$ acts continuously on $T$. Moreover the site $T$ can be identified with the ''quotient'' of $A$ with respect to the natural projection $\pr: \bA\to \bA/F$. Hence the pull-back $\pr^*: \Sh(T, \Ga) \to \Sh(A, \wGa)$ defines an equivalence between the category of $\wGa$-equivariant sheaves on $A$ and the category of $\Ga$-equivariant sheaves on $T$. The basic properties of (equivariant) sheaf cohomology for sites of the type $B$, $A$ and $T$ will be established in section \ref{section:latspaces}. 

The construction of the $\Gamma$-equivariant sheaf $\sD$ on $T$ is based on the notion of a locally polynomial function which we now review. Let $U\subseteq B$ be an open subset and put $\cU := U\cap F$. A map $f: \cU \to \bZ$ is called a {\it locally polynomial function} if there exists $\fa\in \cI$ such that $U$ is $\hatfa$-stable and so that for a $\bZ$-basis $(\omega_1, \ldots, \omega_n)$ of $\fa$ and for every $x\in \cU$ the map $\bZ^n \to \bZ, (x_1, \ldots, x_n) \mapsto f(x+ \sum_{i=1}^n x_i \omega_i)$ is a polynomial function, i.e.\ there exists a polynomial $P\in \bQ[X_1, \ldots, X_n]$ with $f(x+ \sum_{i=1}^n x_i \omega_i)=P(x_1, \ldots, x_n)$ for every $x_1, \ldots, x_n\in \bZ$. The  locally polynomial function $f: \cU \to \bZ$ is said to be of {\it bounded support} if there exists $\fa\in \cI$ such that $\cU$ is $\fa$-stable and so that the support of $f$ is contained in finitely many $\fa$-orbits in $\cU$. We denote by $\Int_{\lpol, b}(\cU, \bZ)$ the ring of locally polynomial functions on $\cU$ of bounded support and put $\cD_{\lpol}(U) := \Hom(\Int_{\lpol, b}(\cU, \bZ), \bZ)$. The assignment $U\mapsto \cD_{\lpol, B}(U)$ defines a $\wGa$-equivariant sheaf on $B$, the sheaf of {\it locally polynomial distributions on $B$}. By pulling it back to $A$ via the morphism $\pr_B: A\to B$ induced by the natural projection $\bA= \bA_f\times F_{\infty}\to \bA_f$ we obtain a $\wGa$-equivariant sheaf $\cD_{\lpol, A}=(\pr_B)^*(\cD_{\lpol, A})$ on $A$. As alluded to above the sheaf $\cD_{\lpol, A}$ decends to a $\Ga$-equivariant sheaf $\sD$ on $T$, i.e.\ we have $\pr^*(\sD) = \cD_{\lpol, A}$.

In section \ref{section:eiscocycle} we investigate the ($\Gamma$-equivariant) cohomology groups of $T$ with coefficients in $\sD$. One of our key results (cf.\ Prop.\ \ref{prop:cohomld}) is that the cohomology $H^{\bu}(T, \sD)$ is concentrated in degree $n$ and that we have $H^n(T, \sD)\cong \bZ$. We then proceed in defining the adelic Eisenstein classes following the blueprint of Beilinson, Kings and Levin. We also establish a direct connection between our classes and theirs (cf.\ Prop.\ \ref{prop:eisencompare}). This will be used in section \ref{section:partzeta} to link the adelic Eisenstein classes to special values of partial zeta functions and to Stickelberger elements.

We now describe our main application.\footnote{It should also be possible to use our adelic Eisenstein classes to give formulas for Brumer-Stark units similar to those obtained using Shintani cocycles (see \cite{dassp}, \S 6 or \cite{dashon}); this question however will not be addressed in this paper.}
 Namely, in section \ref{section:lvalues} we prove certain divisibility properties of Stickelberger elements. For that we fix a finite abelian extension $K/F$, $K\subseteq \bC$ with Galois group $G$. Let $S$ be a finite set of nonarchimedean places of $F$ that contains all places that are ramified in $K$. Recall that the 
partial zeta function associate to an element $\sigma\in G$ is given by
\begin{equation}
\label{partzeta2}
\zeta_S(\sigma, s)\,\, =\,\,\sum_{(\fa,S)=1, \sigma_{\fa}=\sigma}\, \Norm(\fa)^{-s}
\end{equation} 
for $\Real(s) >1$ where $\sigma_{\fa}$ is the image of the ideal $\fa$ under the Artin map. It admits a meromorphic continuation to the whole complex plane with a single simple pole at $s=1$. We package the partial zeta functions into a $\bC[G]$-valued function $\Theta_S(K/F, s)$ -- the Stickelberger element -- defined by
\begin{equation}
\label{stick}
\Theta_S(K/F, s)\, = \, \sum_{\sigma\in G} \zeta_S(\sigma,s) [\sigma^{-1}].
\end{equation} 
By the Theorem of Siegel and Klingen we have $\Theta_S(K/F, -k)\in \bQ[G]$ for $k\in \bZ_{\ge 0}$. It is well-known that in order to obtain integrality results it is necessary to consider a variant of \eqref{stick}, namely the {\it $T$-smoothed Stickelberger element} $\Theta_{S, T}(K/F, s)$ defined by 
\begin{equation}
\label{stick2}
\Theta_{S, T}(K/F, s)\, =\,\prod_{\fq\in T} (1-\Norm(\fq)^{1-s} [\sigma_{\fq}^{-1}]) \Theta_S(K/F, s).
\end{equation} 
Here $T$ is an additional set of nonarchimedean places of $F$ that is disjoint of $S$. Under certain mild conditions on $T$ we have $\Theta_{S, T}(K/F, -k)\in \bZ[G]$ for $k\in \bZ_{\ge 0}$.\footnote{See the remark following Thm.\ \ref{theorem:highervanishing}.}

To state our result we will introduce some ideals in the group ring $\bZ[G]$. For a place $v$ of $F$ let $G_v\subseteq G$ be the decomposition group of $v$. If $v$ is nonarchimedean then we denote by $I_v\subseteq G_v$ the inertia group at $v$ and let $\sigma_v\in G_v/I_v$ be the Frobenius at $v$. If $v$ is archimedean then $\sigma_v$ denotes the generator of $G_v$. For $k\in \bZ_{\ge 0}$ and $v\in S\cup S_{\infty}$ we define ideals $\cI_v^{(k)}\subseteq \bZ[G]$ by
\begin{equation*}
\label{stickidealv}
\cI_v^{(k)}\,:= \left\{ \begin{array}{cc} \ker\left(\bZ[G]\to \bZ[G/I_v]/([\sigma_v^{-1}]-\Norm(v)^k) \right) & \text{if $v\nmid \infty$,}\\
([\sigma_v] + (-1)^{k+1})\, \bZ[G] & \text{if $v\mid \infty$.}
\end{array}\right.
\end{equation*}
Note that for $k=0$ we have $\cI_v^{(0)}= \ker(\bZ[G]\to \bZ[G/G_v])$ for every $v\in S\cup S_{\infty}$. 

We also consider the following obvious map between Galois cohomology groups 
\begin{equation}
\label{galoisres}
H^0(K, \bQ/\bZ(k+1)) \lra \bigoplus_{v\in T_K} H^0(k(v), (\bQ/\bZ)'(k+1))
\end{equation}
where $T_K$ denotes the set of places of $K$ that lie above a place in $T$ and where $k(v)$ is the residue field of $v\in T_K$. The Galois module $\bQ/\bZ(k+1)$ appearing in the source of \eqref{galoisres} is the group $\bQ/\bZ$ with the action of the absolute Galois group $\cG_K = \Gal(\barK/K)$ of $K$ given by the $(k\!+\!1)$-th power of the cyclotomic character $\chi_{\cycl}: \cG_K\to \Aut(\mu(\barQ))\cong \bhatZ^*$.\footnote{For any field $E$ we denote by $\mu(E)$ the roots of unity of $E$.} Similarly, for a place $v\in T_K$ with residue characteristic $p$, the Galois modules $(\bQ/\bZ)'(k+1)$ appearing in the summand $H^0(k(v), (\bQ/\bZ)'(k+1))$ of the target of \eqref{galoisres} 
is the group $(\bQ/\bZ)' : = \{x\in \bQ/\bZ\mid \gcd(\ord(x), p)=1\, \}= \bigoplus_{\ell\ne p} \bQ_{\ell}/\bZ_{\ell}$ where the action of the absolute Galois group $\cG_{k(v)}=\Gal(\barbF_p/k(v))$ is again given by the 
$(k\!+\!1)$-th power of the cyclotomic character $\chi_{\cycl}: \cG_{k(v)}\to \Aut(\mu(\barbF_p))= \Aut(\barbF_p^*)\cong (\bhatZ')^*$ (with $\bhatZ':=\prod_{\ell\ne p} \bZ_{\ell}$). 

Note that if $k=0$ then \eqref{galoisres} is the obvious map $\mu(K) \to \bigoplus_{v\in T_K} \mu(k(v))$. Our main result is 

\begin{theorem}
\label{theorem:highervanishing}
Let $k\in \bZ_{\ge 0}$ and let $\fp\in S$ be a fixed place. Let $T$ be a finite set of nonarchimedean places of $F$ disjoint from $S$ such that the map \eqref{galoisres} is injective. Then we have
\begin{equation}
\label{highvanishing2}
\Theta_{S, T}(K/F, -k)\in \prod_{v\in S \cup S_{\infty}, v\ne \fp} \cI_v^{(k)}.
\end{equation}
\end{theorem}

Hence under the assumption that \eqref{galoisres} is injective we have in particular $\Theta_{S, T}(K/F, -k)\in \bZ[G]$. 

\begin{remarks}
\label{remarks:history} \rm (a) The map \eqref{galoisres} is injective if $T$ contains two primes of different residue characteristics. It is also 
injective if $k=0$ and if $T$ contains one prime of residue characteristic larger than $n + 1$. 
\medskip

\noi (b) Theorem \ref{theorem:highervanishing} is already known in the case $k=0$. It has first been proved in \cite{dassp}, except there it was shown -- under certain mild assumption on $T$ -- that $\Theta_{S, T}(K/F, 0)$ is contained in $\prod_{v\in S \cup S_{\infty}, v\ne v_0} \cI_v^{(0)}$ where $v_0$ is a fixed archimedean place of $F$. Hirose \cite{hirose} later obtained the slightly stronger result \eqref{highvanishing2}. In fact in his work the ''exceptional place'' $\fp$ can be any element of $S\cup S_{\infty}$. It should be noted that in both papers \cite{dassp} and \cite{hirose} the ''Shintani method'' of constructing a degree $(n\!-\!1)$ Eisenstein cohomology class is used. 
\medskip

\noi By work of Burns \cite{burns} it is known that \eqref{highvanishing2} (for $k=0$) can also be deduced from a special case of the equivariant Tamagawa number conjecture. The latter has been established recently in \cite{dks} if $K$ is a CM field. 
\medskip

\noi (c) For $k>0$ the only previously known result towards \eqref{highvanishing2} is the assertion that $\Theta_{S, T}(K/F, -k)$ is contained in the product of the ideals $\cI_v^{(k)}$ taken over all archimedean places $v$ of $F$ except one (see \cite{dassp}, Thm.\ 5.9 (b)).
 \end{remarks}

We now explain the strategy of the proof of Theorem \ref{theorem:highervanishing}. Our method is of course related to that developed in \cite{dassp}. We use however a somewhat different homological algebra machine that is partly inspired by ideas of Hirose \cite{hirose}. Firstly, as in \cite{dassp} we represent the Stickelberger element at $s=-k$ as a cap product of some adelic Eisenstein class $\bEis$ with a homology class that is naturally associated to the global reciprocity map for the extension $K/F$. The key novelty 
in our approach is that we enrich this homology class ''locally'' at each place in $S\setminus \{\fp\}$. This new class does not lie in a homology group anymore but in a certain hyperhomology group that can be capped with $\bEis$ as well. The resulting element then lies in a degree-$0$ hyperhomolgy group that maps naturally into the product of the ideals $\cI_v^{(k)}$ thus implying the divisibility result \eqref{highvanishing2} (firstly though only when $T=\{\fq\}$ and up denominators that are powers of the residue characteristic of $\fq$; by varying $\fq$ and using an argument involving \v{C}ebotarev's density theorem we are then able to get rid of the denominators).

\paragraph{Notation} Throughout this paper we use the following notation. For sets $X$ and $Y$ we let $\Maps(X,Y)$ be the set of maps $X\to Y$. For a partially ordered set $\cI=(\cI, \le)$ we denote by $\cI^{\opp}$ the same set but with the reversed partial ordering.
The torsion subgroup of an abelian group $A$ will be denoted by $A_{\tor}$. 

By $\spaces$ we denote the category of Hausdorff spaces with continuous maps as morphisms. 
If $X, Y\in \spaces$ then we let $C(X,Y)\subseteq \Maps(X,Y)$ denote the subset of continuous maps $X\to Y$. Note that if $Y$ is discrete then $C(X, Y)$ consists of locally constant maps $X\to Y$. If $Y=R$ is a ring (equipped with the discrete topology) then we let $C_c(X,R)$ denote the subset of $C(X, R)$ of locally constant maps with compact support. 

Unless stated otherwise all rings are commutative with $1\ne 0$. For a ring $R$ we denote by $\Mod_R$ the category of left $R$-modules. 
If $A$ is an $R$-algebra and $N$ an $R$-module then we put $N_A= N\otimes_R A$. 
If $M_{\bu}$ is a bounded complex of $R$-modules then we denote by $\sH_n(M_{\bu})$ its $n$-th homology module.
Also if $M_{\bu}$ and $N_{\bu}$ are bounded complexes of $R$-modules then we denote by $M_{\bu}\otimes_R N_{\bu}$ the associated double complex and also -- by abuse of notation -- its total complex. If $G$ is a group and $\chi: G \to R^*$ a character (i.e.\ a homomorphism) then we denote by $R(\chi)$ the $R[G]$-module $R$ with $G$-action given by $\chi$. More generally if $M$ is an $R[G]$-module then $M(\chi)$ denotes the $R[G]$-module $M\otimes_R R(\chi)$.

For a ring $R$ we denote by $\Aff(R)$ the group $R^*\ltimes R$, i.e.\ $\Aff(R)$ is the subgroup of $\GL_2(R)$ of matrices of the form $\begin{pmatrix} a & b\\
0 & 1\end{pmatrix}$ with $a\in R^*$ and $b\in R$. More generally for an $R$-module $M$ we denote by $\Aff_R(M)$ the group 
$\GL(M)\ltimes M$. We often identify an element $\varphi= (\alpha, m)\in \Aff_R(M)$ with the map $\varphi: M\to M$ it induces, i.e.\ the map $M\to M, x \mapsto \alpha(x) + m$.

Let $V$ be a finite-dimensional $\bQ$-vector space. By a {\it lattice} $L$ in $V$ we mean a finitely generated subgroup of $V$ that generates $V$ as a vector space, i.e.\ we have $\rank(L) = \dim(V)$. The set of lattices in $V$ will be denoted by $\Lat(V)$. More generally if $F$ is an algebraic number field with ring of integers $\cO_F$ and if $V$ is a finite-dimensional $F$-vector space then we denote by $\Lat_{\cO_F}(V)$ the set of all lattice $L$ in $V$ that are also $\cO_F$-submodules of $V$. In particular $\Lat_{\cO_F}(F)$ is the set of fractional ideals of $F$. We say that $V$ is {\it oriented} if $V_{\bR}$ is equipped with an orientation. In this case a basis $(v_1, \ldots, v_n)$ of $V$ is called positively oriented if it belongs to the orientation. 

Places of $F$ will be denoted by $v,w$ or also by $\fp, \fq$ etc.\ if they are finite. In the latter case we denote the corresponding prime ideal of $\cO_F$ by $\fp, \fq$ etc.\ as well. The norm of a fractional ideal $\fa$ of $F$ will be denoted by $\Norm(\fa)$.
By $E_F = \cO_F^*$ we denote the group of global units of $F$. More generally for a finite set $S$ of nonarchimedean places of $F$ we let $E_S = E_{F,S}$ be the group of $S$-units of $F$. For a prime number $p$ we denote by $S_p$ the set of primes of $F$ that lie above $p$ and by $S_{\infty}$ the set of archimedean places of $F$. For a place $v$ of $F$ we denote by $F_v$ the completion of $F$ at $v$. Also we let $|\wcdot|_v$ be the associated normalized multiplicative valuation on $F_v$. If $v$ is nonarchimedean then $\cO_v$ denotes the valuation ring of $F_v$, $U_v = \cO_v^*$ its group of units and $\ord_v$ the corresponding the normalized (additive) valuation on $F_v$ (so $\ord_v(\varpi_v) =1$ if $\varpi_v\in \cO_v$ is a local uniformizer at $v$). Moreover if $v=\fp$ is finite then given a non-negative integer $m\ge 0$ we let $U_{\fp}^{(m)}$ be the $m$-th higher unit group, i.e.\ $U_{\fp}^{(m)}=\{x\in U_{\fp}\mid x \equiv 1 \mod \fp^m\cO_{\fp}\}$. 

In sections \ref{section:partzeta} and \ref{section:lvalues} we assume that $F$ is totally real field of degree $n$ over $\bQ$. Recall that  if $v\in S_{\infty}$ corresponds to the embedding $\xi: F \to \bR$ then $|x|_v = |\xi(x)|$ and if $v=\fq$ is finite then $|x|_{\fq} = \Norm(\fq)^{-\ord_{\fq}(x)}$. For $v\in S_{\infty}$ we put $U_v = \bR_+ = \{x\in \bR\mid x>0\}$. We denote by $\bA= \bA_F$ (resp.\ $\bA_f=\bA_{F, f}$) the ring of adeles of $F$ (the ring of finite adeles). For a set $S$ of places of $F$ we let $\bA\!^S$ (resp.\ $\bA_f^S$) denote the ring of $S$-adeles (resp.\ finite $S$-adeles). We also define $U^S =\prod_{v\not\in S} U_v$, $U_f^S =\prod_{v\not\in S, v\nmid \infty} U_v$ and $U_S = \prod_{v\in S} U_v$. If $S$ contains all archimedean places then the factor group $(\bA^S)^*/U^S$ is canonically isomorphic to the group $\cI^S$ of fractional $\cO_F$-ideals that are coprime to $S$. We sometimes view $F$ as a subring of $\bA\!^S$ via the diagonal embedding. If $S$ consists of finitely many nonarchimeden places then we let $\cO_S: = F\cap \bigcap_{v\in S}\cO_v$ be the associated semilocal subring of $F$. 

If $\Sigma$ is a set of places of $\bQ$ then $S_{\Sigma}$ denotes the set of places of $F$ which lie above a place of $\Sigma$. We often write $\bA_{S, \Sigma}$, $\bA\!^{S, \Sigma}$ etc.\ for $\bA_{S\cup S_{\Sigma}}$, $\bA\!^{S\cup S_{\Sigma}}$ etc. We also write $U^{S, p}$, $U_{S, p}$, $U^{S, \infty}$, $U^{p, \infty}$ etc.\ for $U^{S, \{p\}}$, $U_{S, \{p\}}$, $U^{S\cup S_{\infty}}$, $U^{S_p\cup S_{\infty}}$ etc. and use a similar notation for adeles. For example $\bA^{S, \infty}=\bA_f^S$ denotes the ring $S\cup S_{\infty}$-adeles of $F$. If $S=\emptyset$ then we drop it from the notation (e.g.\ $\bA^p$ denotes the set of $S_p$-adeles of $F$ for a prime $p$).   
For $p\in\{ 2,3,5, \ldots, \infty\}$ we put $F_p= F\otimes \bQ_p = \prod_{v\in S_p} F_v$. We shall denote by $F^*_+$, $E_+= E_{F,+}$, $E_{S,+}$ etc.\ the elements of $F$, $E_F$, $E_{S}$ etc.\ that are positive with resp.\ to every embedding $F\hra \bR$.

For an ideal $\fm\subseteq \cO_F$, $\fm\ne (0)$ we let $F^{\fm}$ be the ray class field of $F$ in the narrow sense associated with $\fm$. If $U_{\fm}$ denotes the open subgroup $U_{\fm} =\prod_{\fp\nmid \infty} U_{\fp}^{(m_{\fp})}\times \prod_{v\mid \infty} U_v$ of $\bA^*$ with $m_{\fp}$ being the exponent of $\fp$ occurring in $\fm$, then the Galois group $\Gal(F^{\fm}/F)$ is isomorphic to $\bA^*/F^*U_{\fm}$ via the reciprocity map. Or in terms of ray class groups we have $\Gal(F^{\fm}/F)\cong \cI^{\fm}/\cP^{\fm}$. Here $\cI^{\fm}:=\cI^{S_{\fm}}$ and $S_{\fm}$ is the set of nonarchimedan places of $F$ that divide $\fm$ and $\cP^{\fm}$ is the subgroup of $\cI^{\fm}$ consisting of principal fractional ideals $(x)$ with $x\in F^*_+$ and $x\in U_{\fp}^{(m_{\fp})}$ for all $\fp\in S$. We also denote by $E_{\fm}$ (resp.\ $E_{\fm, +}$) the subgroup of $E_F$ consisting of units $\ep\in E_F$ (resp.\ $\ep\in E_{F,+}$) that are $\equiv 1$ modulo $\fm$. Thus we have $E_{\fm, +}= F^*\cap U_{\fm}$.

\section{Locally polynomial functions and distributions}
\label{section:locpol}

\paragraph{Polynomial functions and polynomial distributions on lattices}

In the following by a {\it lattice} we mean a free-abelian group $L$ of finite rank. For a subring $R$ of $\bQ$ we recall the notion of an $R$-valued polynomial function on $L$ (see e.g.\ \cite{cahen}, Ch.\ XI).

\begin{df} 
\label{df:polfunct}
Let $L$ be a lattice of rank $n$. A map $f: L \to R$ is called an ($R$-valued) polynomial function if for some (hence any) choice of a $\bZ$-Basis $\la_1, \ldots, \la_n$ 
of $L$ there exists a polynomial $P(X_1, \ldots, X_n)\in \bQ[X_1, \ldots, X_n]$ such that 
\[
f(x_1\la_1 + \ldots + x_n \la_n) \, =\, P(x_1, \ldots, x_n)\qquad \forall \, (x_1, \ldots, x_n)\in \bZ^n.
\]
By $\Int(L, R)$ we denote the ring of all $R$-valued polynomial functions $f: L\to R$.
\end{df}

In particular if $L=\bZ^n$ then a polynomial $P(X_1, \ldots, X_n)\in \bQ[X_1, \ldots, X_n]$ -- viewed as a map $\bZ^n\ni (x_1, \ldots, x_n)\mapsto P(x_1, \ldots, x_n)$ -- lies in $\Int(\bZ^n, R)$ if and only if it has values in $R$. For $m\in \bZ_{\ge 0}$ we put $\binom{X}{m} =\frac{\prod_{i=0}^{m-1} (X-i)}{m!}\in \bQ[X]$.
More generally for a multi-index $\um=(m_1, \ldots, m_n)\in (\bZ_{\ge 0})^n$ we set
\begin{equation} 
\label{polyabase}
\binom{\uX}{\um} := \prod_{i=1}^n \binom{X_i}{m_i}\in \bQ[X_1, \ldots, X_n].
\end{equation}
Any polynomial $P\in \bQ[X_1, \ldots, X_n]$ can be uniquely written as
\begin{equation*}
\label{polya1}
P(X_1, \ldots, X_n)\, =\, \sum_{\um\in (\bZ_{\ge 0})^n} a(\um) \binom{\uX}{\um}
\end{equation*}
where the coefficients $a(\um)$ lie in $\bQ$ and vanish for almost all $\um\in (\bZ_{\ge 0})^n$. The polynomial $P$ lies in $\Int(\bZ^n, R)$ if and only if all its coefficients $a(\um)$ lie in $R$, i.e.\ the polynomial functions \eqref{polyabase} form an $R$-basis of $\Int(\bZ^n, R)$ for any subring $R\subseteq \bQ$ (cf.\ \cite{polya, ostrowski}). 

For an arbitrary ring $R$ we define
\begin{equation*}
\label{Aintpol}
\Int(L, R):= \,\Int(L, \bZ)\otimes R.
\end{equation*}
As explained above if $R\subseteq \bQ$ then this definition agrees with the previous one. We can attach to an element $f=\sum_{i=1}^r f_i\otimes a_i \in \Int(L, R)$ a map $L\to R$ given by $\la\ni L \mapsto \sum_{i=1}^r f_i(\la) a_i$. Therefore elements of $\Int(L, R)$ will be called $R$-valued polynomial functions on $L$ and we use the symbolic notation $f:L\to R$. Note though that for general $R$ the passage from an element of $\Int(L, R)$ to the map $L\to R$ ''looses information'' (i.e.\ the obvious map $\Int(L, R)\to \Maps(L, R)$ is in general not injective). 

If $\varphi: L'\to L$ is an affine map between lattices (i.e.\ $L'\to L, \la\mapsto \varphi(\la) -\varphi(0)$ is a homomorphism) then $\varphi^*:\Maps(L, \bZ) \to \Maps(L', \bZ), \, f\mapsto f\circ \varphi$ maps the subring $\Int(L, \bZ)$ into $\Int(L', \bZ)$ hence induces -- by extension of scalars -- a ring homomorphism $\varphi^*: \Int(L, R)\to \Int(L', R)$. 
In particular if $L'\subseteq L$ and $\varphi$ is the inclusion then we denote $\varphi^*$ by $\Int(L, R)\to \Int(L', R), f\mapsto f|_{L'}$ and call it {\it restriction}. If $L=L'$ and $\varphi:L\to L$ is the translation by a fixed element $\la\in L$ then we denote $\varphi^*$ by $\tau_{\la}$ (hence we have $\tau_{\la}(f)(v) = f(v+\la)$ for all $v\in L$). The map $\tau_{\la}: \Int(L, R)\to \Int(L, R)$ will be called translation (by $\la$). 

\begin{lemma}
\label{lemma:convol}
(a) Let $L'\subseteq L$ be a subgroup of the same rank $n$ and assume that the index $d=[L:L']$ is invertible in $R$. Then the restriction 
$\Int(L, R)\to \Int(L', R), f\mapsto f|_{L'}$ is an isomorphism. 
\medskip

\noi (b) For every polynomial function $f: L\to R$ there exists finitely many polynomial functions $f_1, \ldots, f_r, g_1, \ldots, g_r:L\to R$ such that 
\begin{equation}
\label{conv2}
\tau_{\la}(f)\, =\, \sum_{i=1}^r g_i(\la) f_i
\end{equation}
for every $\la\in L$.
\end{lemma}

\begin{proof} (a) We choose basis $\la_1, \ldots, \la_n$ of $L$ and positive integers $d_1, \ldots, d_n$ such that $d_1\la_1, \ldots, d_n\la_n$ is a basis of $L'$ so that $d=d_1 \cdot \ldots \cdot d_n$ . If $P(X_1, \ldots, X_n)\in \bQ[X_1, \ldots, X_n]$ is a polynomial such that $P(d_1 x_1, \ldots, d_n x_n)\in \bZ[1/d]$ for every $(x_1, \ldots, x_n)\in \bZ^n$ then, clearly, we also have $P(x_1, \ldots, x_n)\in \bZ[1/d]$ for every $(x_1, \ldots, x_n)\in \bZ^n$.
Thus the restriction $\Int(L, \bZ[1/d])\to \Int(L', \bZ[1/d]), f\mapsto f|_{L'}$ is an isomorphism. Tensoring it with $R$ yields the assertion.

(b) It suffices to consider the case $L= \bZ^n$, $R=\bZ$ and $f= \binom{\uX}{\um}$. In this case the assertion follows from 
\[
 \binom{\uX+ \uY}{\um} = \sum_{\uk\in (\bZ_{\ge 0})^n, \uk\preceq \um}  \binom{\uX}{\uk}  \binom{\uY}{\um-\uk}.
\] 
Here we have equipped $(\bZ_{\ge 0})^n$ with the following partial order:  for $\uk=(k_1, \ldots, k_n), \ul=(l_1, \ldots, l_n)\in (\bZ_{\ge 0})^n$ we define $\uk \preceq \ul$ if $k_i \le l_i$ for all $i=1, \ldots, n$.
\end{proof}

\begin{df} 
\label{df:dispolfunct}
Let $R$ be a ring. The $R$-module of ($R$-valued) polynomial distributions on $L$ is defined as 
\[
\cD_{\pol}(L, R)=\Hom_{\bZ}(\Int(L, \bZ), R) = \Hom_{R}(\Int(L, R), R).
\]
\end{df}

For $\mu\in \cD_{\pol}(L, R)$ we write $\int_L f(\la) d\mu(\la)$
for the evaluation of $\mu$ at $f\in \Int(L, R)$. Convolution defines a ring structure on $\cD_{\pol}(L, R)$. Concretely, the product $\mu_1\star \mu_2$ of two elements $\mu_1, \mu_2\in \cD_{\pol}(L, R)$ is defined as usual by
\begin{equation}
\label{conv}
\int_L f(\la) d(\mu_1\star \mu_2)(\la) = \int_L \left( \int_L  f(\la_1 + \la_2) d\mu_2(\la_2)\right) d\mu(\la_1)
\end{equation}
for every $f\in \Int(L, R)$. Note that by Lemma \ref{lemma:convol} (b) the map $\la_1\mapsto \int_L  f(\la_1 + \la_2) d\mu_2(\la_2)$ lies again in $\Int(L, R)$, so that \eqref{conv} is well-defined. More precisely, for $f\in \Int(L, R)$ we can choose $f_1, \ldots, f_r, g_1, \ldots, g_r\in \Int(L, R)$ such that \eqref{conv2} holds. Then  
\[
L\lra R, \quad \la_1\mapsto \int_L  f(\la_1 + \la_2) d\mu_2(\la_2)
\]
is the symbolic notation for the element $\sum_{i=1}^r \left(\int_L f_i(\la_2)d\mu_2(\la_2)\right) g_i$ of $\Int(L, R)$. 

\begin{example}
\label{example:distchar0}
\rm Let $L$ be a lattice of rank $n$ and let $R=K$ be a field of characteristic $0$. Let $\xi_1, \ldots, \xi_n: L \to K$ be a $K$-basis of $\Hom(L, K)$. Then $\xi_1, \ldots, \xi_n$ can be viewed as elements of $\Int(L, K)$. In fact elements of $\Int(L, K)$ can be written as polynomials in $\xi_1, \ldots, \xi_n$, i.e.\ the collection of functions $\uxi^{\um} := \prod_{k=1}^n \xi_k^{m_k}$ for $\um\in (\bZ_{\ge 0})^n$ form a $K$-basis of $\Int(L, K)$. Hence there exists unique elements $z_1, \ldots, z_n\in \cD_{\pol}(L, K)$ given by $z_i(\uxi^{\um}) = 1$ if $\um = e_i$ and $z_i(\uxi^{\um}) = 0$ if $\um\ne e_i$. For $\ua=(a_1, \ldots, a_n)\in (\bZ_{\ge 0})^n$ the element $\uz^{\ua}:= z_1^{a_1} \cdot \ldots \cdot z_n^{a_n}\in \cD_{\pol}(L, K)$ is characterised by 
\begin{equation}
\label{normdist}
\uz^{\ua}(\uxi^{\um}) \, =\, \left\{ \begin{array}{cc} \ua! & \mbox{if $\ua = \um$,}\\
                                                      0 & \mbox{otherwise}
                                                      \end{array}\right.
\end{equation}
(with $\ua!=\prod_{k=1}^n a_k!$) for every $\um\in (\bZ_{\ge 0})^n$.
\end{example}

Let $R[L]$ denote the $R$-group algebra of $L$. For $\la_0 \in L$ we let $\delta_{\la_0}\in \cD_{\pol}(L, R)$ be the Dirac distribution, i.e.\ for $f\in \Int(L, R)$ we have $\int_L f(\la) d\delta_{\la_0}(\la) =f(\la_0)$. Since $\delta_{\la_1}\star \delta_{\la_2}= \delta_{\la_1 + \la_2}$ for all $\la_1, \la_2\in L$ the map $L\to \cD_{\pol}(L, R), \la\mapsto \delta_{\la}$ extends to a ring homomorphism 
\begin{equation}
\label{iwasawa}
\delta: R[L] \lra \cD_{\pol}(L, R),\quad \sum_{\la\in L} a_{\la} [\la] \mapsto \sum_{\la\in L} a_{\la} \ \delta_{\la}.
\end{equation} 
The collection of translations $\tau_{\la}$ for $\la\in L$ induce an $R[L]$-module structure on $\Int(L, R)$ given by 
\begin{equation*}
\label{almod}
\star: R[L] \times \Int(L, R) \lra \Int(L, R), \quad  \left(\sum_{\la\in L} a_{\la} [\la]\right)\star f:= \sum_{\la\in L} a_{\la} \cdot\tau_{\la}(f).
\end{equation*}
Note that for $\mu\in \cD_{\pol}(L, R)$, $\alpha\in R[L]$ and $f\in \Int(L, R)$ we have 
\begin{equation}
\label{almod2}
\int_L f(\la) \, d(\delta(\alpha) \star \mu)(\la) \, = \, \int_L (\alpha\star f)(\la)\, d\mu(\la).
\end{equation}
By $R[\![L]\!]$ we denote the completion of $R[L]$ with respect to the kernel $I(L)$ of the augmentation map $\vep: R[L] \to R, \sum_{\la\in L} a_{\la} [\la]\mapsto \sum_{\la\in L} a_{\la}$, i.e.\ we have
\begin{equation*}
\label{iwasawa1}
R[\![L]\!] \, =\, \prolimm R[L]/I(L)^m.
\end{equation*} 

\begin{prop}
\label{prop:poldispower}
The homomorphism \eqref{iwasawa} induces an isomorphism of $R$-algebras
\begin{equation*}
\label{iwasawa2}
R[\![L]\!]\lra \cD_{\pol}(L, R).
\end{equation*} 
\end{prop}

\begin{proof} For $m\in \bZ_{\ge 0}$ we let $\Int^m(L, \bZ)$ be the submodule of $\Int(L, \bZ)$ of integer-valued polynomial functions $f:L\to \bZ$ of degree $\le m$, i.e.\ the total degree of the polynomial $P(X_1, \ldots, X_n)$ in Definition \ref{df:polfunct} is $\le m$. Put $\Int^m(L, R):=\Int^m(L, \bZ)\otimes R$, $\cD_{\pol, m}(L, R):=\Hom(\Int^m(L, \bZ), R)$ and let
\begin{equation}
\label{iwasawa2a}
R[L] \lra  \cD_{\pol, m}(L, R), \quad \sum_{\la\in L} a_{\la} [\la] \mapsto \sum_{\la\in L} a_{\la} \ \res(\delta_{\la})
\end{equation} 
be the composite of \eqref{iwasawa} with the obvious restriction map $\res:\cD_{\pol}(L, R)\to \cD_{\pol, m}(L, R)$. We will show that \eqref{iwasawa2a} is surjective with kernel $I(L)^{m+1}$. Hence it induces an isomorphism 
\begin{equation}
\label{iwasawa3}
R[L]/I(L)^{m+1} \lra  \cD_{\pol, m}(L, R).
\end{equation} 
Passing to the inverse limit over all $m$ yields the assertion. 

For the surjectivity of \eqref{iwasawa2a} it suffices to consider the case $L=\bZ^n$. Let $\mu \in \cD_{\pol, m}(\bZ^n, R)$ and let $\Xi$ be the set of 
$\uk=(k_1, \ldots, k_n)\in (\bZ_{\ge 0})^n$ such that $\sum_{i=1}^n k_i \le m$. To see that $\mu$ lies in the image of \eqref{iwasawa2a} it is enough to show that there exists $\left(a_{\uk}\right)_{\uk\in \Xi}\in R^{\Xi}$ such that
\begin{equation}
\label{iwasawa4}
\mu\, =\, \sum_{\uk\in \Xi} a_{\uk} \res(\delta_{\uk}).
\end{equation} 
Since the polynomials $\binom{\uX}{\ul}$, with $\ul\in \Xi$ form an $R$-basis of $\Int^m(\bZ^n, R)$, an element $\left(a_{\uk}\right)_{\uk\in \Xi}\in R^{\Xi}$ satisfies \eqref{iwasawa4} if and only if it is a solution to the system of linear equations
\begin{equation}
\label{iwasawa5}
\mu\left(\binom{\uX}{\ul}\right)\, =\, \sum_{\uk\in \Xi} a_{\uk} \binom{\uk}{\ul}\qquad \forall\, \ul\in \Xi.
\end{equation} 
For that we equip $\Xi$ with the lexicographic order, i.e.\ for $\uk=(k_1, \ldots, k_n), \ul=(l_1, \ldots, l_n)$ we have $\uk < \ul$ if there exists an index $i\in \{1, \ldots, n\}$ such that $k_j=l_j$ for $j=1, \ldots, i-1$ and $k_i < l_i$. Note that
\[
\binom{\uk}{\ul}\, =\, \left\{ \begin{array}{cc} 0 & \mbox{if $\uk < \ul$,}\\
                                                      1 & \mbox{if $\uk = \ul$.}
                                                      \end{array}\right.
\]
so that $\left( \binom{\uk}{\ul}\right)_{\uk, \ul\in \Xi}$ is an upper triangular matrix with entries in $\bZ$ and diagonal entries $=1$. Hence the system of equations \eqref{iwasawa5} has a unique solution $\left(a_{\uk}\right)_{\uk\in \Xi}\in R^{\Xi}$. 

To show that $I(L)^{m+1}$ lies in the kernel of \eqref{iwasawa2a} we argue as in (\cite{dassp}, Lemma 4.3).
Note that for $\la\in L$ and $f\in \Int^m(L, R)$ we have $([\la] - [0])\star f= \tau_{\lambda}(f) - f\in \Int^{m-1}(L, R)$. It follows that $I(L)\star \Int^m(L, R)\subseteq \Int^{m-1}(L, R)$, hence $I(L)^{m+1}\star \Int^m(L, R)=0$. Since by \eqref{almod2} we have for $\alpha\in I(L)^{m+1}$ and $f\in \Int^m(L, R)$ 
\begin{equation*}
\label{almod3}
\int_L f(\la) \, d \delta(\alpha) (\la) \, = \, \int_L (\alpha\star f)(\la) \, d\delta_0(\la) \, =\, 0
\end{equation*}
we conclude that $I(L)^{m+1}$ lies in the kernel of \eqref{iwasawa2a}. 

We have shown that \eqref{iwasawa2a} induces an epimorphism \eqref{iwasawa3}. To show that the latter is an isomorphism it suffices to remark that both source and target are free $R$-modules of the same rank. Since 
$\Int^m(L, R)$ is a free $R$-module of rank $\sharp(\Xi) = \binom{m+n}{m}$, the same is true for its dual $\cD_{\pol, m}(L, R)$. On the other hand the choice of a basis $\la_1, \ldots, \la_n$ of $L$ yields an isomorphism $R[T_1^{\pm 1}, \ldots, T_n^{\pm 1}]\cong R[L]$
hence an isomorphism 
\begin{equation*}
\label{iwasawa9}
R[t_1, \ldots, t_n]/(t_1, \ldots, t_n)^{m+1} \, \cong \, R[L]/I(L)^{m+1}
\end{equation*}
with $t_i := T_i-1$ for $i=1, \ldots, n$. Thus $R[L]/I(L)^{m+1}$ is a free $R$-module of rank $\binom{m+n}{m}$ as well. 
\end{proof}

The lattice $L$ together with the choice of an orientation (i.e.\ an orientation on $L_{\bR}$) will be called an {\it oriented lattice}. A $\bZ$-basis $(\la_1, \ldots, \la_n)$ of $L$ will be called positively oriented if the induced isomorphism $L_{\bR}\to \bR^n$ preserves the orientation

\begin{coro}
\label{coro:koszul}
Let $L$ be an oriented lattice of rank $n$. Then there exists a canonical isomorphism  
\[
\Ext_{R[L]}^i(R, \cD_{\pol}(L, R))\, =\, \left\{ \begin{array}{cc} R & \mbox{if $i=n$,}\\
                                                      0 & \mbox{if $i \ne n$.}
                                                      \end{array}\right.
\]
for every $i\in \bZ_{\ge 0}$.
\end{coro}

\begin{proof} This follows from Prop.\ \ref{prop:poldispower} and (\cite{bkl}, Thm.\ 3.25). We give a purely algebraic proof. First we assume that $L=\bZ^n$ so that $R[\bZ^n]$ can be identified with the ring of Laurent polynomials in $n$ variable $R[T_1^{\pm 1}, \ldots, T_n^{\pm 1}]$, $I(L)$ with the ideal 
$(t_1, \ldots, t_n)$ with $t_i:=T_i-1$ and  
$\cD_{\pol}(\bZ^n, R)\cong R[\![\bZ^n]\!]$ with the power series ring $R[\![t_1, \ldots, t_n]\!]\cong \prolim_m R[T_1^{\pm 1}, \ldots, T_n^{\pm 1}]/(t_1, \ldots, t_n)^m$. We have to show 
\begin{equation}
\label{koszul1}
\Ext_{R[T_1^{\pm 1}, \ldots, T_n^{\pm 1}]}^i(R, R[\![t_1, \ldots, t_n]\!])\, =\, \left\{ \begin{array}{cc} R & \mbox{if $i=n$,}\\
                                                      0 & \mbox{if $i \ne n$.}
                                                      \end{array}\right.
\end{equation}
Since $\ut:=(t_1, \ldots, t_n)$ form a regular sequence of $R[T_1^{\pm 1}, \ldots, T_1^{\pm 1}]$ we have (see \cite{weibel}, Cor.\ 4.5.5 and Ex.\ 4.5.2) 
\[
\Ext_{R[T_1^{\pm 1}, \ldots, T_n^{\pm 1}]}^i(R, M)\, =\, H^i(\ut, M)\, =\, H_{n-i}(\ut, M)
\]
for any $R[T_1^{\pm 1}, \ldots, T_n^{\pm 1}]$-module $M$ and any $i\in \bZ_{\ge 0}$ (here $H^{\bu}(\ut, M)$ and $H_{\bu}(\ut, M)$ denote Koszul cohomology and homology). Note also that $t_1, \ldots, t_n$ forms a regular sequence of $R[\![t_1, \ldots, t_n]\!]$ and that we have 
\[
R[\![t_1, \ldots, t_n]\!]/(t_1, \ldots, t_n) \cong R[\![t_1, \ldots, t_{n-1}]\!]/(t_1, \ldots, t_{n-1})\cong \ldots \cong R[\![t_1]\!]/(t_1)\cong R.
\]
Thus \eqref{koszul1} follows from (\cite{weibel}, Cor.\ 4.5.4).

Now assume that $L$ is an arbitrary rank $n$ lattice. The choice of a basis $\ula=(\la_1, \ldots, \la_n)$ of $L$ induces an isomorphism 
\begin{equation}
\label{orientedtrace}
\Ext_{R[L]}^n(R, \cD_{\pol}(L, R))\, \cong \, \Ext_{R[\bZ^n]}^n(R, \cD_{\pol}(\bZ^n, R))\, \cong \, R.
\end{equation}
Changing the basis changes the isomorphism by the factor $\sgn(\det(A))$ where $A\in \GL_n(\bZ)$ is the associated transition matrix. Thus if we choose a basis that is positively oriented then the isomorphism \eqref{orientedtrace} does not depend on this choice.
\end{proof}

We finish this section by introducing generalizations of Definitions \ref{df:polfunct} and \ref{df:dispolfunct} to (left) $L$-sets. Recall that the latter is a set $H$ together with an $L$-action $+: L\times H\to H, (\lambda, h) \mapsto \la + h$. The example the reader should have in mind is that $L$ is a subgroup of a rational vector space $V$ and $H$ is an $L$-stable subset of $V$, i.e.\ we have $\la + h\in H$ for every $\la\in L$ and $h\in H$. An $L$-subset of $H$ is a subset $H'\subseteq H$ that is stable under the $L$-action, i.e.\ we have $\la +h\in H'$ for every $h\in H'$ and $\la \in L$. An $L$-set $H$ will be called {\it finite} if there are only finitely many $L$-orbits. Thus a finite $L$-subset $H'$ of an $L$-set $H$ is an $L$-subset that is finite as an $L$-set. 

\begin{df} 
\label{df:phpolfunct} Let $H$ be an $L$-set. 
\medskip

\noi (a) A map $f: H \to \bZ$ is called an integer-valued polynomial function if for every of $h\in H$ and every $\bZ$-Basis $\la_1, \ldots, \la_n$ 
of $L$ there exists a polynomial $P(X_1, \ldots, X_n)\in \bQ[X_1, \ldots, X_n]$ such that 
\[
f(x_1\la_1 + \ldots + x_n \la_n + h) \, =\, P(x_1, \ldots, x_n)\qquad \forall \, (x_1, \ldots, x_n)\in \bZ^n.
\]
The ring of integer-valued polynomial functions $f: H\to \bZ$ will be denoted by $\Int_L(H, \bZ)$ (or by $\Int(H, \bZ)$ for short). More generally for an arbitrary ring $R$ we define $\Int_L(H, R):= \Int_L(H, \bZ)\otimes R$. An element $f\in \Int_L(H, R)$ will be called an ($R$-valued) polynomial function on $H$ and we use again the symbolic notation $f:H\to R$. 
\medskip

\noi (b) An integer-valued polynomial function $f: H \to \bZ$ is said to be of bounded support if its support is contained in a finite $L$-subset $H'\subseteq H$, i.e.\ if there exists finitely many elements $h_1, \ldots, h_m\in H$ such that $\supp(f)=\{h\in H\mid f(h)\ne 0\}\subseteq \bigcup_{i=1}^m L+ h_i$. The ring of integer-valued locally polynomial functions $H\to \bZ$ with bounded support will be denoted by $\Int_{L, b}(H, \bZ)$ (or simply by $\Int_b(H, \bZ)$). For a ring $R$ we put $\Int_{L, b}(H, R):= \Int_{L, b}(H, \bZ)\otimes R$.
\medskip

\noi (c) The dual of $\Int_{L, b}(H, R)$ will be denoted by 
\[
\cD_{\pol}(H, R)\, =\, \cD_{\pol, L}(H, R) \, = \, \Hom_{\bZ}(\Int_{L, b}(H, \bZ), R) \, =\, \Hom_R(\Int_{L, b}(H, R), R).
\] 
It is called the module of $R$-valued polynomial distributions on $H$.
\end{df}

Clearly for a sublattice $L'$ of $L$ we have $\Int_L(H, \bZ)\subseteq \Int_{L'}(H, \bZ)$ and $\Int_{L, b}(H, \bZ)\subseteq \Int_{L', b}(H, \bZ)$. Also for an $L$-subset $H'\subseteq H$ there exists canonical ring homomorphisms
\begin{eqnarray}
\label{polres}
\Int_L(H, \bZ) \lra \Int_L(H', \bZ), & f \mapsto f|_{H'} & \mbox{"Restriction",}\\
\Int_L(H', \bZ) \lra \Int_L(H, \bZ), & f \mapsto f_! & \mbox{"Extension by zero"}
\label{polext0}
\end{eqnarray} 
that map $\Int_{L, b}(H, \bZ)$ (resp.\ $\Int_{L, b}(H', \bZ)$) into $\Int_{L, b}(H', \bZ)$ (resp.\ $\Int_{L, b}(H, \bZ)$). 

\begin{remarks}
\label{remarks:phpolfunc}
\rm (a) Let $\varphi: H_1 \to H_2$ be an $L$-equivariant map between $L$-sets. The map $\varphi^*: \Maps(H_2, \bZ) \to \Maps(H_1, \bZ), f\mapsto f\circ \varphi$ maps $\Int_L(H_2, \bZ)$ into $\Int_L(H_1, \bZ)$. Moreover if $\varphi$ is injective then $\varphi^*$ maps $\Int_{L,b}(H_2, \bZ)$ into $\Int_{L,b}(H_1, \bZ)$, i.e.\ it induces a homomorphism of $L$-algebras $\varphi^*: \Int_b(H_2, R)\to \Int_b(H_1, R)$. We denote the dual map by 
\begin{equation*}
\label{cofunc}
\varphi_*: \cD_{\pol, L}(H_1, R) \lra \cD_{\pol, L}(H_2, R).
\end{equation*} 

\noi (b) In particular for an $L$-set $H$  and $\la \in L$ the map $\varphi: H \to H, h \mapsto \la + h$ is $L$-equivariant. The induced homomorphism of $R$-algebras $\tau_{\la}:= \varphi^*: \Int(H, R)\to \Int(H, R)$ will be called translation by $\la$. Similar to Lemma \ref{lemma:convol} (b) for every polynomial function $f: H\to R$ there exists finitely many polynomial functions $f_1, \ldots, f_N: H\to R$ and $g_1, \ldots, g_N:L\to R$ such that 
\begin{equation}
\label{conv2a}
\tau_{\la}(f)\, =\, \sum_{i=1}^r g_i(\la) f_i
\end{equation}
for every $\la\in L$. Moreover if $f$ has  bounded support then $f_1, \ldots, f_N$ can be chosen to have bounded support as well. The collection of translations $\tau_{\la}$ for $\la\in L$ induce an $R[L]$-module structure on $\Int(H, R)$. It will be denoted by $\star: R[L] \times \Int(H, R) \to \Int(H, R)$.
\medskip

\noi (c) By using property \eqref{conv2a} for $\mu\in \cD_{\pol}(L, R)$ and $\nu\in \cD_{\pol, L}(H, R)$ one can define the convolution $\mu\star \nu\in \cD_{\pol, L}(H, R)$ similar to \eqref{conv} namely we have
\begin{equation*}
\label{conv3}
\int_H f(h) d(\mu\star \nu)(h) = \int_H \left( \int_L  f(\la + h) d\mu(\la)\right) d\nu(h)
\end{equation*}
for every $f\in \Int_{L,b}(H, R)$. Thus the convolution product defines on $\cD_{\pol, L}(H, R)$ a $\cD_{\pol}(L, R)$-module structure. 
\medskip

\noi (d) Let $H$ be an $L$-set and let $H= \bigcup_{i\in I} H_i$ be a covering by disjoint $L$-subsets. The family of maps $\Int_{\loc, b}(H_i, \bZ) \to \Int_{\loc, b}(H, \bZ), f_i \mapsto (f_i)_!$ induces an isomorphism
\begin{equation}
\label{intcosheaf}
\bigoplus_{i\in I} \Int_{L, b}(H_i, \bZ) \lra \Int_{L, b}(H, \bZ). 
\end{equation} 
For that it suffices to verify surjectivity. Let $f\in  \Int_{L, b}(H, \bZ)$ there exists a finite $L$-subset $H'$ of $H$ such that $\supp(f)\subseteq H'$. Note that the set $J:=\{i\in I\mid H_i\cap H'\ne \emptyset\}$ is finite since $H'$ has only finitely many $L$-orbits. If we put $f_i := f|_{H_i}$ for $i\in I$ then we have $f_i=0$ except possibly for $i\in J$. It follows that $(f_i)_{i\in I}\in \bigoplus_{i\in I} \Int_{L, b}(H_i, \bZ)$ is mapped to $f$ under \eqref{intcosheaf}.
\enddemo
\end{remarks}

\paragraph{Locally polynomial functions and distributions} In the following $V$ denotes a $\bQ$-vector space of dimension $n$ (in Prop.\ \ref{prop:koszul2} we assume additionally that $V$ is oriented). Recall that $\Lat(V)$ is the set of all subgroups $L\subseteq V$ that are free-abelian of rank $n$. Elements of $\Lat(V)$ will be called lattices. We fix a non-empty subset $\sL$ of $\Lat(V)$ such that we have 
\begin{equation}
\label{sLclosed}
 L_1\cap L_2, \quad L_1 + L_2 \in \sL \qquad \text{for all} \quad L_1, L_2\in \sL.
\end{equation}
For such $\sL$ we put
\begin{equation*}
\label{sLgen}
\Lambda\,=\, \Lambda(\sL):= \, \bigcup_{L\in \sL} L.
\end{equation*}
Note that $\Lambda$ is a subgroup of $V$. For such data we consider the group
\begin{equation*}
\label{autsl}
\Aut_{\sL}(\Lambda)\, =\, \{\alpha\in \GL(V)\mid \alpha(L), \alpha^{-1}(L)\in \sL \, \forall \, L\in \sL\}.
\end{equation*}
The elements of $\Aut_{\sL}(\Lambda)$ map $\Lambda$ onto itself, i.e.\ we have $\Aut_{\sL}(\Lambda)\subseteq \GL(\Lambda)$. Let $\Ga$ be a subgroup of $\GL(V)$. The set $\sL$ will be called $\Ga$-stable if $\Ga\subseteq \Aut_{\sL}(\Lambda)$.

\begin{examples}
\label{examples:fracideals} \rm (a) Let $F$ be a number field of degree $n$ over $\bQ$ with ring of integers $\cO_F$. The set of all fractional ideals $\cI=\Lat_{\cO_F}(F)$ of $F$ is a closed $F^*$-stable subset of $\Lat(F)$ and we have $\Lambda= F$. More generally, the set of fractional ideals $\cI^S$ of $F$ that are coprime to a fixed finite set $S$ of nonarchimedean places of $F$ is a closed subset of $\Lat_{\cO_F}(F)$. If $\cO_S=\{x\in F\mid \ord_{\fp}(x) \ge 0\,\forall \,\fp \in S\}$ is the semilocal subring of $F$ associated to $S$ then we have $\Lambda(\cI^S) = \cO_S$ and $\cI^S$ is $\cO_S^*$-stable.
\medskip

\noi (b) More generally let $n=dm$, let $F$ be a number field of degree $d$ over $\bQ$ and let $V$ be a $m$-dimensional $F$-vector space. Let $S$ as above be a finite set of nonarchimedean places of $F$ and let $\cM$ be a finitely generated $\cO_S$-submodule of $V$ with $\rank_{\cO_S} \cM = m$, i.e.\ we have $\cM\otimes_{\cO_S} F=V$. Then we consider the set $\sL=\sL(\cM)$ of finitely generated $\cO_F$-submodules $L\subseteq \cM$ satisfying $L\otimes_{\cO_F} \cO_S = \cM$. It is a closed $\GL_{\cO_S}(\cM)$-stable subset of $\Lat_{\cO_F}(V)$ and we have $\Lambda(\sL) = \cM$.

For later use we note that the set $\sL(\cM)$ admits an adelic description (similar to the description of $\cI^S$ as $\cI^S\cong (\bA_f^S)^*/U_f^S$). More precisely there is a canonical transitive action of the group $\GL_{\bA_f^S}(\bA_f^S\otimes_{\cO_S} V)$ on $\sL(\cM)$: given $L\in \sL(\cM)$ and $g\in \GL_{\bA_f^S}(\bA_f^S\otimes_{\cO_S} V)$ we define $g\cdot L$ by $g\cdot L = \, g (L\otimes_{\cO_F}\widehat{\cO}^S) \cap \cM$
where $\widehat{\cO}^S:=\prod_{v\not\in S, v\nmid \infty} \cO_v$ (the intersection is taken within $V\otimes_F \bA_f^S$). Note that the stabilizer of $L\in \sL(\cM)$ is the group $\GL_{\widehat{\cO}^S}(\widehat{L}^S)$ where 
$\widehat{L}^S = L\otimes_{\cO_F}\widehat{\cO}^S$.
\enddemo
\end{examples}

Let $\cH$ be a (left) $\Lambda$-set, i.e.\ $\cH$ is a set together with an $\Lambda$-action $\Lambda\times \cH \to \cH, (\la, h) \mapsto \la +h$. 
For $L\in \sL$ we say that a subset $H\subseteq \cH$ is $L$-stable if $\la +x \in H$ for every $x\in H$ and $\la \in L$. We put
\begin{equation*}
\label{lstablesets}
\fOpen(\cH) \,=\, \fOpen_{\sL}(\cH)\, :=\, \{H\subseteq \cH\mid \mbox{$H$ is $L$-stable for some $L\in \sL$}\}.
\end{equation*}
An element $H\in \fOpen(\cH)$ will be called {\it finite} if $H$ is $L$-stable and has only finitely many $L$-orbits for some $L\in \sL$. The collection of finite elements of $\fOpen(\cH)$ will be denoted by $\fOpen_{\fin}(\cH)$. Also for $H\in \fOpen(\cH)$ we let $\fOpen_{\fin}(H)$ be the collection 
of subsets $H'\subseteq H$ with $H'\in \fOpen_{\fin}(\cH)$. 

We are now in the position of defining $R$-valued locally polynomial functions and locally polynomial distributions with respect to the family of lattices $\sL$. 

\begin{df} 
\label{df:locpolfunct}
Let $\cH$ be a $\Lambda$-set and let $H\in \fOpen(\cH)$.
\medskip

\noi (a) A map $f: H \to \bZ$ will be called an integer-valued locally polynomial function (with respect to $\sL$) if there exists $L\in \sL$ such that $H$ is $L$-stable and so that $f\in \Int_L(H, \bZ)$. The ring of all integer-valued locally polynomial functions $H\to \bZ$ with respect to $\sL$ will be denoted by $\Int_{\sL\mloc}(H, \bZ)$ or simply by $\Int_{\loc}(H, \bZ)$. For a ring $R$ we define $\Int_{\loc}(H, R):= \Int_{\loc}(H, \bZ)\otimes R$. Elements of $\Int_{\loc}(H, R)$ will be called $R$-valued locally polynomial functions on $H$.
\medskip

\noi (b) An element $f\in \Int_{\loc}(H, \bZ)$ is called an integer-valued polynomial function with bounded support if there exists a subset $H'\subseteq H$ with $H'\in \fOpen_{\fin}(\cH)$ and $\supp(f)\subseteq H'$. The abelian group of locally polynomial functions $H\to \bZ$ with bounded support will be denoted by $\Int_{\sL\mloc, b}(H, \bZ)=\Int_{\loc, b}(H, \bZ)$. For a ring $R$ we define $\Int_{\loc, b}(H, R):= \Int_{\loc, b}(H, \bZ)\otimes R$. 
\medskip

\noi (c) The dual of $\Int_{\loc, b}(H, R)$ will be denoted by 
\[
\cD_{\lpol}(H, R)\, =\,\cD_{\sL\mlpol}(H, R)\, =\, \Hom_{\bZ}(\Int_{\loc}(H, \bZ), R) = \Hom_R(\Int_{\loc}(H, R), R).
\]
It is called the module of locally polynomial $R$-valued distributions on $H$ with respect to $\sL$.
\medskip

\noi (d) An additive map $\mu: \Int_{\loc}(H, \bZ)\to R$ will be called an $R$-valued locally polynomial distribution on $H$ with bounded support, if there exists a subset $H'\subseteq H$ with $H'\in \fOpen_{\fin}(\cH)$ such that $\mu(f)$ depends only on $f|_{H'}$, i.e.\ there exists an additive map $\mu': \Int_{\loc}(H', \bZ)\to R$ such that $\mu(f) =\mu'(f|_{H'})$ for every $f\in \Int_{\loc}(H, \bZ)$. We denote by $\cD_{\lpol, b}(H, R)$ the module of $R$-valued locally polynomial distributions on $H$ with bounded support.
\end{df}

Note that $\Int_{\loc, b}(H, R)$ is an ideal in $\Int_{\loc}(H, R)$ (that $\Int_{\loc, b}(H, R)\subseteq \Int_{\loc}(H, R)$ follows from the fact that $\Int_{\loc}(H, \bZ)/\Int_{\loc, b}(H, \bZ)$ is torsionfree). Let $H\in \fOpen(\cH)$ and assume that $H$ is $L$-stable for $L\in \sL$. If we put $\sL(L):= \{L'\in \sL\mid L'\subseteq L\}$ and consider it as a partially ordered set (ordered by inclusions) then 
we have
\[
\Int_{\loc, b}(H, \bZ) = \dlim_{L'\in \sL(L)} \Int_{L', b}(H, \bZ) \qquad \mbox{and} \qquad \Int_{\loc}(H, \bZ) = \dlim_{L'\in \sL(L)} \Int_{L'}(H, \bZ).
\]
It follows that for a ring $R$ we have
\begin{equation}
\label{locpoldist}
\cD_{\lpol}(H, R) \, =\, \prolim_{L'\in \sL(L)^{\opp}} \cD_{\pol, L'}(H, R).
\end{equation}
Note that for an arbitrary coefficient ring $R$ an element $\Int_{\loc}(H, R)$ defines again a function $H\to R$.  

For $\mu\in \cD_{\lpol}(H, R)$ (resp.\ $\mu\in \cD_{\lpol, b}(H, R)$) we write $\int_H f(h) d\mu(h)$ for the evaluation of $\mu$ at $f\in \Int_{\loc, b}(H, R)$ (resp.\ at $f\in \Int_{\loc}(H, R)$). More generally, for an $R$-Algebra $A$ we have an evaluation pairing that will be denoted again by
\begin{equation}
\label{eval}
\cD_{\lpol}(H, R)\times \Int_{\loc, b}(H, A)\, \lra\, A, \quad (\mu, f) \mapsto \int_H f(h)d\mu(h).
\end{equation}
Let $H', H\in \fOpen(\cH)$ with $H'\subseteq H$. Again there exists canonical ring homomorphisms
\begin{eqnarray}
\label{locpolres}
\Int_{\loc}(H, \bZ) \lra \Int_{\loc}(H', \bZ), & f \mapsto f|_{H'} & \mbox{"Restriction",}\\
\Int_{\loc}(H', \bZ) \lra \Int_{\loc}(H, \bZ), & f \mapsto f_! & \mbox{"Extension by zero"}
\label{polext1}
\end{eqnarray} 
that map $\Int_{\loc, b}(H, \bZ)$ (resp.\ $\Int_{\loc, b}(H', \bZ)$) into $\Int_{\loc, b}(H', \bZ)$ (resp.\ $\Int_{\loc, b}(H, \bZ)$). 
Dually, we obtain maps 
\begin{eqnarray}
\label{distext0}
\cD_{\lpol}(H', R) \lra \cD_{\lpol}(H, R), & \mu \mapsto \mu_! & \mbox{"Extension by zero"}\\
\cD_{\lpol}(H, R) \lra \cD_{\lpol}(H', R), & \mu \mapsto \mu|_{L'} & \mbox{"Restriction"} 
\label{distres}
\end{eqnarray} 
i.e.\ \eqref{distext0} and \eqref{distres} are characterized by 
\begin{equation*}
\label{resexta}
\int_H f(h) \, d\mu_! (h) = \int_{H'} (f|_{H'})(h')\, d\mu(h'), \qquad \int_{H'} g(h')\, d(\nu|_{H'}) (\la')=\int_H g_{!}(h) \, d\nu (h)
\end{equation*}
for all $\mu\in \cD_{\lpol, L}(H', R)$, $f\in \Int_{\loc, b}(H, R)$ and $\nu\in \cD_{\lpol, L}(H, R)$, $g\in \Int_{\loc, b}(H', R)$. Moreover 
note that \eqref{locpolres} also induces a map 
\begin{equation}
\label{bdistext0}
\cD_{\lpol, b}(H', R) \lra \cD_{\lpol, b}(H, R), \quad \mu \mapsto \mu_! \qquad \mbox{"Extension by zero"}
\end{equation} 
Note that for $H\in \fOpen(\cH)$ we have
\begin{equation*}
\label{locintb}
\Int_{\loc, b}(H, \bZ) \, =\, \dlim_{H'\in \fOpen_{\fin}(H)} \Int_{\loc}(H', \bZ) 
\end{equation*}
where the transition maps in the limit are given by \eqref{polext1}. It follows  
\begin{equation}
\label{locpoldist2}
\cD_{\lpol}(H, R) \, =\, \prolim_{H'\in \fOpen_{\fin}(H)^{\opp}} \cD_{\lpol}(H', R).
\end{equation}
Moreover, we note that for the $R$-module $\cD_{\lpol, b}(H, R)$ we have 
\begin{equation*}
\label{locpoldistb}
\cD_{\lpol, b}(H, R) \, =\, \dlim_{H'\in \fOpen_{\fin}(H)} \cD_{\lpol}(H', R).
\end{equation*}
Here the transition map $\cD_{\lpol}(H_1', R)\to \cD_{\lpol}(H_2', R)$ for $H_1' \subseteq H_2'\subseteq H$ with  $H_1', H_2'\in \fOpen_{\fin}(\cH)$ is the map \eqref{distext0}.

Let $H\in \fOpen(\cH)$ and let $H= \bigcup_{i\in I} H_i$ be a covering of $H$ by disjoint subsets of $\fOpen(\cH)$. The family of maps $\Int_{\loc, b}(H_i, \bZ) \to \Int_{\loc, b}(H, \bZ), f_i \mapsto (f_i)_!$ induces a homomorphism 
\begin{equation}
\label{intloccosheaf}
\bigoplus_{i\in I} \Int_{\loc, b}(H_i, \bZ) \lra \Int_{\loc, b}(H, \bZ)
\end{equation} 
and by passing to duals a homomorphism 
\begin{equation}
\label{cdsheaf}
\cD_{\lpol}(H, R) \lra \prod_{i\in I} \cD_{\lpol}(H_i, R), \quad \mu \mapsto (\mu|_{H_i})_{i\in I} 
\end{equation}
for any ring $R$.

\begin{lemma}
\label{lemma:cosheaf}
Let $H\in \fOpen(\cH)$ and let $H= \bigcup_{i\in I} H_i$ be a covering by disjoint subsets of $\fOpen(\cH)$. If there exists $L\in \sL$ such that $H_i$ is $L$-stable for every $i\in I$ then \eqref{intloccosheaf} and \eqref{cdsheaf} are isomorphisms. 
\end{lemma}

\begin{proof} It suffices to prove the surjectivity of \eqref{intloccosheaf}. For a given $f\in \Int_{\loc, b}(H, \bZ)$ we can choose $L\in \sL$ sufficiently small such that all $H_i$ are $L$-stable and such that $f\in \Int_{L, b}(H, \bZ)$. Now the assertion follows from Remark \ref{remarks:phpolfunc} (d).
\end{proof}

\begin{remark}
\label{remark:cosheaflocpoldistb} 
\rm Note that under the assumptions of the Lemma the collection of maps $\cD_{\lpol, b}(H_i, R) \to \cD_{\lpol, b}(H, R), \mu \mapsto \mu_!$ for $i\in I$ 
induces an isomorphism  
\begin{equation*}
\label{dpolbcosheaf}
\bigoplus_{i\in I}\cD_{\lpol, b}(H_i, R) \lra \cD_{\lpol, b}(H, R).
\end{equation*} 
\end{remark}

To discuss further functorial properties of the above Definitions \ref{df:locpolfunct} we introduce the notion of an $\sL$-affine map.

\begin{df} 
\label{df:affineslmaps}
Let $\cH_1$, $\cH_2$ be $\Lambda$-sets. A map $\varphi:\cH_1 \to \cH_2$ will be called an $\sL$-affine map if there exists $\alpha\in \Aut_{\sL}(\Lambda)$ such that $\varphi(\la + h) = \alpha(\la) + \varphi(h)$ for every $\la\in \Lambda$ and $h\in \cH_1$.
\end{df}

\begin{remark}
\label{remark:funct}
\rm Let $\cH =\Lambda$ be equipped with the obvious structure as a $\Lambda$-set. An $\sL$-affine map $\varphi:\Lambda\to \Lambda$ is an automorphism composed with a translation. Thus every $\sL$-affine map $\varphi:\Lambda\to \Lambda$ is a bijection. Hence the $\sL$-affine maps $\varphi:\Lambda\to \Lambda$ form a group -- denote by $\Aff_{\sL}(\Lambda)$ -- and we have $\Aff_{\sL}(\Lambda)=\Aut_{\sL}(\Lambda)\ltimes \Lambda$.
\end{remark}

Let $\varphi:\cH_1\to \cH_2$ be an $\sL$-affine map between $\Lambda$-sets. Note that for $H \in \fOpen(\cH_2)$ we have $\varphi^{-1}(H)\in \fOpen(\cH_1)$ and that $\varphi^*: \Maps(\varphi^{-1}(H), \bZ) \to \Maps(H, \bZ), f\mapsto f\circ \varphi$ maps the subring $\Int_{\loc}(\varphi^{-1}(H), \bZ)$ into $\Int_{\loc}(H, \bZ)$. Now assume that $\varphi:\cH_1\to \cH_2$ is injective. In this case we have $\varphi^{-1}(H)\in \fOpen_{\fin}(\cH_1)$ for every $H\in \fOpen_{\fin}(\cH_2)$. It follows that the homomorphism $\varphi^*: \Maps(\varphi^{-1}(H), \bZ) \to \Maps(H, \bZ)$ for $H \in \fOpen(\cH_2)$ maps the subring $\Int_{\loc, b}(\varphi^{-1}(H), \bZ)$ into $\Int_{\loc, b}(H, \bZ)$. By passing from $\varphi^*: \Int_{\loc, b}(\varphi^{-1}(H), \bZ)\to \Int_{\loc, b}(H, \bZ)$ to duals we obtain a map  
\begin{equation}
\label{funct2}
\varphi_*: \cD_{\loc}(H, R) \lra \cD_{\loc}(\varphi^{-1}(H), R).
\end{equation} 

Now assume that $\cH= \Lambda$ and let $H=L\in \sL$. Then $\cD_{\lpol}(L, R)$ carries a natural ring structure given by the convolution \eqref{conv} and the map
\begin{equation}
\label{respol}
\cD_{\lpol}(L, R) \lra \cD_{\pol}(L, R)
\end{equation} 
dual to the inclusion $\Int(L, R)\hra \Int_{\loc}(L, R)$ is a ring homomorphism. Since the Dirac distributions $\delta_{\la}$, $\la\in L$ lie in $\cD_{\lpol}(L, R)$ the map $L\to \cD_{\lpol}(L, R), \la\mapsto \delta_{\la}$ extends to an $R$-algebra homomorphism $R[L]\to \cD_{\lpol}(L, R)$ and \eqref{respol} is a homomorphism of $R[L]$-algebras. 

The $R$-module $\cD_{\lpol}(\Lambda, R)$ is equipped with a left action of the group $\Aff_{\sL}(\Lambda)$ defined by the 
homomorphisms \eqref{funct2}, i.e.\ it is equipped with a natural $R[\Aff_{\sL}(\Lambda)]$-module structure. Similarly one defines an $R[\Aff_{\sL}(\Lambda)]$-module structure on $\cD_{\lpol, b}(\Lambda, R)$ so that the canonical map $\cD_{\lpol,b }(\Lambda, R)\to \cD_{\lpol}(\Lambda, R)$
is $\Aff_{\sL}(\Lambda)$-equivariant. 

Now assume that $V$ is oriented. This provides also each lattice $L\subset V$ with the structure of an oriented lattice so we can apply Corollary \ref{coro:koszul}. The $\Aff(L)$-action on $\Ext_{R[L]}^n(R, \cD_{\lpol}(L, R)) \cong R$ is not trivial but is given by the {\it sign character}, i.e.\ the homomorphism
\begin{equation}
\label{sign}
\ep: \Aff(V)  \lra \{\pm 1\} = \bZ^*, \quad (\alpha, v)\mapsto  \sgn(\det(\alpha)).
\end{equation}
Thus Corollary \ref{coro:koszul} can be restated as 
\begin{equation*}
\label{koszul2}
\Ext_{R[L]}^i(R, \cD_{\pol}(L, R)(\ep))\, =\, \left\{ \begin{array}{cc} R & \mbox{if $i=n$,}\\
                                                      0 & \mbox{if $i \ne n$.}
                                                      \end{array}\right.
\end{equation*}
when taking into account the $\Aff(L)$-actions. Similarly, we have

\begin{prop}
\label{prop:koszul2}
For every $i\ge 0$,
\[
\Ext_{R[\Lambda]}^i(R, \cD_{\lpol}(\Lambda, R)(\ep))\, =\, \left\{ \begin{array}{cc} R & \mbox{if $i=n$,}\\
                                                      0 & \mbox{if $i \ne n$.}
                                                      \end{array}\right.
\]
as $R[\Aff_{\sL}(\La)]$-modules.
\end{prop}

For the proof we need the following obvious

\begin{lemma}
\label{lemma:basechange}
For a ring $R$ and $L', L\in \sL$ with $L'\subseteq L$ the canonical map
\begin{equation*}
\label{resext1}
\bZ[L]\otimes_{\bZ[L']} \Int(L', \bZ) \lra \Int_{L'}(L, \bZ), \quad \alpha\otimes f \mapsto  \alpha \star f_{!}
\end{equation*}
and its dual
\begin{equation*}
\label{resext3}
\cD_{\pol, L'}(L, R) \lra \Hom_{R[L']}(R[L], \cD_{\pol}(L', R)) 
\end{equation*}
are isomorphisms.
\end{lemma}

\begin{proof}[Proof of Prop.\ \ref{prop:koszul2}] Firstly, note that by \eqref{locpoldist2} we have 
\begin{equation*}
\label{locpoldistring}
\cD_{\lpol}(\Lambda, R) \, =\, \prolim_{L\in \sL^{\opp}} \cD_{\lpol}(L, R).
\end{equation*}
We fix a lattice $L_0$ in $\sL$. The homomorphism of $R[L_0]$-modules $\cD_{\lpol}(L_0, R) \to \cD_{\lpol}(\Lambda, R)$ (see \eqref{distext0}) induces a homomorphism of $R[\Lambda]$-modules 
\begin{equation}
\label{coind}
\cD_{\lpol}(\Lambda, R) \,\lra \, \Hom_{R[L_0]}(R[\Lambda], \cD_{\lpol}(L_0, R)). 
\end{equation}
It is an isomorphism by Lemma \ref{lemma:cosheaf}. Indeed if $\cR$ is a system of representatives for the $L_0$-cosets in $\Lambda$ then the map \eqref{coind} can be identified with the map \eqref{cdsheaf} for the covering $\Lambda = \bigcup_{\la\in \cR} \la +L_0$. Moreover by \eqref{locpoldist} and Lemma \ref{lemma:basechange} we have 
\begin{equation*}
\label{koszul5}
\cD_{\lpol}(L_0, R) \, =\, \prolim_{L\in \sL(L_0)^{\opp}} \cD_{\pol, L}(L_0, R)\, =\, \prolim_{L\in \sL(L_0)^{\opp}} \Hom_{R[L]}(R[L_0], \cD_{\pol}(L, R))
\end{equation*}
hence 
\begin{eqnarray}
\label{koszul6}
\cD_{\lpol}(\Lambda, R) & \cong & \Hom_{R[L_0]}(R[\Lambda], \prolim_{L\in \sL(L_0)^{\opp}} \Hom_{R[L]}(R[L_0], \cD_{\pol, L}(L, R)))\\
& \cong &\prolim_{L\in \sL^{\opp}}\Hom_{R[L]}(R[\Lambda], \cD_{\pol}(L, R)).\nonumber
\end{eqnarray}
By Cor.\ \ref{coro:koszul} and Shapiro's Lemma we have
\begin{equation}
\label{koszul7}
\Ext_{R[\Lambda]}^i(R, \Hom_{R[L]}(R[\Lambda], \cD_{\pol}(L, R))) \, = \,  \Ext_{R[L]}^i(R, \cD_{\lpol}(L, R)) \, =\, \left\{ \begin{array}{cc} R(\ep) & \mbox{if $i=n$,}\\
                                                      0 & \mbox{if $i \ne n$}
                                                      \end{array}\right.
\end{equation}
for every $i\ge 0$ and $L\in \sL$. On the other hand \eqref{koszul6} and the first equality in \eqref{koszul7} together with (\cite{weibel}, Thm.\ 3.5.8) imply that there exists a short exact sequence 
\begin{eqnarray}
\label{koszul8}
&& 0 \lra \prolim_{L\in \sL^{\opp}}^{(1)} \Ext_{R[L]}^{i-1}(R, \cD_{\lpol}(L, R))(\ep) \lra \Ext_{R[\Lambda]}^i(R, \cD_{\lpol}(\Lambda, R)(\ep))\\
&& \hspace{6cm}\lra \prolim_{L\in \sL^{\opp}} \Ext_{R[L]}^i(R, \cD_{\lpol}(L, R))(\ep) \lra 0.\nonumber
\end{eqnarray}
From the second equality in \eqref{koszul7} we deduce that the first term in \eqref{koszul8} vanishes for every $i$ and the third for every $i$ except for $i=n$ when it is $\cong R$.
\end{proof}

\paragraph{Change of $\sL$} We also need to investigate the effects of changing the subset $\sL\subseteq \Lat(V)$ in specific circumstances. For that let $\sL_1$ and $\sL_2$ be non-empty closed subsets of $\Lat(V)$ and put $\Lambda_i:= \Lambda(\sL_i)$ for $i=1,2$. 
We assume that $\Lambda_2 \cap L_1 \in \sL_2$ for every $L_1\in \sL_1$ and that 
\begin{equation*}
\label{changesl}
\sL_1\lra \sL_2, \quad L_1\mapsto \Lambda_2 \cap L_1
\end{equation*}
is a bijection. This implies in particular that $\Lambda_1\supseteq \Lambda_2$. We consider a $\Lambda_1$-set $\cH_1$ and a $\Lambda_2$-subset $\cH_2\subseteq \cH_1$. Let $H_1\in \fOpen(\cH_1)$ and let $L_1\in \sL_1$ such that $H_1$ is $L_1$-stable. Then $H_2:=H_1\cap \cH_2$ is $L_2:=L_1\cap\Lambda_1$-stable. Moreover if $f\in \Int_{L_1}(H_1, \bZ)$ (resp.\ 
$f\in \Int_{L_1, b}(H_1, \bZ)$) then $f|_{H_2} \in \Int_{L_2}(H_2, \bZ)$ (resp.\ $f|_{H_2}\in \Int_{L_2, b}(H_2, \bZ)$). 

It follows that the map 
\begin{equation}
\label{slchangeres1}
\Int_{\sL_1\mloc}(H_1, \bZ) \lra \Int_{\sL_2\mloc}(H_2, \bZ), \qquad f \mapsto f|_{H_2} 
\end{equation} 
is well-defined for every $H_1\in \fOpen(\cH_1)$ and $H_2\in \fOpen(\cH_2)$ with $H_2\subseteq H_1$. Moreover \eqref{slchangeres1} maps the subring $\Int_{\sL_1\mloc, b}(H_1, \bZ)$ of $\Int_{\sL_1\mloc}(H_1, \bZ)$ into $\Int_{\sL_2\mloc, b}(H_2, \bZ)$. Thus, dually, we obtain maps 
\begin{eqnarray}
\label{slchangeext1}
&& \cD_{\sL_2\mlpol, b}(H_2, R) \lra \cD_{\sL_1\mlpol, b}(H_1, R), \qquad \mu \mapsto \mu_! \\
&& \cD_{\sL_2\mlpol}(H_2, R) \lra \cD_{\sL_1\mlpol}(H_1, R), \qquad \mu \mapsto \mu_! 
\label{slchangeext2}
\end{eqnarray} 
that are characterized by 
\begin{equation*}
\label{resext}
\int_{H_1} f(h_1) \, d\mu_! (h_1) \,=\, \int_{H_2} (f|_{H_2})(h_2)\, d\mu(h_2)
\end{equation*}
for all $\mu\in \cD_{\sL_2\mlpol, b}(H_2, R)$ (resp.\ $\mu\in \cD_{\sL_2\mlpol}(H_2, R)$) and $f\in \Int_{\sL_1\mloc}(H_1, \bZ)$ (resp.\ $f\in \Int_{\sL_1\mloc, b}(H_1, \bZ)$). In particular for $H_1=\cH_1$ and $H_2=\cH_2$ we obtain maps
\begin{eqnarray}
\label{slchangeext1a}
&& \cD_{\sL_2\mlpol, b}(\cH_2, R) \lra \cD_{\sL_1\mlpol, b}(\cH_1, R), \qquad \mu \mapsto \mu_! \\
&& \cD_{\sL_2\mlpol}(\cH_2, R) \lra \cD_{\sL_1\mlpol}(\cH_1, R), \qquad \mu \mapsto \mu_!.
\label{slchangeext2a}
\end{eqnarray} 
Under certain conditions on $\sL_1$, $\sL_2$, $\cH_1$, $\cH_2$ and $R$ these are isomorphisms. 

\begin{prop}
\label{prop:changesl} 
Assume that 
\medskip

\noi (i) The index $d:=[L_1 : \Lambda_2 \cap L_1]$ is independent of the choice of $L_1\in \sL_1$.
\medskip

\noi (ii) The action of $\La_i$ on $\cH_i$ is faithful for $i=1,2$ and the map $\Lambda_2\backslash \cH_2 \to \Lambda_1\backslash \cH_1, \Lambda_2+h\mapsto \Lambda_1+h$ is bijective.
\medskip

\noi (iii) The index $d=[\Lambda_1 : \Lambda_2]$ is invertible in $R$.
\medskip

\noi Then the maps \eqref{slchangeext1} and \eqref{slchangeext2} for $H_1=\cH_1$ and $H_2=\cH_2$
are isomorphisms.
\end{prop} 

\begin{proof} Note that (i) implies $d=[\Lambda_1 : \Lambda_2]$. By Lemma \ref{lemma:cosheaf} and Remark \ref{remark:cosheaflocpoldistb} it suffices to consider the case when $\Lambda_2$ acts freely and transitively on $\cH_2$ (hence $\Lambda_1$ acts freely and transitively on $\cH_1$ as well by (ii)). Thus we may assume $\cH_i=\Lambda_i$ for $i=1,2$. 

We fix $L_1\in \sL_1$ and put $L_2=L_1\cap \Lambda_2$. As remarked in the proof of Prop.\ \ref{prop:koszul2} we have 
\begin{equation*}
\label{coind2}
 \cD_{\sL_i\mlpol}(\Lambda_i, R) \,\cong \, \Hom_{R[L_i]}(R[\Lambda_i], \cD_{\sL_i\mlpol}(L_i, R))
\end{equation*}
for $i=1, 2$. Similarly, using Remark \ref{remark:cosheaflocpoldistb} one can show that 
\begin{equation*}
\label{ind1}
 \cD_{\sL_i\mlpol, b}(\Lambda_i, R) \,\cong R[\Lambda_i]\otimes_{R[L_i]} \cD_{\sL_i\mlpol}(L_i, R).
\end{equation*}
Thus it suffices to show that the map \eqref{slchangeext1} for $H_1=L_1$, $H_2=L_2$, i.e.\ the map 
\begin{equation*}
\label{distext2}
\cD_{\sL_2\mlpol}(L_2, R) \lra \cD_{\sL_1\mlpol}(L_1, R), \quad \mu \mapsto \mu_! 
\end{equation*} 
is an isomorphism. 

Firstly, note that (iii) and Lemma \ref{lemma:convol} (a) imply that $\cD_{\pol}(L_2, R)\to \cD_{\pol}(L_1, R), \mu\mapsto \mu_!$ is an isomorphism. Note also that the map $\sL_1(L_1)\to \sL_2(L_2), L_1'\mapsto L_2 \cap L_1'$ is a bijection. Moreover note that we have by \eqref{locpoldist} we have 
\begin{equation*}
\label{distext3}
\cD_{\sL_i\mlpol}(L_i, R) \, =\, \prolim_{L_i'\in \sL(L_i)^{\opp}} \cD_{\pol, L_i'}(L_i, R)
\end{equation*}
for $i=1, 2$. Thus if we fix $L_1'\in \sL_1(L_1)$ and put $L_2':=L_2\cap L_1'$ then it suffices to show that the canonical map
\begin{equation}
\label{distext4}
\cD_{\pol, L_2'}(L_2, R) \lra \cD_{\pol, L_1'}(L_1, R)
\end{equation} 
is an isomorphism. From (i) we deduce that $R[L_2]\otimes_{R[L_2']} R[L_1']\cong R[L_1]$. Since -- again by Lemma \ref{lemma:convol} (a) -- the map $\cD_{\pol}(L_2', R)\to \cD_{\pol}(L_1', R), \mu\mapsto \mu_!$ is bijective, it induces an isomorphism   
\begin{eqnarray}
\label{distext5}
 \Hom_{R[L_2']}(R[L_2], \cD_{\pol}(L_2', R)) &  \stackrel{\cong}{\lra} &  \Hom_{R[L_2']}(R[L_2], \cD_{\pol}(L_1', R))\\
& \cong & \Hom_{R[L_1']}(R[L_2]\otimes_{R[L_2']} R[L_1'], \cD_{\pol}(L_1', R))\nonumber\\
& \cong & \Hom_{R[L_1']}(R[L_1], \cD_{\pol}(L_1', R))\nonumber
\end{eqnarray}
According to Lemma \ref{lemma:basechange} it can be identified with \eqref{distext4}. 
\end{proof}

We apply Prop.\ \ref{prop:changesl} in the following situation. As in Example \ref{examples:fracideals} (b) let $F$ be a number field of degree $d$ over $\bQ$, $V$ an $m$-dimensional $F$-vector space, $S$ a finite set of nonarchimedean places of $F$ and $\cM\subseteq V$ a finitely generated $\cO_S$-submodule of $V$ with $\rank_{\cO_S} \cM = \dim_F(V)$. Let $\sL\subseteq \Lat_{\cO_F}(V)$ be the set of $L\in \Lat_{\cO_F}(V)$ with $L\subseteq \cM$ and $L\otimes_{\cO_F} \cO_S= \cM$, i.e.\ $L$ generates $\cM$ as an $\cO_S$-module. 

\begin{coro}
\label{coro:ordtinv}
Let $v\in V$ and let $h$ be the order of $v+\cM$ in $V/\cM$, i.e.\ $h$ is the minimal positive integer such that $h v\in \cM$. If $h$ is invertible in $R$ then there exists natural isomorphisms 
\begin{equation*}
\label{slchangeextstalk}
 \cD_{\sL\mlpol, b}(v+ \cM, R) \,\cong\, \cD_{\sL\mlpol, b}(\cM, R), \qquad \cD_{\sL\mlpol}(v+ \cM, R) \,\cong\, \cD_{\sL\mlpol}(\cM, R)
 \end{equation*} 
\end{coro}

\begin{proof} Let $\cM'\subseteq V$ be another finitely generated $\cO_S$-submodule of $V$ with $\cM'\supseteq \cM$ and put $\sL'=\{L'\in \Lat_{\cO_F}(V)\mid L'\subseteq \cM', L'\otimes_{\cO_F} \cO_S = \cM'\}$. We remark that the map
$\sL' \lra \sL, \quad L'\mapsto L'\cap \cM$ is a bijection and that $L'/L \to \cM'/\cM, \la + L\mapsto \la + \cM'$
is an isomorphism for every $L'\in \sL'$ and $L:=L'\cap \cM$ (these facts can be easily seen using the adelic descriptions of the sets $\sL'$ and $\sL$ given in \ref{examples:fracideals} (b)). If we choose $\cM' := \frac{1}{h} \cM : = \{w\in V\mid h w\in \cM\}$ then we see that conditions (i) and (ii) of Prop.\ \ref{prop:changesl} hold for $\sL_1= \sL'$, $\sL_2:= \sL$, $\cH_1 = \cM'$ and $\cH_2:= v+ \cM$. Moreover the index $[\cM':\cM]$ is a divisor of $h^n$ with $n=dm$ hence is invertible in $R$ so (iii) holds as well. Thus we can apply $\ref{prop:changesl}$ and obtain
\[
\cD_{\sL\mlpol, b}(v + \cM, R) \,\cong\, \cD_{\sL'\mlpol, b}(\cM', R), \qquad \cD_{\sL\mlpol}(v + \cM, R) \,\cong\, \cD_{\sL'\mlpol}(\cM', R).
\]
Since $v\in \frac{1}{h} \cM$ was arbitrary the assertion follows.
\end{proof}

\begin{remark}
\label{remark:ordVmodcM}
\rm Note that the only possible prime divisors of $h$ are the primes numbers lying below the places in $S$. 
\end{remark}

\section{Lattice topology, sheaves and cohomology}
\label{section:latspaces}

\paragraph{$(\sL, \wGa)$--spaces and lattice topology} As in section \ref{section:locpol} we fix a finite-dimensional $\bQ$-vector space $V$ and a non-empty subset $\sL\subseteq \Lat(V)$ with the property \eqref{sLclosed}. Recall that $\Lambda = \Lambda(\sL)=\bigcup_{L\in \sL} L$ is a subgroup of $V$. We also fix a subgroup $\Gamma\subseteq \GL(V)$ such that $\sL$ is $\Gamma$-stable and put $\wGa:= \Gamma\ltimes \Lambda\subseteq \Aff_{\bZ}(\Lambda)$. Elements of $\wGa$ will be denoted by $\wga= (\ga, \la)$ where $\ga\in \Gamma$ and $\la\in \Lambda$. 

We associate to this data the following category $\cC(\sL, \wGa)$. Its set of objects is $\sL$. For $L_1$, $L_2\in \sL$ a morphism $\varphi: L_1 \to L_2$ is a triple $\varphi=(\wga, L_1, L_2)\in \Gamma\times \sL\times \sL$ with $\gamma(L_1) \subseteq L_2$  if $\wga= (\gamma, \la)$. The composition of two morphisms $\varphi_1=(\wga_1, L_1, L_2)$ and $\varphi_2=(\wga_2, L_2, L_3)$ is the morphism $\varphi_2\circ \varphi_1:= (\wga_2\cdot\wga_1, L_1, L_3)$. Given a morphism $\varphi = (\wga, L_1, L_2)$ 
with $\wga= (\gamma, \la)$ we define its degree by 
\begin{equation*}
\label{deg}
\deg(\varphi):=[L_2:\gamma(L_1)].
\end{equation*}
The degree is multiplicative, i.e.\ we have $\deg(\varphi_2\circ \varphi_1) = \deg(\varphi_2)\cdot \deg(\varphi_1)$ for two morphism 
$\varphi_1: L_1 \to L_2$ and $\varphi: L_2 \to L_3$ in $\cC(\sL, \wGa)$. Morphisms of the form $\varphi=(1, L_1, L_2)$ (i.e.\ when the first component of $\varphi$ is the neutral element in $\Gamma$) will be denoted by 
\begin{equation}
\label{equiv3}
\pi =\pi_{L_1, L_2}: L_1 \hra L_2.
\end{equation}
For these morphisms we have $\deg(\pi_{L_1, L_2}) = [L_2: L_1]$. For $L\in \sL$ and $\wga=(\gamma, \la)\in \Gamma$ we denote the morphism $(\wga, L , \gamma(L))$ by 
\begin{equation}
\label{equiv4}
\wga_*: L \lra \gamma(L).
\end{equation}
Consider the special case when $\gamma=1$, i.e.\ when $\wga=(1, \la)$. Then the morphism \eqref{equiv4} will be denoted by 
\begin{equation*}
\label{equiv4a}
\la_*: L \lra L.
\end{equation*}
We note that an arbitrary morphism $\varphi=(\wga, L_1, L_2)$ with $\wga=(\gamma, \la)\in \Gamma$ can be factored as a composition of morphisms of the type \eqref{equiv3} and \eqref{equiv4}. Indeed, we have
\begin{equation*}
\label{equiv5}
\varphi =  \wga_* \circ \pi_{L_1, \gamma^{-1}(L_2)} = \pi_{\gamma(L_1), L_2} \circ \wga_* .
\end{equation*}

We consider a functor $X: \cC(\sL, \wGa)\to \spaces$. For $L\in \sL$ we will write $X_L$ instead of $X(L)$. The image of a morphism $\varphi: L_1\lra L_2$ under $X$ will be denote by $\varphi: X_{L_1} \to X_{L_2}$ as well. In particular for $L, L_1, L_2\in \sL$ with $L_1\subseteq L_2$ and $\wga=(\gamma, \la)\in \Gamma$ we obtain morphisms
\begin{equation*}
\label{equiv6}
\pi =\pi_{L_1, L_2}: X_{L_1} \lra X_{L_2}\qquad \mbox{and}\qquad \wga_*: X_L\lra X_{\gamma(L)}.
\end{equation*}

\begin{df}
\label{df:sspace}
(a) A functor $X: \cC(\sL, \wGa)\to \spaces$ will be called a $(\sL, \wGa)$-space if the map $\varphi: X_{L_1}\to X_{L_2}$ is a covering of degree $\deg(\varphi)$ for every morphism $\varphi: L_1\to L_2$ in $\cC(\sL, \wGa)$. 
\medskip

\noi (b) A $(\sL, \Ga)$-space is a $(\sL, \wGa)$-space such that the map $X_L \to X_L$ induced by the morphism $\lambda_*: L\to L$ is the identity for every $L\in \sL$ and $\la\in \Lambda$. 
\medskip

\noi (c) A morphism of $(\sL, \wGa)$-spaces $f: X\to Y$ is a morphism of functors such that for every morphism $\varphi: L_1\to L_2$ in $\cC(\sL, \wGa)$ the diagram 
\[
\begin{CD}
X_{L_1} @> f_{L_1} >> Y_{L_1}\\
@VV \varphi V @VV \varphi V\\
X_{L_2} @> f_{L_2} >> Y_{L_2}
\end{CD}
\]
is cartesian. The set of morphisms $X\to Y$ will be denoted by $\Hom_{\wGa}(X,Y)$.
\end{df}

For a $(\sL, \wGa)$-space $X$ we are going to associate a certain site which allows us to consider $\wGa$-equivariant sheaf on $X$. For that we put 
\begin{equation*}
\label{Xhat}
\hX : = \prolim_{L\in \sL} X_L
\end{equation*}
i.e.\ as a set $\hX$ is the projective limit of the inverse system consisting 
of the sets $X_L$ for $L\in \sL$ and the maps $\pi_{L_1, L_2}$ for $L_1, L_2\in \sL$ with $L_1\subseteq L_2$. 
For $L\in \sL$ we let $\pi_L: \hX \to X_L$ denote the canonical projection. A subset $S\subseteq \hX$ will be called {\it $L$-stable} if it is of the form $S =\pi_L^{-1}(S')$ for a subset $S'\in X_L$. The set $\hX$ is equipped with a natural $\wGa$-action induced by the collection of maps $\{\wga_*:X_L\to X_{\gamma(L)}\}_{L\in \sL, \wga\in \wGa}$. Indeed, given $\wga=(\gamma, \la)\in \wGa$ we define the map $\wga\cdot : \hX \to \hX, x\mapsto \wga \cdot x
$ by
\begin{equation}
\label{Xhga}
\wga\cdot  : =\prolim_L (\wga, L, \gamma(L)): \hX  = \prolim_{L\in \sL} X_L \lra \prolim_{L\in \sL} X_{\gamma(L)} = \hX.
\end{equation}
Rather than equipping $\hX$ with the projective limit topology we consider a much coarser notion of open sets and coverings on $\hX$. Namely, 
a subset $U$ of $\hX$ will be called {\it open} if and only if there exists $L\in \sL$ and an open subset $W$ of $X_L$ with $U=\pi_L^{-1}(W)$. Equivalently, $U$ is open if and only if it is $L$-stable for some $L\in \sL$ and if $\pi_L(U)$ is open in $X_L$. We let $\fOpen(\hX)$ be the collection of all open subsets of $\hX$. For an open set $U$ of $\hX$, a family of open subsets $\{U_i\}_{i\in I}$ of $U$ will be called a covering of $U$ if $\{U_i\}_{i\in I}$ is a covering in the naiv sense (i.e.\ we have $U = \bigcup_{i\in I} U_i$) and if there exists $L\in \sL$ such that every $U_i$ is $L$-stable. Equivalently, there exists an open subset $W$ of $X_L$ and an open covering $\{W_i\}_{i\in I}$ of $W$ such that  $U=\pi_L^{-1}(W)$ and $U_i=\pi_L^{-1}(W_i)$ for all $i\in I$. The collection of all coverings of $U$ will be denoted by $\Cov(U)$. 

\begin{lemma}
\label{lemma:gtop}
The triple $(\hX, \fOpen(X), \Cov)$ is a site (in the sense of Def.\ \ref{df:coverage} (b) of the appendix).
\end{lemma}

\begin{proof} One easily checks that for open subsets $U, V\subseteq \hX$ both $U\cup V$ and $U\cap V$ are open and that 
for $\{U_i\}_{i\in I}\in \Cov(U)$ we have $\{U_i\cap V\}_{i\in I}\in \Cov(U\cap V)$. 
\end{proof}

We refer to $\hX=(\hX, \fOpen(X), \Cov)$ as the {\it adelic space} associated to $X$ and say that $\hX$ is equipped with the {\it lattice topology}. Note that the map $\pi_L: \hX\to X_L$ is continuous as a morphisms between sites for every $L\in \sL$ (in the sense of Def.\ \ref{df:coverage} (c)) and that we have chosen $\fOpen(X)$ as well as $\Cov$ to be minimal with this property. Note also that the map \eqref{Xhga} is continuous for every $\wga\in \wGa$. Thus $\hX$ is equipped with a continuous $\wGa$-action. Moreover if $X$ is a $(\sL, \Ga)$-space then 
the $\Lambda$-action on $\hX$ is trivial so that $\hX$ is just equipped with a $\Ga$-action.

\begin{remarks}
\label{remarks:nogtop}
\rm (a) We remark that in general the lattice topology is not a topology in the usual sense, i.e.\ the site $(\hX, \fOpen(X), \Cov)$ is not a topological space. In fact although the union of a finite collection of open subsets of $\hX$ is open as well, this does not hold for the union of an arbitrary collection of open subsets. Moreover $(\hX, \fOpen(X), \Cov)$ is in general not a site in the sense of (\cite{stacks}, Tag 00VH). Indeed, albeit conditions (3) and (4) of (\cite{FvdP}, Def.\ 2.4.1) hold, condition (5) usually does not. Namely, if an open subset $U\subseteq \hX$ together with $\{U_i\}_{i\in I}\in \Cov(U)$ and $\{U_{ij}\}_{j\in J_i}\in \Cov(U_i)$ for each $i\in I$ are given, then the collection of subsets $\{U_{ij}; i\in I, j\in J_i\}$ of $U$ is in general not a covering. 
\medskip

\noi (b) We note that a morphism $f: X\to Y$ of $(\sL, \wGa)$-spaces induces a $\wGa$-equivariant continuous morphism 
\begin{equation*}
\label{fhat}
\hf := \prolim_{L\in \sL} f_L : \hX \lra \hY.
\end{equation*}
\end{remarks}

The obvious way to produce examples of $(\sL, \wGa)$-spaces is as follows. Let $\cX$ be a locally compact Hausdorff space equipped with a free and continuous left $\wGa$-action (i.e.\ every $\wga\in \wGa$ acts as a homeomorphism on $\cX$). We denote the action of the subgroup $\La\subseteq \wGa$ on $\cX$ additively. Assume that every $L\in \sL$ -- viewed as a subgroup of $\La\subseteq \wGa$ -- acts properly discontinuously on $\cX$. For $L\in \sL$ we let $X_L =\{ L + x\mid x\in \cX\}$ be the set of $L$-orbits in $\cX$ equipped with the quotient topology. For a morphism $\varphi=(\wga, L_1, L_2)$ in $\cC(\sL, \wGa)$ with $\wga=(\ga, \la)\in \wGa$ we consider the induced map $\varphi_*: X_{L_1} \to X_{L_2}, \quad x + L_1\mapsto \wga\cdot x + L_2$. It is easy to see that the assignment $L\mapsto X_L, \varphi\mapsto \varphi_*$ is a $(\sL, \wGa)$-space $X$. Also any $\wGa$-equivariant continuous map $f: \cX\to \cY$ between locally compact Hausdorff space with such $\wGa$-actions induces a morphism of the associated $(\sL, \wGa)$-spaces. 

In this example the site $(\hX, \fOpen(X), \Cov)$ admits a rather simple description. Namely the set $\fOpen(X)$ can be identified with the collection of open subsets $U\subseteq \cX$ that are stable under the action of some $L\in \sL$ (i.e.\ we have $L+U=U$). Moreover the coverings of the site can be identified with open coverings $U = \bigcup_{i} U_i$ in the usual sense such that there exists $L\in \sL$ so that every $U_i$ is $L$-stable. 

In principle all the specific examples of $(\sL, \wGa)$-spaces considered below and that are relevant to us arise in this way. However in the construction of the adelic Eisenstein classes in section \ref{section:eiscocycle} we have to consider certain morphism between $(\sL, \wGa)$-spaces that do not stem from $\wGa$-equivariant continuous maps $f: \cX\to \cY$ as discussed above. Also we will later consider the notion of stalks for a sheaf on the site $\hX$ at truly adelic points. For these reasons it is not sufficient to work entirely within the realm of spaces $\cX$ with a $\wGa$-action as described above. Our construction of Eisenstein classes is of a genuine adelic nature. 

\begin{examples}
\label{examples:torpt} \rm (a)  We consider the above example for $\cX=\La$ equipped with the discrete topology and the obvious $\wGa$-action. Thus for $L\in \sL$ we consider the quotient $B_L:= \Lambda/L$ as a discrete space and for a morphism 
$\varphi=\varphi_{\wga, L_1, L_2}: L_1\to L_2$ in $\cC(\sL, \wGa)$ we consider the map
\begin{equation*}
\label{equiv1a}
\varphi=\varphi_{\wga, L_1, L_2}:B_{L_1}=\Lambda/L_1\lra B_{L_2}=\Lambda/L_2, \quad \la + L_1 \mapsto \wga(\la) + L_2.
\end{equation*}
This is the simplest example of a $(\sL, \wGa)$-space. It will be denoted by $B = B^{\sL}$ and will be called the {\it basic $(\sL, \wGa)$-space}.
Note that we have 
\begin{equation}
\label{Bhat}
\hB \,= \,\prolim_{L\in \sL} \Lambda/L \, = : \hLa.
\end{equation}
We describe the lattice topology on $\hLa$. For $L\in \sL$ define $\hL = \hL^{\sL}:=\ker(\pi_L: \hB=\hLa\to B_L=\Lambda/L)$ so that 
\begin{equation*}
\label{Lhat}
\hL \,:=\, \prolim_{L'\in \sL, L'\subseteq L} L/L'.
\end{equation*} 
A subset $U \subseteq \hLa$ is open if and only if there exists $L\in \sL$ such that $U$ is a union of $\hL$-cosets. Equivalently, $U$ is open if it is $\hL$-stable for some $L\in \sL$, i.e.\ we have $x + U = U$ for all $x\in \hL$. A covering of an open subset $U\subseteq \hLa$ consists of a covering $\bigcup_{i\in I} U_i = U$, so that there exists $L\in \sL$ such that $U_i$ is $\hL$-stable for all $i\in I$. Note that since $\La/L\cong \hLa/\hL$, an $\hL$-stable subset of $\hLa$ is necessarily of the form $\cU + \hL= \bigcup_{\la\in \cU} \la + \hL$ for a unique $L$-stable subset $\cU\subseteq \Lambda$ (namely $\cU=U\cap \Lambda$). Thus the map
\begin{equation*}
\label{Bhatopen}
\fOpen(B) \lra \{\cU\subseteq \La\mid \mbox{$\cU$ is $L$-stable for some $L\in \sL$}\,\}, \quad U \mapsto U\cap \Lambda
\end{equation*}
is bijective. 

\medskip

\noi (b) If $X$ is $(\sL, \wGa)$-space and $Y$ a topological space equipped with a $\Gamma$-action then we denote by $X\times Y$ the following $(\sL, \wGa)$-space. For $L\in \sL$ we put $(X\times Y)_L:= X_L\times Y$. Also for a morphism 
$\varphi=\varphi_{\wga, L_1, L_2}: L_1\to L_2$ in $\cC(\sL, \wGa)$ with $\wga=(\ga, \la)$ we define the induced map $\varphi: (X\times Y)_{L_1} \to (X\times Y)_{L_2}$ by 
\begin{equation*}
\label{equiv1b}
\varphi: X_{L_1}\times Y \lra X_{L_2}\times Y, \quad (x, y) \mapsto (\varphi(x), \gamma\cdot y).
\end{equation*}
Note that the projection onto the first factor $\pr_X: X\times Y\to X$ is a morphism of $(\sL, \wGa)$-spaces.
Also if $Y_1$ and $Y_2$ are topological spaces equipped with a $\Gamma$-action and if $f: Y_1 \to Y_2$ is a continuous $\Gamma$-equivariant map  then $\id\times f: X\times Y_1\to X\times Y_2$ is a morphism of $(\sL, \wGa)$-spaces.

For example if we equip $\Lambda_{\bR} = V_{\bR}$ with usual topology and the obvious $\Gamma$-action then the product 
\begin{equation}
\label{covtorus}
A\,=\,A^{\sL}\,:=\, B\times V_{\bR}: \cC(\sL, \wGa)\lra \spaces, \quad L\mapsto A_L = B_L \times V_{\bR}=  \Lambda/L\times V_{\bR}
\end{equation}
is a $(\sL, \wGa)$-space and $\xi: B\cong B\times \{0\} \hra B\times V_{\bR} = A$ and $\pr_B : A= B\times V_{\bR}\to B$ are morphisms of $(\sL, \wGa)$-spaces. We call $A$ the {\it adelic $(\sL, \wGa)$-space}. We have
\begin{equation*}
\label{ahat} 
\hA =\,\hLa \times V_\bR 
\end{equation*}
Note that the morphism $\widehat{\pr_B}: \hA\to \hB$ induced by the projection $\pr_B: A \to B$ is given by $\hLa \times V_\bR\to \hLa, (\la, v)\mapsto \la$. 
\medskip

\noi (c) In example \eqref{covtorus} each space $A_L=B_L \times V_{\bR}$ is equipped with a natural $\Lambda$-action given by $(v + L, v_{\infty}) + \la := (v+ \la + L, v_{\infty} +\la)$ for $(v + L, v_{\infty})\in B_L \times V_{\bR}$ and $\la\in \Lambda$. We define the $(\sL, \wGa)$-space $T$ as  $T= A/\Lambda$. More precisely for $L\in \sL$ we let $T_L:= A_L/\Lambda \cong V_{\bR}/L$ and let $\pr_L : A_L\to T_L$ be the quotient map. Note that each is an $n$-dimensional real torus. 
It is clear that a morphism $\varphi: L_1\to L_2$ in $\cC(\sL, \wGa)$ induces a canonical map $\varphi: T_{L_1}\to T_{L_2}$ such that the diagram 
\begin{equation*}
\label{covtorus2}
\begin{CD}
B_{L_1} \times V_{\bR} @> \pr_{L_1} >> T_{L_1}\\
@VV \varphi V@VV \varphi V\\
B_{L_2} \times V_{\bR} @> \pr_{L_2} >> T_{L_2}
\end{CD}
\end{equation*}
commutes. Thus 
\begin{equation}
\label{torus}
T= T^{\sL}: \cC(\sL, \wGa)\lra \spaces, \quad L\mapsto T_L, \quad (\varphi: L_1\to L_2)\mapsto (\varphi: T_{L_1}\to T_{L_2})
\end{equation}
is a $(\sL, \wGa)$-space. It is in fact a $(\sL, \Ga)$-space, called the {\it torus $(\sL, \Ga)$-space}. Note that each $T_L \cong V_{\bR}/L$ is an $n$-dimensional real torus. For $\hT$ we get
\begin{equation*}
\label{torushat} 
 \hT  \, =\, \left(\hLa \times V_\bR\right)/\Lambda.
\end{equation*}
The collection of maps $\pr_L: A_L\to T_L$, $L\in \sL$ is a morphism of $(\sL, \Ga)$-spaces $\pr: A \to T$. The induced morphism 
$\hA\to \hT$ is the map
\begin{equation}
\label{univcov}
\pr: \hLa \times V_\bR\lra \left(\hLa \times V_\bR\right)/\Lambda, \quad x \mapsto x+ \Lambda.
\end{equation}
It is easy to see that the lattice topology on $\hT$ can also be defined as the quotient topology of the lattice topology on $\hA$ with respect to \eqref{univcov}. More precisely a subset $U\subseteq \hT$ is open if and only if $\pr^{-1}(U)$ is open in 
$\hA$ and a family of subsets $\{U_i\}_{i\in I}$ of an open subset $U\subseteq \hT$ is a covering if and only if $\{\pr^{-1}(U_i)\}_{i\in I}$ is a covering of $\pr^{-1}(U)$. 

Consider the composition of morphisms
\begin{equation}
\label{zerosec}
\io:=\pr\circ \xi: B\lra T
\end{equation} 
It has the following concrete description. For $L\in \sL$ we identify $T_L$ with $V_{\bR}/L$. Then the map $\io_L: \Lambda/L\lra V_{\bR}/L$ is induced by the inclusion $\Lambda\hra V_{\bR}$.
\medskip

\noi (d) Consider the case of Example \ref{examples:fracideals} (a) above, i.e.\ where $V=F$ is a number field and $\sL=\cI^S\subseteq \Lat(F)$ is the set of fractional ideals in $F$ that are coprime to a finite set of nonarchimedan places $S$ of $F$ and $\wGa=\cO_S^*\ltimes \cO_S= \Aff(\cO_S)$. For $\fa\in \cI^S$ we have
\[
B_{\fa} \, =\, \cO_S/\fa, \qquad A_{\fa}\,=\, (B\times F_{\infty})_{\fa} = \cO_S/\fa \times F_{\infty} \qquad \mbox{and}\qquad T_{\fa}\, = \, F_{\infty}/\fa
\]
where $F_{\infty}:= F\otimes \bR$. It will be useful to have adelic descriptions of $B$, $A$ and $T$. 
Recall that $\bA_f^S=\prod_{v\not\in S, v\nmid\infty}' F_v$ (resp.\ $\bA^S= \prod_{v\not\in S}' F_v$) are the finite prime-to-$S$ (resp.\ the prime-to-$S$) adeles of $F$. By the strong approximation theorem we have 
\begin{equation}
\label{btadel}
B_{\fa}\, =\, \bA_f^S/\hatfa^{S}, \qquad A_{\fa}\, =\, \bA^S/\hatfa^{S}\qquad \mbox{and}\qquad T_{\fa}\, =\, \bA^S/(\cO_S+ \hatfa^{S})
\end{equation}
where $\hatfa^{S} :=  (\prod_{v\not\in S, v\nmid\infty}\cO_v)\, \fa\subseteq \bA_f^S$. 

Using \eqref{btadel} we obtain
\begin{equation*}
\label{starhat}
\hB \,= \,\prolim_{\fa\in \cI^S} \cO_S/\fa \, =\, \prolim_{\fa\in \cI^S} \bA_f^S/\hatfa^{S}\, =\, \bA_f^S.
\end{equation*}
Note that for $\fa\in \cI^S$ we have 
\begin{equation*}
\label{fahat}
\hatfa^{\cI^S}:=\, \prolim_{\fa'\in \cI^S\, \fa'\subseteq \fa} \fa/\fa' \, =\, \hatfa^{S}.
\end{equation*} 
Thus a subset $U\subseteq \bA_f^S$ is open if and only if there exists $\fa\in \cI^S$ such that $U$ is a union of $\hatfa^{S}$-cosets. 

For $\hA$ and $\hT$ we get 
\begin{equation}
\label{torushat1} 
\hA =\,\bA^S \qquad \mbox{and}\qquad  \hT  \, =\,  \bA^S/\cO_S
\end{equation}
and $\pr: \hA\to \hT$ (resp.\ $\pr_B: \hA\to \hB$) is the natural projection $\bA^S \to \bA^S/\cO_S$ (resp.\ $\bA^S \to \bA_f^S$). We also remark that under the identifications \eqref{btadel}, \eqref{torushat1} the map $\pi_{\fa}: \hT\to T_{\fa}$ is given by the projection
\begin{equation*}
\label{torushat2}
\pi_{\fa}: \bA^S/\cO_S \lra \bA^S/(\cO_S+ \hatfa^{S}).
\end{equation*}
Note that the lattice topology on $\hB$, $\hA$ and $\hT$ are coarser than the usual topologies on $\bA_f^S$, $\bA^S$ and $\bA^S/\cO_S$ respectively. For example a subset $U\subseteq \bA^S$ is open in the lattice topology if and only if it is open in $\bA^S$ with respect to the usual topology and if it is $\hatfa^{S}$-stable for some $\fa\in \cI^S$.
\medskip

\noi (e) More generally we consider the framework of Example \ref{examples:fracideals} (b), i.e.\ we have $n=dm$, $F/\bQ$ is an extension of degree $m$, $V$ is an $d$-dimensional $F$-vector space, $\cM$ is a finitely generated $\cO_S$-submodule of $V$ with $\cM\otimes_{\cO_S} F=V$ and $\sL=\sL(\cM)$ is the set of finitely generated $\cO_F$-submodules $L$ of $M$ satisfying $L\otimes_{\cO_F} \cO_S = \cM$. In this case the associated $(\sL, \wGa)$-spaces $B$, $A$ and $T$ (for $\wGa=\Aff_{\cO_S}(\cM)$) admit the adelic descriptions 
\[ 
B_L \, =\, (\cM\otimes_{\cO_S} \bA_f^S)/\widehat{L}^S , \qquad A_L \, =\, (\cM\otimes_{\cO_S} \bA^S)/\widehat{L}^S
\qquad \mbox{and}\qquad T_L \, = \, (\cM\otimes_{\cO_S} \bA^S)/(\cM + \widehat{L}^S)
\]
for every $L\in \sL$ (where $\widehat{L}^{S}$ is the closure of $L$ in $\bA_f^S$, i.e.\ $\widehat{L}^{S} =  (\prod_{v\not\in S, v\nmid\infty}\cO_v) L$). Moreover we have
\begin{equation*}
\label{torushat3}
\hB \,= \, V \otimes_F \bA_f^S, \qquad \hA\, =\, V\otimes_F \bA^S \qquad \mbox{and}\qquad  \hT  \, =\, (V\otimes_F \bA^S)/\cM.
\end{equation*}
Again $\pr: V\otimes_F \bA^S \to (V\otimes_F \bA^S)/\cM$ and $\pr_B: V\otimes_F \bA^S\to V \otimes_F \bA_f^S$ are the natural projections.
\end{examples}

\paragraph{Sheaves on $(\sL, \wGa)$-spaces} Let $R$ be a ring and let $X$ be a $(\sL, \wGa)$-space. We denote the category of $R$-sheaves on the site $X=(\hX, \fOpen(X), \Cov)$ by $\Sh(X, R)$. The objects of $\Sh(X, R)$ will be called a $R$-sheaves on $X$. The notion of an $R$-presheaf on $X$ is defined similarly. By Prop.\ \ref{prop:sheafnice} the category $\Sh(X, R)$ is $R$-linear, abelian and has enough injectives. Also for $L\in \sL$ there is a pair of adjoint functors $\pi_L^*: \Sh(X_L, R)\to \Sh(X, R)$ and $(\pi_L)_*: \Sh(X, R)\to \Sh(X_L, R)$ associated to the projection $\pi_L: X\to X_L$. We will denote the second by 
\begin{equation}
\label{res1}
\Sh(X, R) \lra \Sh(X_L, R), \qquad \sF \mapsto \sF_L.
\end{equation}
Recall that for $\sF\in \Sh(X, R)$ the sheaf $\sF_L$ is given by $\sF_L(W) := \sF((\pi_L)^{-1}(W))$ for every open subset $W\subseteq X_L$. Note that for $L_1, L_2\in \sL$ with $L_1\subseteq L_2$ we have $\sF_{L_2} \cong (\pi_{L_1, L_2})_*\sF_{L_1}$. 

More generally if $\wGa'\subseteq\wGa$ is a subgroup then we define a $\wGa'$-equivariant $R$-sheaf $\sF$ on $X$ to be a $\wGa'$-equivariant $R$-sheaf on $\hX$ (see Def.\ \ref{df:equisheaf}). The category of $\wGa'$-equivariant $R$-sheaves on $X$ will be denoted by $\Sh(X, \wGa', R)$. 

We note that a sheaves on $X$ can be completely described in terms of collections of sheaves $\sF_L$ on $X_L$ for $L\in \sL$ together with the collections of isomorphisms $\sF_{L_2} \cong (\pi_{L_1, L_2})_*\sF_{L_1}$ for every pair $L_1\subseteq L_2$ in $\sL$. More precisely let $\Sh'(X, R)$ denote the category whose objects 
\begin{equation*}
\label{sheaf2b}
\sF' = \{\sF'_L, \rho_{L_1, L_2}; L, L_1, L_2\in \sL, L_1\subseteq L_2\}
\end{equation*}
consists of a collection of sheaves $\sF'_L\in \Sh(X_L, R)$ for $L\in \sL$ together with a collection of isomorphism $\rho_{L_1, L_2}: \sF_{L_2} \to (\pi_{L_1, L_2})_*\sF_{L_1}$ that satisfy the cocycle condition 
\begin{equation*}
\label{sheaf2}
 \left(\pi_{L_2, L_3}\right)_*\left(\rho_{L_1, L_2}\right)\circ \rho_{L_2, L_3} \,= \, \rho_{L_1, L_3}
\end{equation*}
for every triple of lattices $L_1 \subseteq L_2\subseteq L_3$ in $\sL$. A morphism $\alpha: \sF'\to \sG'$ in $\Sh'(X, R)$ between two objects $\sF'= \{\sF'_L\}_{L\in \sL}$ and $\sG'= \{\sG'_L\}_{L\in \sL}$
consists of a collection of morphisms $\alpha_L: \sF'_L \to \sG'_L$ in $\Sh(X_L, R)$ for every $L\in \sL$ that are compatible with the isomorphisms $\rho_{L_1, L_2}$ in the obvious sense. We have the obvious

\begin{lemma}
\label{lemma:sheafcomparison}
The functor 
\begin{equation*}
\label{sheaf3}
\Sh(X, R) \lra \Sh'(X, R),\quad \sF\mapsto \{\sF_L\}_{L\in \Lat}
\end{equation*}
is an equivalence of categories. 
\end{lemma}

As a first application we obtain

\begin{prop}
\label{prop:layer}
The functor \eqref{res1} is exact and preserves injectives. In fact a sequence of $R$-sheaves $\sF_1 \lra \sF_2 \lra \sF_3$ on $X$ is exact if and only if the sequence $\sF_{1, L} \lra \sF_{2, L} \lra \sF_{3, L}$ of $R$-sheaves on $X_L$ is exact for every $L\in \sL$.
\end{prop}

\begin{proof} This can be easily deduced from the exactness of the functor $(\pi_{L_1, L_2})_*:  \Sh(X_{L_1}, R)\to  \Sh(X_{L_2}, R)$ for every pair of lattices $L_1\subseteq L_2$ in $\sL$.
\end{proof}

Let $f: X\to Y$ be a morphism of $(\sL, \wGa)$-spaces. As mentioned before $f$ induces a $\wGa$-equivariant continuous morphism of sites $\hf: \hX\to \hY$. We will denoted the functor $\hf_*$ and its left adjoint $\hf^*$ (see \eqref{fimage} and \eqref{finv})
by 
\begin{equation*}
\label{pushout}
f_*: \Sh(X, R) \lra \Sh(Y, R) \qquad \mbox{and} \qquad f^*: \Sh(Y, R) \lra \Sh(X, R).
\end{equation*}
Thus $f_*$ is given $f_*(\sF)(V) = \sF(\hf^{-1}(V))$ for every open subset $V\subseteq \hY$. 

\begin{prop}
\label{prop:pullbackex}
The functor $f^*: \Sh(Y, R) \to \Sh(X, R)$ is exact. Moreover for $\sG\in \Sh(Y, R)$ and $L\in \sL$ we have
\begin{equation}
\label{pullback2}
f^*(\sG)_L\,=\, (f_L)^*(\sG_L).
\end{equation}
Also for the right derived functors of $f_*$ we have 
\begin{equation}
\label{derivedfstar}
(R^i f_*(\sF))_L\,=\, (R^i f_L)_*(\sF_L) 
\end{equation}
for every $i\ge 0$ and $\sF\in \Sh(X, R)$.
\end{prop}

\begin{proof} The first assertion and \eqref{pullback2} can be seen by interpreting the functor $f_*$ in terms of the categories $\Sh'(X, R)$ and $\Sh'(Y, R)$ (compare Lemma \ref{lemma:sheafcomparison}). Namely, by the definitions we have $(f_*(\sF))_L=(f_L)_*(\sF_L)$ for every $L\in \sL$ and $\sF\in \Sh(X, R)$. Thus the functor $f_*: \Sh(X, R) \to \Sh(Y,R)$ corresponds under the equivalences $\Sh(X, R)\simeq \Sh'(X, R)$, $\Sh(Y, R)\simeq \Sh'(Y, R)$ to 
\begin{equation}
\label{sheaf11a}
f_*: \Sh'(X, R) \lra \Sh'(Y, R),\quad \{\sF'_L, \rho_{L_1, L_2}\} \mapsto  \{(f_L)_*(\sF'_L), (f_{L_2})_*(\rho_{L_1, L_2})\}
\end{equation}
Since the maps $\pi_{L_1, L_2}$ are coverings we have $(\pi_{L_1, L_2})_* \circ (f_{L_1})^* \cong (f_{L_2})^*\circ (\pi_{L_1, L_2})_*$
for every pair $L_1\subseteq L_2$ in $\sL$. From this it follows easily that the functor 
\begin{equation*}
\label{sheaf12}
f^*: \Sh'(Y, R) \lra \Sh'(X, R),\quad \{\sG'_L, \rho_{L_1, L_2}\} \mapsto  \{(f_L)^*(\sG'_L), (f_{L_2})^*(\rho_{L_1, L_2})\}
\end{equation*}
is well-defined, exact and left adjoint to \eqref{sheaf11a}. 

For \eqref{derivedfstar} note that $\pi_L \circ f = f_L \circ \pi_L$ hence $(\pi_L)_* \circ R^i f_* = R^i (\pi_L \circ f)_* = R^i(f_L \circ \pi_L)_* = R^i(f_L)_* \circ (\pi_L)_*$ by Prop.\ \ref{prop:layer}.
\end{proof}

\begin{example}
\label{example:torsh} \rm Let $\pr: A\to T$ be the natural morphism between the $(\sL, \wGa)$-spaces defined in Examples \ref{examples:torpt} (b), (c), i.e.\ we have $\pr_L: A_L =\Lambda/L\times V_\bR\to T_L = \left(\Lambda/L\times V_\bR\right)/\Lambda \approx V_{\bR}/L$ is the natural projection for every $L\in \sL$ with $\Lambda =\bigcup_{L\in \sL} L$. Since the lattice topology on $\hT$ can also be defined as the quotient topology of the lattice topology on $\hA$ with respect to \eqref{univcov} we see that the functor \begin{equation*} 
\label{unicov2} 
\Sh(T, R) \lra \Sh(A, \Lambda, R), \quad \sF \mapsto \pr^*(\sF)
\end{equation*}
is an equivalence of categories. 
\end{example}

\begin{df}
\label{df:latspace}
A $(\sL, \wGa)$-space $X$ will be called discrete if $X_L$ is a discrete topological space for every $L\in \sL$. 
\end{df}

In the case of a discrete $(\sL, \wGa)$-space it is easy to characterize when an $R$-presheaf $\sF$ on $X$ is a sheaf.

\begin{lemma}
\label{lemma:sheafdisc} 
Let $X$ be a discrete $(\sL, \wGa)$-space and let $\sF$ be an $R$-presheaf $\sF$ on $\hX$.
\medskip

\noi (a) $\sF$ is a sheaf if and only if the map 
\begin{equation*} 
\label{sheafdisc}
\begin{CD}
\sF(U) @> s\mapsto  (s|_{U_i})_{i\in I} >> \prod_{i\in I} \cF(U_i)
\end{CD}
\end{equation*}
is an isomorphism for every $L\in \sL$, every $L$-stable subset $U\subseteq \hX$ and every covering $U=\bigcup_{i\in I} U_i$ of $U$ by disjoint $L$-stable subsets.
\medskip

\noi (b) Assume that $\sF$ is a sheaf and let $U_1\subseteq U_2 \subseteq \hX$ be open subsets. Then the restriction $\res_{U_2, U_1}: \sF(U_2) \to \sF(U_1)$ has a canonical section 
\begin{equation*} 
\label{ressection}
\io_{U_1, U_2}: \sF(U_1) \, \lra\,  \sF(U_2)
\end{equation*}
i.e.\ we have $\res_{U_2, U_1}\circ \io_{U_1, U_2} =\id_{\sF(U_1)}$. 
\end{lemma}

\begin{proof} (b) Since $X$ is discrete the complement $U_2\setminus U_1$ is again open and $U_2=U_1 \cup (U_2\setminus U_1)$ is an open covering. Hence 
\begin{equation*} 
\label{sheafdisc2}
\begin{CD}
\sF(U_2) @> s\mapsto  (s|_{U_1}, s|_{U_2\setminus U_1}) >> \sF(U_1) \oplus \sF(U_2\setminus U_1)
\end{CD}
\end{equation*}
is an isomorphism.
\end{proof}

\paragraph{Stalks} Let $X$ be a $(\sL, \wGa)$-space and let $\sF$ be an $R$-sheaf on $X$. For a point $x\in \hX$ we introduce two types of $R$-modules $\sF_x$ and $\sF^x$ which may be viewed both as certain kind of stalk of $\sF$ in $x$. For $L\in \sL$ put $x_L:=\pi_L(x)$. Firstly, we consider the case when $X$ is discrete. Then for $L_1, L_2\in \sL$ with $L_1\subseteq L_2$ the fibers $\pi_{L_1}^{-1}(x_{L_1})\subseteq \pi_{L_2}^{-1}(x_{L_2})$ are open subsets of $\hX$ containing $x$.  By part (b) of the Lemma the restriction map 
\begin{equation} 
\label{resstalk}
p_{x, L_2, L_1}:= \res_{\pi_{L_2}^{-1}(x_{L_2}), \pi_{L_1}^{-1}(x_{L_1})}: \sF(\pi_{L_2}^{-1}(x_{L_2})) \lra \sF(\pi_{L_1}^{-1}(x_{L_1})), \quad s\mapsto s|_{\pi_{L_1}^{-1}(x_{L_1})}
\end{equation}
has a canonical section
\begin{equation} 
\label{ressection2}
\io_{x, L_1, L_2}: \sF(\pi_{L_1}^{-1}(x_{L_1}))\lra \sF(\pi_{L_2}^{-1}(x_{L_2})).
\end{equation}
The collection of $R$-modules $\{\sF(\pi_L^{-1}(x_L))\}_{L\in \sL}$ together with the maps \eqref{ressection2} form a direct system of $R$-modules over the ordered set $\sL= (\sL, \subseteq)$ so we can consider the direct limit
\begin{equation*}
\label{stalklow}
\sF_x := \, \dlim_{L\in \sL} \sF(\pi_L^{-1}(x_L)). 
\end{equation*}
Similarly the collection of $R$-modules $\{\sF(\pi_L^{-1}(x_L))\}_{L\in \sL}$ together with the maps \eqref{resstalk} form an inverse system over $\sL^{\opp}$ and we can define 
\begin{equation*}
\label{stalkup}
\sF^x := \, \prolim_{L\in \sL^{\opp}} \sF(\pi_L^{-1}(x_L)). 
\end{equation*}
For an arbitrary $(\sL, \wGa)$-space $X$ we define $\sF_x$ and $\sF^x$ by first pulling $\sF$ back to the associated discrete $(\sL, \wGa)$-space $X_{\disc}$ (i.e.\ $(X_{\disc})_L$ is the set $X_L$ equipped with the discrete topology for every $L\in \sL$). More precisely we define

\begin{df}
\label{df:stalk}
Let $\sF$ be an $R$-sheaf on a $(\sL, \wGa)$-space $X$. For a point $x\in \hX$ we define the $R$-modules $\sF_x$ and $\sF^x$ by 
\begin{equation}
\label{stalk}
\sF_x := \io^*(\sF)_x \qquad \mbox{and} \qquad \sF^x := \io^*(\sF)^x
\end{equation}
where $\io: X_{\disc}\to X$ is the identical map viewed as a morphism of $(\sL, \wGa)$-spaces. We call $\sF^x$ the upper and $\sF_x$ the lower stalk of $\sF$ at $x$.
\end{df}

\begin{example}
\label{example:stalkb} \rm Let $B$ be a discrete $(\sL, \wGa)$-space defined in Example \ref{examples:torpt} (a) and let $\sF$ be an $R$-sheaf on $B$. We consider the upper stalk of $\sF$ at $0\in \hB=\hLa$. For that put $0_L: =\pi_L(0) = L\in \Lambda/L=B_L$ so that $\pi_L^{-1}(0_L)=\hL$ for every $L\in \sL$. Note that we have $\hB= \bigcup_{L\in \sL} \pi_L^{-1}(0_L)$. Moreover for a fixed lattice $L_0\in \sL$ every element of the family $\cU= \{\pi_L^{-1}(0_L)\}_{L\in \sL, L\supseteq L_0}$ is $\hL_0$-stable hence $\cU$ is an open covering of $\hB$. Since $\{ L\in \sL\mid L\supseteq L_0\}^{\opp}$ is a filtered partially ordered set the sheaf property implies the last equality in
\begin{equation*}
\label{stalkb0}
\sF^0 \,= \, \prolim_{L\in \sL^{\opp}} \sF(\pi_L^{-1}(x_L))\, =\, \prolim_{\{ L\in \sL\mid L\supseteq L_0\}^{\opp}} \sF(\pi_L^{-1}(0_L))\, =\, \sF(\hB).
\end{equation*}
\end{example}

\begin{remarks}
\label{remarks:stalks}
\rm (a) \rm If $\sF$ is a $\wGa$-equivariant sheaf then the $\wGa$-actions on $\sF$ and $\hX$ induce canonical homomorphisms 
\begin{equation*}
\label{stalk5}
\wga: \sF_x \lra \sF_{\wga^{-1}(x)}, \qquad \wga: \sF^x \lra \sF^{\wga^{-1}(x)}.
\end{equation*}
for every $\wga\in \wGa$. In particular if $x$ is a $\wGa$ fixed-point then $\sF_x$ and $\sF^x$ are $R[\wGa]$-modules. A similar remark holds for $\Ga$-equivariant sheaves. 
\medskip

\noi (b) One could also consider the stalk of $\sF\in \Sh(X, R)$ at $x\in \hX$ defined in the usual sense, namely as the direct limit $\dlim_{U\in \cU_x} \sF(U)$ where $\cU_x$ consists of all open subsets of $\hX$ that contain $x$ and where the transition maps are the restrictions. The above two $R$-modules \eqref{stalk} differ in general from this "naiv" stalk. For that assume that $X$ is discrete so that the fibers $\{\pi_L^{-1}(x_L)\}_{L\in \sL}$ form a cofinal subset of $\cU_x$. Thus we have 
$\dlim_{U\in \cU_x} \sF(U)= \dlim_{L\in \sL^{\opp}} \sF(\pi_L^{-1}(x_L))$ where the transition maps are the maps \eqref{resstalk}. 
\end{remarks}
 
\begin{lemma}
\label{lemma:stalklim}
Let $X$ be a $(\sL, \wGa)$-space and let $\sF$ be an $R$-sheaf on $X$. Let $x\in \hX$ and put $x_L=\pi_L(x)$ for $L\in \sL$. Then we have 
\begin{equation}
\label{stalk3}
\sF_x \,= \, \dlim_{L\in \sL} (\sF_L)_{x_L} \qquad \mbox{and} \qquad \sF^x \,= \,\prolim_{L\in \sL^{\opp}} (\sF_L)_{x_L}.
\end{equation}
\end{lemma}

\begin{proof} Assume first that $X$ is discrete. Then the $R$-module $\sF(\pi_L^{-1}(x_L))=\sF_L(\{x_L\})$ is equal to stalk of the sheaf $\sF_L$ on $X_L$ at the point $x_L$ for every $L\in \sL$ (since $X_L$ is discrete). Therefore \eqref{stalk3} holds in this case. If $X$ is not discrete then by applying \eqref{pullback2} to the morphism $\io: X_{\disc}\to X$ we obtain
\begin{equation*} 
\label{stalk3b}
\sF_x \,= \, \io^*(\sF)_x \, =\, \dlim_{L\in \sL} (\io^*(\sF)_L)_{x_L}  \, =\, \dlim_{L\in \sL} (\io_L)^*(\sF_L)_{x_L}\, =\, \dlim_{L\in \sL} (\sF_L)_{x_L}
\end{equation*}
The proof of the second equality in \eqref{stalk3} is analogous.
\end{proof}

An immediate consequence of Prop.\ \ref{prop:pullbackex} and Lemma.\ \ref{lemma:stalklim} is

\begin{prop}
\label{prop:pullbackstalk}
Let $f:X\to Y$ be a morphism of $(\sL, \wGa)$-spaces and let $\sF\in \Sh(Y, R)$. We have 
\begin{equation*}
\label{stalk4}
f^*(\sF)_x \, =\,  \sF_{\hf(x)}\qquad \mbox{and} \qquad f^*(\sF)^x \, =\,  \sF^{\hf(x)}
\end{equation*}
for every $x\in \hX$. 
\end{prop}

\paragraph{$\wGa$-stable closed subspaces} We introduce the notion of a $\wGa$-stable closed subspace of a $(\sL, \wGa)$-space.
The main example we have in mind are a finite set of torsion points of the torus $(\sL, \wGa)$-space.

\begin{df}
\label{df:closedsubspace}
Let $X$ be a $(\sL, \wGa)$-space. A $\wGa$-stable closed subspace of $X$ is a pair $(C, \io)$ consisting of a topological space $C$ equipped with a $\wGa$-action and a $\wGa$-equivariant morphism $\io: C\to \hX$ of sites such that $\io_L:= \pi_L\circ \iota: C\to X_L$ is a closed embedding for every $L\in \sL$, i.e.\ we have 
 \medskip
 
\noi (i) $\io_L: C\to X_L$ is injective and $C_L:= \io_L(C)$ is closed in $X_L$.
\medskip

\noi (ii) The induced map $\io_L: C\to C_L$ is a homeomorphism.
\end{df}

Note that the fact that $\io$ is $\wGa$-equivariant implies that we have $\wga(C_L) = C_{\gamma(L)}$ for every $L\in \sL$ and $\wga=(\gamma, \la)\in \wGa$. Note also that condition (i) implies (ii) if $C$ is finite. If $X$ be a $(\sL, \wGa)$-space and $(C, \io)$ is a $\wGa$-stable closed subspace of $X$ then the $\wGa$-action on $X$ factors through $\Ga$ (i.e.\ $\La$-acts trivially on $C$). In this case we call $(C, \io)$ is a $\Ga$-stable closed subspace of $X$. 

\begin{examples}
\label{examples:closedsub} 
\rm Let $F$ be a number field, $S$ a finite set of nonarchimedean places of $F$ and $\bA^S$ (resp.\ $\bA_f^S$) the prime-to-$S$ (resp.\ finite prime-to-$S$) adeles of $F$. We give examples of $\wGa$-stable closed subspaces of the $(\sL, \wGa)$-spaces $A$ and $T$ considered in Examples \ref{examples:torpt} (d) and (e). 
\medskip

\noi (a) Let $V=F$ and $\sL = \cI^S$ be as in Examples \ref{examples:torpt} (d). We fix a subgroup $\Ga$ of $\cO_S^*$ and put $\wGa = \Ga\ltimes \cO_S\subseteq \Aff(\cO_S)$. Firstly, we describe certain $\Ga$-stable closed subspaces of the $(\sL, \Ga)$-space $T$. For that note that for $\fa\in \cI^S$ the restriction of the projection $\pi_{\fa}: \hT= \bA^S/\cO_S\to \bA^S/(\cO_S+\hatfa^{S}) \cong F_{\infty}/\fa$ to the subset $F/\cO_S\subseteq \bA^S/\cO_S$ is injective. Indeed, 
if we put $S^{-1}\fa := (S^{-1}\cO_F)\cdot \fa$ then we have $F=\cO_S+S^{-1}\fa$ and $\cO_S\cap S^{-1}\fa =\fa$ hence $S^{-1} \fa/\fa \cong F/\cO_S$. Under this identification the restriction of $\pi_{\fa}$ to $F/\cO_S$ is given by the inclusion $S^{-1} \fa/\fa\hra F/\fa\hra F_{\infty}/\fa$. Thus if $C$ is a finite $\Ga$-stable subset of $V/\cO_S\subseteq \hT = \bA^S/\cO_S$ then the pair $(C, \io)$, where $\io=\incl: C\hra  V/\cO_S\hra \bA^S/\cO_S$ is the inclusion, is a $\Ga$-stable closed subspaces of $T$. For example if $\fA$ is an ideal of $\cO_S$ and if we put $T[\fA]:= \fA^{-1}/\cO_S$ then the pair $(T[\fA], \io)$ a $\cO_S^*$-stable closed subspaces of $T$. 
\medskip

\noi (b) Let $F$, $S$, $V$, $\cM$ and $\sL$ be as in Examples \ref{examples:fracideals} (b) and \ref{examples:torpt} (e).
Again it is easy to see that the restriction of the projection $\pi_L : \hT=\cM\otimes_{\cO_S} \bA^S/\cM\to T_L = (\cM\otimes_{\cO_S} \bA^S)/(\cM + \widehat{L}^S)$ 
to the subset $V/\cM = \cM\otimes_{\cO_S} F/\cM \subseteq \cM\otimes_{\cO_S} \bA^S/\cM$ is injective. Hence if $\Ga$ is a subgroup of $\GL_{\cO_S}(\cM)$ and $C \subseteq V/\cM$ a finite $\Ga$-stable subset then the pair $(C, \incl)$ is a $\Ga$-stable closed subspace of $T$. For example if $\cN\supseteq \cM$ is another finitely generated $\cO_S$-submodule of $V$ and if $\Ga$ is a subgroup of $\GL_F(V)$ leaving both $\cM$ and $\cN$ invariant then the pair $(\cN/\cM, \io)$ is a $\wGa$-stable closed subspace of $T$. 
\end{examples}

\paragraph{Sheaf Cohomology} Let $X$ be a $(\sL, \wGa)$-space and let $R$ be a ring. For $\sF\in \Sh(X, R)$ (resp.\ $\sF\in \Sh(X, \wGa, R)$) we will write $H^i(X, \sF)$ (resp.\ $H^i(X, \wGa, \sF)$) for the cohomology groups $H^i(\hX, \sF)$ (resp.\ $H^i(\hX, \wGa, \sF)$). Also if $X$ is a $(\sL, \Ga)$-space and if $\sF\in \Sh(X, \Ga, R)$ then the cohomology groups $H^i(\hX, \Ga, \sF)$ will be denoted by $H^i(X, \Ga, \sF)$. We have

\begin{prop}
\label{prop:layercoho}
For every $L\in \sL$ and $\sF\in \Sh(X, R)$ there exists a canonical isomorphism 
\[
H^i(X, \sF)\, \cong \, H^i(X_L, \sF_L).
\]
\end{prop}

\begin{proof} This follows immediately from Prop.\ \ref{prop:layer} and the fact that we have $\sF(\hX) = \sF_L(X_L)$ for every $\sF\in \Sh(X, R)$.
\end{proof}

In particular we obtain 

\begin{coro}
\label{coro:cohodisc}
Let $X$ be a discrete $(\sL, \wGa)$-space. We have 
\medskip

\noi (a) $H^i(X, \sF)=0$ for every $\sF\in \Sh(X, R)$ and $i\ge 1$.
\medskip

\noi (b) $H^i(X, \wGa, \sF)\, =\, H^i(\wGa, H^0(X, \sF))$ for every $\sF\in \Sh(X, \wGa, R)$ and $i\ge 0$.
\end{coro}

\begin{proof} (a) Let $L\in \sL$. By Prop.\ \ref{prop:layercoho} we have 
$H^i(X, \sF) = H^i(X_L, \sF_{L}) = 0$
for every $i\ge 1$ since $X_L$ is discrete. 

For (b) note that (a) implies that the spectral sequence \eqref{hsscomp} degenerates. 
\end{proof}

Assume now that there exists $L\in \sL$ that is stabilized by the action of $\Ga$ so that the topological space $X_L$ becomes equipped with a $\wGa=\Ga\ltimes \Lambda$-action and that $\sF_L$ for $\sF\in \Sh(X, \wGa, R)$ becomes a $\wGa$-equivariant sheaf on $X_L$. In this case we can compare the cohomology groups $H^{\bu}(X, \wGa, \sF)$ with the equivariant cohomology groups $H^{\bu}(X_L, \wGa, \sF_L)$ in the usual sense. 

\begin{coro}
\label{coro:classicequiv}
Let $X$ be a $(\sL, \wGa)$-space. Assume that $L\in \sL$ is $\Ga$-stable, i.e.\ that we have $\bga(L) = L$ for every $\bga\in \Ga$. 
Then there are canonical isomorphisms 
\begin{equation}
\label{compequivcoh}
H^{\bu}(X, \wGa, \sF)\, \cong\, H^{\bu}(X_L, \wGa, \sF_L).
\end{equation}
for every $\sF\in \Sh(X, \wGa, R)$. A similar statement holds for $\Ga$-equivariant cohomology if 
$X$ is a $(\sL, \Ga)$-space and $\sF$ is a $\Ga$-equivariant sheaf on $X$.
\end{coro}

\begin{remark}
\label{remark:slfinal}
\rm Let $X$ be a $(\sL, \Ga)$-space and assume that $\sL$ contains a maximal element $L_0$ (so that $\Lambda =L_0$ and $\Ga\subseteq \GL(L_0)$). Then Prop.\ \ref{prop:layercoho} and Cor.\ \ref{coro:classicequiv} imply
\[
H^{\bu}(X, \sF)\, =\, H^{\bu}(X_{L_0}, \sF_{L_0}) \quad  \mbox{and}\quad   H^{\bu}(X, \Ga, \sF)\, =\, H^{\bu}(X_{L_0}, \Ga, \sF_{L_0})
\]
for $\sF\in \Sh(X, R)$ and $\sF\in \Sh(X, \Ga, R)$ respectively.
\end{remark}

\begin{proof} We show that there exists a morphism of spectral sequences for $\wGa$-equivariant cohomology (see Prop.\ \ref{prop:covss})
\begin{eqnarray}
\label{morhsscomp}
&& \left(E_2^{rs} = H^r(\wGa, H^s(X, \sF)) \, \Longrightarrow\, E^{r+s}=H^{r+s}(X, \wGa, \sF)\right) \, \lra \, \\
&& \hspace{3cm}  \left(E_2^{rs} = H^r(\wGa, H^s(X_L, \sF_L)) \, \Longrightarrow\, E^{r+s}=H^{r+s}(X_L, \wGa, \sF_L)\right) .
\nonumber
\end{eqnarray}
By Prop.\ \ref{prop:layercoho} it is an isomorphism on the $E_2$-page hence also on the limit terms. 

To define \eqref{morhsscomp} we argue as in (\cite{spiess2}, Prop.\ 3.42). Let $0\to \sF \to \sI^{\bu}$ be an injective resolution of $\sF$ in $\Sh(X, \wGa, R)$. By Prop.\ \ref{prop:injequiv} (b) the sequence $0\to \sF_L \to \sI_L^{\bu}$ is still exact. Let $0\to \sF_L \to \sJ^{\bu}$ be an injective resolution of $\sF_L$ in the category $\Sh(X_L, \wGa, R)$ and let 
$\alpha:  \sI_L^{\bu}\to \sJ^{\bu}$ be a morphism such that 
\[
\begin{CD}
0 @>>> \sF_L @>>> \sI_L^{\bu}\\
@.@VV \id V @VV \alpha V\\
0 @>>> \sF_L @>>> \sJ_L^{\bu}
\end{CD}
\]
commutes. Passing to global sections yields a homomorphism of complexes of $R[\wGa]$-modules 
\begin{equation*}
\label{analytic6f}
\begin{CD} 
\sI^{\bu}(X) \,=\, \sI_L^{\bu}(X_L) @> \alpha >> \sJ^{\bu}(X_L)
\end{CD} 
\end{equation*}
that induces a morphism between $\wGa$-hypercohomology spectral sequences
\begin{eqnarray*}
&&  \left(\,E_2^{rs} = H^r(\wGa, H^s(\sI^{\bu}(X)))\, \Longrightarrow\, E^{r+s}=H^{r+s}(\wGa, \sI^{\bu}(X)) \,\right) \, \lra \, \\
&& \hspace{3cm} \left( \,E_2^{rs} = H^r(\wGa, H^s(\sJ^{\bu}(X_L)))\, \Longrightarrow\, E^{r+s}=H^{r+s}(\wGa, \sJ^{\bu}(X_L))
\, \right).
\end{eqnarray*}
The first spectral sequence can be identified with the source (by Prop.\ \ref{prop:injequiv} (d)) and the second with the target of \eqref{morhsscomp}.
\end{proof}

We can also define cohomology with support in a $\wGa$-stable closed subspace $(C, \io)$. As before for $L\in \sL$ we put $\io_L= \pi_L\circ \io: C\hra X_L$ and set $C_L:= \io_L(C)$. For $\sF\in \Sh(X, R)$ we set
\begin{equation*}
\label{suppc}
H_{\brC}^0(X, \sF) \, =\, \{ s\in \sF(\hX)\mid \exists L\in \sL, s|_{\hX\setminus \pi_L^{-1}(C_L)} =0\}.
\end{equation*}
For $L_1, L_2\in \sL$ with $L_1\subseteq L_2$ we have $\hX\setminus \pi_{L_2}^{-1}(C_{L_2}) \subseteq \hX\setminus \pi_{L_1}^{-1}(C_{L_1})$ hence
\begin{eqnarray}
\label{suppc1}
&& H^0_{C_{L_1}}(X_{L_1}, \sF_{L_1}) =\ker(\res: \sF(\hX) \lra \sF(\hX\setminus \pi_{L_1}^{-1}(C_{L_1})))\\
&& \hspace{4cm} \subseteq H^0_{C_{L_2}}(X_{L_2}, \sF_{L_2})=\ker(\res: \sF(\hX) \lra \sF(\hX\setminus \pi_{L_2}^{-1}(C_{L_2}))).
\nonumber
\end{eqnarray}
It follows
\begin{equation}
\label{suppc2}
H_{\brC}^0(X, \sF) \, =\, \bigcup_{L\in \sL} \ker(\res: \sF(\hX) \lra \sF(\hX\setminus \pi_L^{-1}(C_L))) \, =\, \dlim_{L\in \sL} H_{C_L}^0(X_L, \sF_L).
\end{equation}
where we view $\sL$ again as a partially ordered set with respect to the inclusion.  We also define 
\begin{equation}
\label{suppc3}
H^0(X\setminus \brC, \sF) : =\, \dlim_{L\in \sL} \sF(\hX\setminus \pi_L^{-1}(C_L)) \, =\, \dlim_{L\in \sL} H^0(X_L\setminus C_L, \sF_L)
\end{equation}
where the transition maps in the direct limit are restriction maps. We note that if $\sF$ is a $\wGa$-equivariant $R$-sheaf then $H_{\brC}^0(X, \sF)$ and $H^0(X\setminus \brC, \sF)$ carry natural $\wGa$-actions, so we may consider their fixmodules
\begin{equation}
\label{equivsupp1}
H_{\brC}^0(X, \wGa, \sF) := \,H_{\brC}^0(X, \sF)^{\wGa},\qquad  H^0(X\setminus \brC, \wGa, \sF) := \,H^0(X\setminus \brC, \sF)^{\wGa}.
\end{equation}

\begin{df}
\label{df:cohsupp}
(a) The $i$-th right derived functor of the functor $H_{\brC}^0(X, \wcdot)$ will be denote by
\begin{equation*}
\label{suppc4}
H_{\brC}^i(X, \wcdot): \Sh(X, R) \lra \Mod_R, \quad \sF\mapsto H_{\brC}^i(X, \sF)
\end{equation*}
and the $i$-th right derived functor of $H^0(X\setminus \brC, \wcdot)$ by 
\begin{equation*}
\label{suppc5}
H^i(X\setminus \brC, \wcdot): \Sh(X, R) \lra \Mod_R, \quad \sF\mapsto H^i(X\setminus \brC, \sF).
\end{equation*}
(b) The derived functors of the two functors \eqref{equivsupp1} will be denoted by
\begin{eqnarray*}
\label{equivsupp2}
&& H_{\brC}^i(X, \wGa, \wcdot): \Sh(X, \wGa, R) \lra \Mod_R, \quad \sF\mapsto H_{\brC}^i(X, \wGa, \sF)\hspace{2cm} \mbox{and}\\
\label{equivsupp3}
&& H^i(X\setminus \brC, \wGa, \wcdot): \Sh(X, \wGa, R) \lra \Mod_R, \quad \sF\mapsto H^i(X\setminus \brC, \wGa, \sF)
\end{eqnarray*}
respectively.
\end{df}

\begin{prop}
\label{prop:cohsupp1}
Let $X$ be a $(\sL, \wGa)$-space and let $(C, \io)$ be a $\wGa$-stable closed subspace of $X$.
\medskip

\noi (a) We have 
\begin{eqnarray*}
H_{\brC}^i(X, \sF) & = & \dlim_{L\in \sL} H_{C_L}^i(X_L, \sF_L)\\
H^i(X\setminus \brC, \sF) & = & \dlim_{L\in \sL} H^i(X_L\setminus C_L, \sF_L)
\end{eqnarray*}
for every $i\in \bZ_{\ge 0}$ and $\sF\in \Sh(X, R)$. Moreover if $\sF$ is a $\wGa$-equivariant sheaf then the $R$-modules $H_{\brC}^i(X, \sF)$ and $H^i(X\setminus \brC, \sF)$ carry a natural $\wGa$-action. 
\medskip

\noi (b) There exists a long exact sequence 
\begin{equation}
\label{supp6}
\ldots\lra  H_{\brC}^i(X, \sF) \lra H^i(X, \sF) \lra H^i(X\setminus \brC, \sF) \lra H_{\brC}^{i+1}(X, \sF)\lra \ldots 
\end{equation}
for every $\sF\in \Sh(X, R)$.
\end{prop}

\begin{proof} (a) follows immediately from \eqref{suppc2}, \eqref{suppc3}, Prop.\ \ref{prop:layer} (b), (c) and the exactness of the direct limit. (b) Follows from (a) by passing in the long exact sequence 
\begin{equation*}
\label{supp7}
\ldots\lra  H_{C_L}^i(X_L, \sF_L) \lra H^i(X_L, \sF_L) \lra H^i(X_L\setminus C_L, \sF_L) \lra H_{C_L}^{i+1}(X_L, \sF_L)\lra \ldots 
\end{equation*}
to the direct limit over $L\in \sL$. 
\end{proof}

\begin{remark}
\label{remark:layercohsupp}
\rm Note that Prop.\ \ref{prop:cohsupp1} implies that there are canonical homomorphisms
\begin{equation}
\label{supp5}
H_{C_L}^i(X_L, \sF_L)\, \lra \, H_{\brC}^i(X, \sF), \hspace{1cm}   H^i(X_L\setminus C_L, \sF_L)\, \lra \, H^i(X\setminus \brC, \sF)
\end{equation}
for every $L\in \sL$, $\sF\in \Sh(X, R)$ and $i\ge 0$. If $L$ is $\wGa$-stable lattice and $\sF$ is a $\wGa$-equivariant sheaf then it is easy to see that the maps \eqref{supp5} are $\wGa$-equivariant.
\end{remark} 

\begin{coro}
\label{coro:gpt}
Let $X$ be a $(\sL, \wGa)$-space and let $(C, \io)$ be a $\wGa$-stable closed and discrete subspace of $X$.

\noi (a) Let $\sF\in \Sh(X, R)$ and assume that the following conditions hold
\medskip

\noi (i) $X_L$ is an oriented $n$-dimensional manifold for every $L\in \sL$.\\
\noi (ii) The covering $\pi_{L_1, L_2}: X_{L_1}\subseteq X_{L_2}$ preserves the orientations for every pair $L_1\subseteq L_2$ in $\sL$.\\
\noi (iii) $\sF_L$ is a locally constant sheaf.
\medskip

\noi Then we have
\begin{equation}
\label{supp8}
H_{\brC}^i(X, \sF) \, \cong \, \left\{\begin{array}{cc} \bigoplus_{c\in C} \sF_c & \mbox{if $i=n$,}\\
0 & \mbox{if $i\ne n$.}
\end{array}
\right.
\end{equation}
\noi (b) Let $\sF\in \Sh(X, \wGa, R)$. Assume that (i), (iii) and 
\medskip

\noi (ii') The covering $\varphi: X_{L_1} \to X_{L_2}$ preserves orientations for every morphism $\varphi: L_1 \to L_2$ in $\cC(\sL, \wGa)$ 
\medskip

\noi holds. Then \eqref{supp8} is a $\wGa$-equivariant isomorphism.
\end{coro}

\begin{proof} By (\cite{ks2}, 3.2.3) we have
\begin{eqnarray*}
H_{\brC}^i(X, \sF) & = & \dlim_{L\in \sL} H_{C_L}^i(X_L, \sF_L)
\, = \, \left\{\begin{array}{cc} \dlim_{L\in \sL} (\sF_L)|_{C_L} & \mbox{if $i=n$,}\\
0 & \mbox{if $i\ne n$.}
\end{array}
\right.\\
& = & \left\{\begin{array}{cc} \bigoplus_{c\in C} \sF_c & \mbox{if $i=n$,}\\
0 & \mbox{if $i\ne n$.}
\end{array}
\right.
\end{eqnarray*}
\end{proof}

In the examples we are interested in (namely the examples \eqref{covtorus} and \eqref{torus}) the strong condition (ii') usually does not hold. Instead of (ii') we assume that there exists a homomorphism $\ep: \Gamma\to \{\pm 1\}$ that indicates whether the covering $\varphi: X_{L_1} \to X_{L_2}$ preserves the orientation or not. 

\begin{coro}
\label{coro:gpt2}
Let $X$ be a $(\sL, \wGa)$-space, let $(C, \io)$ be a $\wGa$-stable closed and discrete subspace of $X$ and let $\ep: \Gamma\to \{\pm 1\}$ be a homomorphism. Let $\sF\in \Sh(X, \wGa, R)$ and assume that the conditions (i), (iii) of Cor.\ \ref{coro:gpt} as well as 
\medskip

\noi (ii'') The covering $\varphi: X_{L_1} \to X_{\gamma(L_2)}$ preserves (resp.\ reverses) orientations for every morphism $\varphi = (\wga, L_1, L_2)$, $\wga=(\gamma, \la)$ in $\cC(\sL, \wGa)$ with $\ep(\gamma) = 1$ (resp.\ $\ep(\gamma) = -1$)
\medskip

\noi holds. Then there exists an isomorphism of $\wGa$-modules \footnote{See Remark \ref{remark:twists} for the definition of $\sF(\ep)$.}
\begin{equation*}
\label{supp8a}
H_{\brC}^n(X, \sF(\ep)) \, \cong \, \bigoplus_{c\in C} \sF_c.
\end{equation*}
\end{coro}

We briefly address functorial properties for some of the cohomology groups introduced in this section. Let $f: X\to Y$ be a morphism of $(\sL, \wGa)$-spaces. By \eqref{apppullbackcoh} and \eqref{apppullbackcoh2} there exists natural morphisms of $\delta$-functors
\begin{eqnarray}
\label{pullbackcoh}
&& H^i(Y, \sF)\lra H^i(X, f^*(\sF))\qquad i\ge 0, \,\sF\in  \Sh(Y, R), \\
\label{pullbackcoh2}
&& H^i(Y, \wGa, \sF)\lra H^i(X, \wGa, f^*(\sF)) \qquad i\ge 0, \,\sF\in  \Sh(Y, \wGa, R).
\end{eqnarray}
Now assume that $(C, \io)$ is a $\wGa$-stable closed subspace of $Y$ disjoint from the image of $f$, i.e.\ we assume that $C_L$ is disjoint of $f_L(X_L)$ for every $L\in \sL$.\footnote{Note that this condition is equivalent to requiring that $\hat{f}(\hX)$ and $\io(C)$ are disjoint} It follows immediately from the definitions that 
the functor \eqref{pullbackcoh} for $i=0$ factors canonically in the form 
\begin{equation*}
\label{pullbacksections2}
H^0(Y, \sF)\lra H^0(Y\setminus \brC, \sF) \lra H^0(X, f^*(\sF)).
\end{equation*} 
The second homomorphism extends to a morphism of $\delta$-functors 
\begin{equation*}
\label{pullbackcoh3}
H^i(Y\setminus \brC, \sF)\lra H^i(X, f^*(\sF))\qquad i\ge 0,\, \sF\in  \Sh(Y, R)
\end{equation*}
so that the composition $H^i(Y, \sF) \to H^i(Y\setminus \brC, \sF)\to H^i(X, f^*(\sF))$ is the morphism \eqref{pullbackcoh}. Similarly, for $\wGa$-equivariant sheaves there exists a natural morphism of $\delta$-functors 
\begin{equation}
\label{pullbackcoh4}
H^i(Y\setminus \brC, \wGa, \sF)\lra H^i(X, \wGa, f^*(\sF)) \qquad i\ge 0, \,\sF\in  \Sh(Y, \wGa, R)
\end{equation}
so that \eqref{pullbackcoh2} factors in the form $H^i(Y, \wGa, \sF)\to H^i(Y\setminus \brC, \wGa, \sF)\to H^i(X, \wGa, f^*(\sF))$.

\section{Adelic Eisenstein Classes}
\label{section:eiscocycle}

\paragraph{The sheaf of locally polynomial distributions} In this section we let $V$ be an oriented $\bQ$-vector space of dimension $n$, let $\sL$ be a non-empty subset of $\Lat_V$ such that \eqref{sLclosed} holds and put $\Lambda= \Lambda(\sL)= \, \bigcup_{L\in \sL} L$. Moreover we fix a subgroup $\Ga$ of $\GL_{\bQ}(V)$ such that $\sL$ is $\wGa$-stable and put $\wGa := \Ga\ltimes \La\subseteq \Aff_{\bQ}(V)$. In the following we denote by $B$, $A$ and $T$ the $(\sL, \wGa)$-spaces of examples \ref{examples:torpt} (a), (b), (c). Recall that we have $\hB = \hLa := \prolim_{L\in \sL} \Lambda/L$ (see \eqref{Bhat})
and that a subset $U\subseteq \hB$ is open if and only if there exists $L\in \sL$ such that $U$ is $\hL$-stable. 

For a fixed ring $R$ we define an $R$-presheaf $\cD_{\lpol}=\cD_{\lpol, B}$ on $B$ by   
\begin{equation*}
\label{locpoldistsh1}
\cD_{\lpol}(U): =\, \ \cD_{\sL-\lpol}(U\cap \Lambda, R) \, =\, \Hom(\Int_{\sL\mloc, b}(U\cap \Lambda, \bZ), R)
\end{equation*}
for $U\subseteq \hLa$ open. For a pair of open subsets $U_1 \subseteq U_2\subseteq \hLa$ the restriction is the map 
\begin{equation*}
\label{locpoldistsh2}
\cD_{\lpol}(U_2)\lra \cD_{\lpol}(U_1), \quad \mu\mapsto \mu|_{U_1\cap \La}
\end{equation*}
(compare \eqref{distres}). By Lemmas \ref{lemma:cosheaf} and \ref{lemma:sheafdisc} $\cD_{\lpol}$ is a sheaf. Moreover it carries a canonical $\wGa$-action given by the maps \eqref{funct2}. 

Recall that $A= B\times V_{\bR}$ and that $T = A/\Lambda$ and that there are canonical morphisms $\pr: A\to T$ and $\pr_B : A\to B$. The pull-back of $\cD_{\lpol, B}$ under $\pr_B: A\to B$ will be denoted by $\cD_{\lpol}=\cD_{\lpol, A}$ as well. By Example \ref{example:torsh} there exists a natural $\Ga$-equivariant $R$-sheaf $\sD=\sD_{\lpol}$ on $T$ such that 
\begin{equation*}
\label{locpoldistsh3}
\pr^*(\sD)\, =\, \cD_{\lpol, A} \,=\, (\pr_B)^*(\cD_{\lpol, B}). 
\end{equation*}
Concretely, the section of $\sD$ over an open subset $U\subseteq \hT = \left(\hLa \times V_\bR\right)/\Lambda$ are the $\Lambda$-invariant elements of $\cD_{\lpol, A}(\pr^{-1}(U))$. 

\begin{df}
\label{df:locpoldistsh}
The sheaves $\cD_{\lpol, B}$, $\cD_{\lpol, A}$ and $\sD_{\lpol}$ will be called the sheaf of locally polynomial $R$-valued distributions on $B$, $A$ and $T$ respectively.
\end{df}

We are going to describe the sheaf $\sD_L$ on $T_L\cong V_{\bR}/L$ for $L\in \sL$. Let 
\begin{equation*}
\label{polmeas}
j : V_{\bR} \lra \hT = \left(\hLa \times V_\bR\right)/\Lambda, \quad v \mapsto (0, v) + \Lambda.
\end{equation*}
The composition $\pi_L\circ j$ is the universal covering of $T_L$, namely it is the map 
\begin{equation*}
\label{polmeas2}
q_L : V_{\bR} \lra V_{\bR}/L, \quad v \mapsto v + L.
\end{equation*}
Since the group of deck transformations of the covering $q_L$ is the group $L$, any $R[L]$-module $M$ defines a local system $\wM$ on $T_L$. Recall that its sections over an open subset $U\subseteq T_L$ are given by 
\begin{equation}
\label{locsys}
\wM(U)\, =\, \{ f\in C(q_L^{-1}(U), M)\mid  f(\la + v) = \la\cdot f(v)\,\, \forall\, \la \in L, v\in V_{\bR}\}.
\end{equation}

\begin{lemma}
\label{lemma:sdlayer}
We have $\sD_L= \widetilde{\cD_{\lpol}(L, R)}$, i.e.\ $\sD_L$ is the sheaf associated to the $R[L]$-module $\cD_{\lpol}(L, R)=\Hom(\Int_{\loc}(L, \bZ), R)$ on $T_L$. 
\end{lemma}

\begin{proof} Firstly, note that $(\cD_{\lpol, A})_L\in \Sh(A_L, R)$ is the pull-back of the sheaf $(\cD_{\lpol, B})_L$ on the discrete space $B_L=\Lambda/L$ under the projection $\pr_1: A_L= \Lambda/L \times V_{\bR}\to \Lambda/L$. 
Define $\io: V_{\bR} \to A_L= \Lambda/L \times V_{\bR}, v\mapsto (0, v)$ so that $\pr_L\circ \io = q_L$. Both maps 
$\pr_L: A_L \to T_L$ and $q_L:V_{\bR}\to T_L$ are Galois coverings with group of deck transformation $\Lambda$ and $L$ respectively. 
Hence the functors $(\pr_L)^*: \Sh(T_L, R) \to \Sh(A_L, \Lambda, R)$ and $(q_L)^*: \Sh(T_L, R) \to \Sh(V_{\bR}, L, R)$ are equivalences of categories. Since the pull-back of $\sD_L$ under $\pr_L:A_L\to T_L$ is the $\Lambda$-equivariant sheaf  $(\cD_{\lpol, A})_L$ on $A_L$, we see that the pull-back of $\sD_L$ under $q_L:V_{\bR}\to T_L$ is the $L$-equivariant sheaf  $\io^*(\cD_{\lpol, A})_L= (\pr_1 \circ \io)^*((\cD_{\lpol, B})_L)$. Note that $\pr_1 \circ \io$ is the constant map $V_{\bR} \to \Lambda/L, v\mapsto 0$. Hence $(q_L)^*(\sD_L)= \io^*(\cD_{\lpol, A})_L$ is the constant $L$-equivariant sheaf on $V_{\bR}$ associated to the $L$-module $(\cD_{\lpol, B})_L(\{0\})= \cD_{\lpol}(L, R)$.
\end{proof}

Next we determine the stalks of $\cD_{\lpol, B}$, $\cD_{\lpol, A}$ and $\sD_{\lpol}$.

\begin{lemma}
\label{lemma:stalksdl}
(a) There exists canonical isomorphisms 
\begin{equation*}
\label{stalkdb1} 
\beta_b: (\cD_{\lpol, B})_b \lra \cD_{\lpol, b}(\Lambda, R) \qquad \mbox{and} \qquad \beta^b: (\cD_{\lpol, B})^b \lra \cD_{\lpol}(\Lambda, R) 
\end{equation*} 
for every $b\in \hB$. Moreover the diagram 
\begin{equation}
\label{stalkdb2} 
\begin{CD}
(\cD_{\lpol, B})_b @>\beta_b >> \cD_{\lpol, b}(\Lambda, R) @.\hspace{1.5cm} (\cD_{\lpol, B})^b @>\beta^b >> \cD_{\lpol}(\Lambda, R)\\
@VV \wga V @VV \wga V\hspace{2cm} @VV \wga V @VV \wga V\\
(\cD_{\lpol, B})_{\wga^{-1}(b)} @>\beta_{\wga(b)} >> \cD_{\lpol, b}(\Lambda, R) @. \hspace{1.5cm} (\cD_{\lpol, B})^{\wga^{-1}(b)} @>\beta^{\wga^{-1}(b)} >> \cD_{\lpol}(\Lambda, R)
\end{CD}
\end{equation}
commutes for every $\wga\in \wGa$ and $b\in \hB$. 
\medskip

\noi (b) Similarly there exists canonical isomorphisms 
\begin{equation*}
\label{stalkda1} 
\beta_a: (\cD_{\lpol, A})_a \lra \cD_{\lpol, b}(\Lambda, R) \qquad \mbox{and} \qquad \beta^a: (\cD_{\lpol, A})^a \lra \cD_{\lpol}(\Lambda, R) 
\end{equation*} 
for every $a\in \hA$. Also the diagram analogous to \eqref{stalkdb2} for $\cD_{\lpol, A}$ commutes for every $a\in \hA$ and 
$\wga\in \wGa$. 
\medskip

\noi (c) There exists a canonical isomorphisms
\begin{equation}
\label{stalkdt}
\beta_t: \sD_t\, \lra \, \cD_{\lpol, b}(\pr^{-1}(t), R)\qquad \mbox{and} \qquad \beta^t:  \sD^t\, \lra \, \cD_{\lpol}(\pr^{-1}(t), R)
\end{equation}
for every $t\in \hT$. Moreover the diagram 
\begin{equation*}
\label{stalkdt2} 
\begin{CD}
\sD_t @>\beta_t >> \cD_{\lpol, b}(\pr^{-1}(t), R)@. \hspace{1.5cm} \sD^t @>\beta^t >> \cD_{\lpol}(\pr^{-1}(t), R)\\
@VV \gamma V @VV \gamma V \hspace{2cm} @VV \gamma V @VV \gamma V\\
\sD_{\gamma^{-1}(t)} @>\beta_{\gamma^{-1}(t)} >> \cD_{\lpol, b}(\pr^{-1}(\gamma(t)), R) @. \hspace{1.5cm} \sD^{\gamma^{-1}(t)} @>\beta^{\gamma^{-1}(t)} >> \cD_{\lpol}(\pr^{-1}(\gamma^{-1}(t)), R)
\end{CD}
\end{equation*}
commutes for every $\gamma\in \Ga$ and $t\in \hT$. In particular if $t\in \hT$ is a $\Ga$ fixed-point then the maps \eqref{stalkdt} are isomorphisms of $R[\Ga]$-modules.
\end{lemma}

\begin{proof} (a) Let $b\in \hB= \hLa$ and put $b_L :=\pi_L(b) \in \Lambda/L$ for $L\in \sL$. Since $\hLa=\bigcup_{L\in \sL} \hL$ we have $b\in \hL$ for $L\in \sL$ sufficiently large. It follows 
\begin{equation*}
\label{stalkdb3}
(\cD_{\lpol, B})_b =\dlim_{L\in \sL} \cD_{\lpol, B}(\pi_L^{-1}(b_L))= \dlim_{L\in \sL, b\in \hL} \cD_{\lpol, B}(\hL)\,= \, \dlim_{L\in \sL} \cD_{\lpol}(L) = \cD_{\lpol, b}(\Lambda, R)
\end{equation*}
and similarly $(\cD_{\lpol, B})^b= \cD_{\lpol}(\Lambda, R)$. (b) follows immediately from (a) and Prop.\ \ref{prop:pullbackstalk}. 

For (c) note that for the stalk $\sD_t$ for $t\in \hT$ can be identified with the $\Lambda$-invariant elements of $\prod_{a\in \pr^{-1}(t)} (\cD_{\lpol, A})_a$. Thus according to (b), the $R$-module $\sD_t$ can (and will) be identified with the set of maps $\mu: \pr^{-1}(t)\to \cD_{\lpol, b}(\Lambda, R), a\mapsto \mu_a$ that satisfy $\mu_{\la + a} = (\tau_{-\la})_*(\mu_a)$ for every $\la\in \Lambda$ and $a\in \pr^{-1}(t)$. For such $\mu$ choose $a\in \pr^{-1}(t)$ and define $\widetilde{\mu} := (\psi_a)_*(\mu_a)\in\cD_{\lpol, b}(\pr^{-1}(t), R)$ where $\psi_a$ is the $\Lambda$-equivariant bijection
$\psi_a: \Lambda\to \pr^{-1}(t), \la \mapsto \la +a$. We show that $\widetilde{\mu}$ is independent of the choice of $a$. For that let 
$a'\in \pr^{-1}(t)$ and let $\la\in \Lambda$ with $a'=\la +a$. Then we have 
\[
 (\psi_a')_*(\mu_a') =  (\psi_{\la + a})_*(\mu_{\la + a}) = ((\psi_{\la + a})_*\circ (\tau_{-\la})_*) (\mu_a) =  (\psi_a)_*(\mu_a).
 \]
It is easy to see that the map $\mu \mapsto \widetilde{\mu}$ is the desired isomorphism $\beta_t$. The existence of second isomorphism $\beta^t$ is proved similarly. 
\end{proof}

\begin{remark}
\label{remark:stalksdlayer}
\rm For $L\in \sL$ the stalks of the sheaf $\sD_L$ on $T_L$ admit a description similar to that for $\sD$ in Lemma \ref{lemma:stalksdl} (c). Namely, there exists a canonical isomorphism 
\begin{equation}
\label{stalksdlayer}
(\sD_L)_x \, \lra \, \cD_{\lpol}(q_L^{-1}(x), R)
\end{equation}
for every $x\in T_L$. Indeed, by Lemma \ref{lemma:sdlayer} and \eqref{locsys} the stalk $(\sD_L)_x$ can be identified with the set of maps $\mu: q_L^{-1}(x) \to \cD_{\lpol}(L, R), h\mapsto \mu_h$ that satisfy $\mu_{\la + h} = \la \cdot \mu_h$ for every $h\in \pr_L^{-1}(t)$ and $\la\in L$. Given such $\mu$ we define its image under \eqref{stalksdlayer} as the unique $\widetilde{\mu}\in \cD_{\lpol}(q_L^{-1}(t), R)$ that satisfies 
\begin{equation*}
\label{stalksdlayer2} \int_{q_L^{-1}(t)} f(h) d\widetilde{\mu} (h)\, =\, \int_L f(\la +h') d(\mu_h')(\la).
\end{equation*}
for every $f\in \Int_{\loc}(q_L^{-1}(t), \bZ)$ and $h'\in q_L^{-1}(t)$.
\end{remark}

Next we are going to determine the cohomology of $T$ with coefficients in the twisted sheaf $\sD(\ep)$ (see Remark \ref{remark:twists}). Here $\ep$ denotes the sign character of $\Aff(V)$ introduced in \eqref{sign}, i.e.\ we have $\ep((\gamma, v))=\sgn(\det(\gamma))$ for every $(\gamma, v)\in \Aff(V)=\GL(V)\ltimes V$. Firstly, note that according to Cor.\ \ref{coro:cohodisc} we have 
\begin{equation}
\label{cohdlab}
H^i(A, \cD_{\lpol, A}) \,=\, H^i(B, \cD_{\lpol, B}) \,=\, \left\{\begin{array}{cc} \cD_{\lpol}(\Lambda, R) & \mbox{if $i=0$,}\\
0 & \mbox{if $i\ge 1$.}
\end{array}
\right.\\
\end{equation}
The second equality follows from Cor.\ \ref{coro:cohodisc}. The first equality can be seen using the Leray spectral sequence \eqref{lerayss} associated to the morphism $\pr_B : A\to B$. For that note that according to Prop.\ \ref{prop:pullbackex} we have $(\pr_B)_*\cD_{\lpol, A} = \cD_{\lpol, B}$ and $R^i (\pr_B)_*\cD_{\lpol, A}=0$ for $i\ge 1$ since the fibers of 
$\pr_{B, L} : A_L \approx \bR^n \times B_L\to B_L$ are contractible for every $L\in \sL$. 

\begin{prop}
\label{prop:cohomld} 
(a) For $i\ge 0$ we have 
\begin{equation}
\label{cohpolmeas}
H^i(T, \sD(\ep)) \,=\, \left\{\begin{array}{cc} R & \mbox{if $i=n$,}\\
0 & \mbox{if $i\ne n$}
\end{array}
\right.\\
\end{equation}
and $H^i(T, \Ga, \sD(\ep)) \,=\, \Ext_{R[\Ga]}^{i-n}(R, R)$. 
\medskip

\noi (b) Let $(C, \io)$ be a finite $\Ga$-stable closed subspace of $T$. Then, 
\begin{equation}
\label{cohsupppolmeas1}
H_{\brC}^i(T, \sD(\ep)) \,=\, \left\{\begin{array}{cc} \cD_{\lpol, b}(\pr^{-1}(C), R) & \mbox{if $i=n$,}\\
0 & \mbox{if $i\ne n$}
\end{array}
\right.\\
\end{equation}
and 
\begin{equation}
\label{cohsupppolmeas2}
H_{\brC}^i(T, \Ga, \sD(\ep)) \,=\, \left\{\begin{array}{cc} \cD_{\lpol, b}(\pr^{-1}(C), R)^{\Ga} & \mbox{if $i=n$,}\\
0 & \mbox{if $i\le n-1$.}
\end{array}
\right.\\
\end{equation}
\medskip

\noi (c) Under the assumptions of (b) we have 
\begin{equation}
\label{cohsupppolmeas3}
H^i(T\setminus \brC, \sD(\ep)) \,=\, \left\{\begin{array}{cc} \ker \left(\aug: \cD_{\lpol, b}(\pr^{-1}(C), R) \to R\right) & \mbox{if $i=n-1$,}\\
0 & \mbox{if $i\ne n-1$}
\end{array}
\right.
\end{equation}
and 
\begin{equation}
\label{cohsupppolmeas4}
H^i(T\setminus \brC, \Ga, \sD(\ep)) \,=\,\left\{\begin{array}{cc} \ker \left(\aug: \cD_{\lpol, b}(\pr^{-1}(C), R)^{\Ga} \to R\right) & \mbox{if $i=n-1$,}\\
0 & \mbox{if $i\ne n-1$.}
\end{array}
\right.
\end{equation}
\end{prop}

On the right side of \eqref{cohsupppolmeas1} -- \eqref{cohsupppolmeas4} we have identified $C$ with $\io(C)\subseteq \hT$ so that $\pr^{-1}(C)$ is a $\Lambda$-subset of $\hA$. The map $\aug: \cD_{\lpol, b}(\pr^{-1}(C), R) \to R$ is given by evaluating $\mu\in \cD_{\lpol, b}(\pr^{-1}(C), R)$ 
at the constant function $\equiv 1$ on $\pr^{-1}(C)$.

\begin{proof} (a) Consider the spectral sequence \eqref{hsscomp} for the $\Lambda$-equivariant cohomology of $A$ with coefficients in $\cD_{\lpol, A}(\ep)$
\begin{equation*}
\label{hsscompa}
E_2^{rs} = \Ext_{R[\Lambda]}^r(R, H^s(A, \cD_{\lpol}(\ep))) \, \Longrightarrow\, E^{r+s}=H^{r+s}(A, \Lambda, \cD_{\lpol}(\ep)).
\end{equation*}
By \eqref{cohdlab} the spectral sequence degenerates, i.e.\ we have $E_2^{i0}= E^i$ for every $i\ge 0$. Moreover the fact that the functor \eqref{unicov2} is an equivalence of categories implies $H^{\bu}(T, \sD) \cong H^{\bu}(A, \Lambda, \cD_{\lpol, A})$. Together with Prop.\ \ref{prop:koszul2} we conclude
\begin{equation*}
\label{cohpolmeas2}
H^i(T, \sD(\ep)) \,=\,\Ext_{R[\Lambda]}^i(R, H^0(A, \cD_{\lpol}(\ep))) \,= \,\Ext_{R[\Lambda]}^i(R, \cD_{\lpol}(\Lambda, R)(\ep))\, =\, \left\{\begin{array}{cc} R & \mbox{if $i=n$,}\\
0 & \mbox{if $i\ne n$.}
\end{array}
\right.\\
\end{equation*}
The last assertion follows from \eqref{cohpolmeas} by applying the spectral sequence \eqref{hsscomp} to $X=\hT$ and $\cF=\sD$.

(b) The first assertion follows from Cor.\ \ref{coro:gpt2} and Lemma \ref{lemma:stalksdl} (c)
and the second from the first and Lemma \ref{lemma:nogrss}.

(c) The equality \eqref{cohsupppolmeas3} follows from from (a) and (b) using the long exact sequences \eqref{supp6} and \eqref{cohsupppolmeas4} follows from \eqref{cohsupppolmeas3} and Lemma \ref{lemma:nogrss}.
\end{proof}

Let $t\in \hT$ be a $\Ga$ fixed-point and put $t_L:=\pi_L(t)\in T_L$ for every $L\in \sL$. We define a morphism of $(\sL, \wGa)$-spaces 
\begin{equation*}
\label{iotat}
\io_t :=  t + \io : B \lra T
\end{equation*}
where $\io: B\to T$ was defined in \eqref{zerosec}. More precisely we define $\io_{t, L}: B_L \to T_L$ by $\io_{t, L}(b) = \io_L(b) + t_L$ for every $L\in \sL$ and $b\in B_L$. For the $\Ga$-equivariant cohomology of $B$ with coefficients in $\io_t^*(\sD(\ep))$ we obtain 

\begin{prop}
\label{prop:cohomld2} 
For $i\ge 0$ we have 
\begin{equation}
\label{cohbtdl}
H^i(B, \io_t^*(\sD(\ep))) \,=\, \left\{\begin{array}{cc}  \cD_{\lpol}(\pr^{-1}(t), R)(\ep)& \mbox{if $i=0$,}\\
0 & \mbox{if $i\ge 1$}
\end{array}
\right.\\
\end{equation}
and
\begin{equation*}
\label{cohbtdl2}
H^i(B, \Ga, \io_t^*(\sD(\ep))) \,=\, H^i(\Ga, \cD_{\lpol}(\pr^{-1}(t), R)(\ep)).
\end{equation*}
\end{prop}

\begin{proof} By Cor.\ \ref{coro:cohodisc} (a) we have $H^i(B, \io_t^*(\sD(\ep)))=0$ if $i\ge 1$. Moreover by Example \ref{example:stalkb}, Prop.\ \ref{prop:pullbackstalk} and Lemma \ref{lemma:stalksdl} (c) we have 
\begin{equation*}
\label{cohbtdl3}
H^0(B, \io_t^*(\sD(\ep))) \, =\, (\io_t^*(\sD(\ep)))^0 \, =\, \sD^t(\ep)\, =\, \cD_{\lpol}(\pr^{-1}(t), R)(\ep).
\end{equation*}
This proves \eqref{cohbtdl}. The second assertion follows from the first and Cor.\ \ref{coro:cohodisc} (b).
\end{proof}

\paragraph{Adelic Eisenstein Classes} In this section we continue with the set-up of the last section but for specific $V$ and $\sL$. Namely we choose the set-up of Example \ref{examples:fracideals} (b), i.e.\ $F$ denotes a number field of degree $d$ over $\bQ$, $S$ a finite set of nonarchimedean places of $F$, $V$ and $F$-vector space of dimension $m$ (so that $n=dm$), $\cM$ a finitely generated $\cO_S$-submodule of $V$ with $\cM\otimes_{\cO_S} F=V$ and $\sL$ the set of finitely generated $\cO_F$-submodules $L\subseteq \cM$ satisfying $L\otimes_{\cO_F} \cO_S = \cM$ (note then that we have $\Lambda(\sL) = \cM$). Moreover we fix a subgroup $\Ga$ of $\GL_{\cO_S}(\cM)$. Note that $\sL$ is $\Ga$-stable and that it satisfies \eqref{sLclosed}. According to Example \ref{examples:torpt} (e) we have $\hT = (V\otimes_F \bA^S)/\cM$. As in the last section let $R$ be a ring and let $\sD=\sD_{\lpol}$ be the associated sheaf of locally polynomial $R$-valued distributions on $T$.

Let $v\in V$ and put $t=v+ \cM \in V/\cM \subseteq V\otimes_F \bA^S/\cM= \hT$. Recall that by Lemma \ref{lemma:stalksdl} (c) there exists canonical isomorphisms
\begin{equation*}
\label{stalkdt2a}
\beta_t: \sD_t\, \lra \, \cD_{\lpol, b}(v + \cM, R)\qquad \mbox{and} \qquad \beta^t:  \sD^t\, \lra \, \cD_{\lpol}(v + \cM , R).
\end{equation*}
According to Cor.\ \ref{coro:ordtinv} if the order of $t$ in $V/\cM$ is invertible in $R$ then the targets of these maps can be identified with $\cD_{\lpol, b}(\cM, R)$ and $\cD_{\lpol}(\cM, R)$ respectively, i.e.\ in this case there exists canonical isomorphisms 
\begin{equation}
\label{stalkdt3}
\beta_t: \sD_t\, \lra \, \cD_{\lpol, b}(\cM, R)\qquad \mbox{and} \qquad \beta^t:  \sD^t\, \lra \, \cD_{\lpol}(\cM , R).
\end{equation}
Let $C\subseteq V/\cM$ be a finite $\Ga$-stable subset. By Example \ref{examples:closedsub} (b) the pair $(C, \io)$ is a closed subspaces of $T$ (where $\io=\incl: C\hra  V/\cM\hra \hT$ is the inclusion). Thus by \eqref{stalkdt3} and Prop.\ \ref{prop:cohomld} (c) we obtain

\begin{prop}
\label{prop:cohomldtor} 
Let $C$ be a finite $\Ga$-stable subset of $V/\cM$ and assume that the order of every element in $C$ is invertible in $R$. Then,
\begin{equation*}
\label{cohsupppolmeas4a}
H^i(T\setminus \brC, \Ga, \sD(\ep)) \,=\,\left\{\begin{array}{cc} \ker\left(\aug: \Maps_{\Ga}(C, \cD_{\lpol, b}(\cM, R)) \to  R\right) & \mbox{if $i=n-1$,}\\
0 & \mbox{if $i\ne n-1$.}
\end{array}
\right.
\end{equation*}
\end{prop}

Here the homomorphism $\aug$ is given by mapping $\mu\in \Maps_{\Ga}(C, \cD_{\lpol, b}(\cM, R))$ to the sum (over $c\in C$) of the evaluation of $\mu(c)$ at the constant function $\equiv 1$ on $\cM$.

Let $v_0\in V$ and assume that $t=v_0+\cM \in V/\cM \subseteq \hT$ is a $\Ga$ fixed-point. According to Prop.\ \ref{prop:cohomld2} for the $\Ga$-equivariant cohomology of $B$ with coefficients in $\io_t^*(\sD(\ep))$ we get
\begin{equation}
\label{cohbtdl2tor}
H^i(B, \Ga, \io_t^*(\sD(\ep))) \,=\, H^i(\Ga, \cD_{\lpol}(v_0+ \cM, R)(\ep))
\end{equation}
for every $i\ge 0$. If the order of $t$ is invertible in $R$ then this simplifies to
\begin{equation*}
\label{cohbtdl2tor2}
H^i(B, \Ga, \io_t^*(\sD(\ep))) \,=\, H^i(\Ga, \cD_{\lpol}(\cM, R)(\ep)).
\end{equation*}

With this preparation we are able to define our {\it adelic Eisenstein classes}. The construction is modelled after that of the Eisenstein classes of Beilinson, Kings and Levin (\cite{bkl}, Def.\ 3.32). We fix a finite $\Ga$-stable non-empty subset $C\subseteq V/\cM$. In the following we assume that the order of every element of $C$ is invertible in $R$.\footnote{Note that by Remark \ref{remark:ordVmodcM} this holds if the residue characteristic of every prime in $S$ is invertible in $R$.}
Let $t\in V/\cM$ be a $\Ga$ fixed-point that does not lie in $C$. We need the following simple

\begin{lemma}
\label{lemma:tcdisjoint}
The image of $\io_t$ is disjoint from $C$, i.e.\ the image of $\io_{t, L}: B_L \to T_L$ is disjoint from $\pi_L(C)$ for every $L\in \sL$. 
\end{lemma}

\begin{proof} It suffices to see that the image of $\widehat{\io_t}: \hB\to \hT$ is disjoint from $C$. For that we write elements of 
$V\otimes_F \bA^S$ as pairs $(x,y)$ with $x\in V\otimes_F \bA_f^S$ and $y\in V\otimes_F F_{\infty} = V_{\bR}$. Concretely, $\widehat{\io_t}$ is the map $V\otimes_F \bA_f^S \to V\otimes_F \bA^S/\cM, x \mapsto (x+v, v) \mod \cM$. Assume that $c\in C$, $c= w + \cM\in V/\cM$ lies in the image of $\widehat{\io_t}$. Thus there exists 
$x\in V\otimes_F \bA_f^S$ with $(x+v, v) = (w, w) \mod \cM$, i.e.\ $(x+v, v) = (w, w) + (m, m)$ for some $m\in \cM$. 
It follows $x=0$ and $v-w = m\in \cM$ hence $t=c$, a contradiction.
\end{proof}

\begin{remark}
\label{remark:imiotl}
\rm The image of $\io_{t,L}: \cM/L \subseteq V_{\bR}/L$ can be characterized as follows. It consists of those elements $t_0\in V/L\subseteq V_{\bR}/L$ that are mapped to $t$ under the canonical projection $V/L\to V/\cM$. 
\end{remark}
 
By \eqref{pullbackcoh4} the morphism $\io_t$ induces a homomorphisms 
\begin{equation*}
\label{pullbacktorus}
H^i(T\setminus \brC, \Ga, \sD(\ep)) \lra H^i(B, \Ga, \io_t^*(\sD(\ep)))
\end{equation*}
for every $i\ge 0$. Together with Prop.\ \ref{prop:cohomldtor} and \eqref{cohbtdl2tor} we obtain for $i=n-1$ the homomorphism
\begin{equation}
\label{pullbacktorus2}
 \ker\left(\aug: \Maps_{\Ga}(C, \cD_{\lpol, b}(\cM, R)) \to  R\right) \lra H^{n-1}(\Ga, \cD_{\lpol}(v+ \cM, R)(\ep)). 
\end{equation}
By viewing $R\subseteq R[\cM]$ as submodules of $\cD_{\lpol, b}(\cM, R)$ we obtain
\[
R[C]\,=\,\Maps(C, R)\subseteq \Maps(C, R[\cM])\,\subseteq \,\Maps_{\Ga}(C, \cD_{\lpol, b}(\cM, R))
\]
hence $(R[C]^0)^{\Ga}\subseteq \ker(\aug: \Maps_{\Ga}(C, \cD_{\lpol, b}(\cM, R)) \to R)$
where $R[C]^0:=\ker(\deg: R[C]\to R)$. 

\begin{df}
\label{df:adeleis}
Let 
\begin{equation}
\label{eisenmap}
\bEis(t): (R[C]^0)^{\Ga}\lra H^{n-1}(\Ga, \cD_{\lpol}(v+ \cM, R)(\ep)) \quad \alpha \mapsto \bEis_{\alpha}(t)
\end{equation}
be the restriction of the map \eqref{pullbacktorus2} to the subgroup $(R[C]^0)^{\Ga}$. The element 
\begin{equation}
\label{eisenclass}
\bEis_{\alpha}(t)\in H^{n-1}(\Ga, \cD_{\lpol}(v+ \cM, R)(\ep))
\end{equation}
will be called adelic Eisenstein class associated to $t$ and $\alpha$. 
\end{df}

\begin{remark}
\label{remark:vareisenclass} \rm If the order of $t$ is invertible in $R$ as well then the coefficients of the cohomology group in \eqref{eisenclass} can be identified with $\cD_{\lpol}(\cM, R)$, i.e.\ in this case 
we have 
\begin{equation*}
\label{eisenclassa}
\bEis_{\alpha}(t)\in H^{n-1}(\Ga, \cD_{\lpol}(\cM, R)(\ep)).
\end{equation*}
\end{remark}

We want to relate our adelic Eisenstein classes \eqref{eisenclass} to the classes (\cite{bkl}, Def.\ 3.32) of Beilinson, Kings and Levin. For $L\in \sL$ put $t_L = \pi_L(t)\in T_L$. Let $t_0\in T_L$ be contained in the image of the morphism $\io_t:B\to T$, i.e.\ we assume that $t_0$ is of the form $t_0= t_L + t_1\in V/L\subseteq T_L$ for some $t_1 = m + L \in \cM/L$ so that $t_0 = v_0 + L$ with $v_0 := v + m\in V$. Let $\Ga_0$ be a subgroup of $\Ga$ that stabilises 
  $L$ and $t_0$, i.e.\ we have $\ga(L) = L$ and $\ga(t_0) = t_0$ for every $\ga\in \Ga_0$. Since $q_L^{-1}(t_0) = v_0 + L \subseteq v_0 + \cM = v+ \cM\subseteq V$ we can consider the $\Ga_0$-equivariant restriction (see \eqref{distres})
\begin{equation}
\label{eisenclasscomp}
\cD_{\lpol}(v+ \cM, R)\, \lra\, \cD_{\lpol}(v_0 + L, R), \quad \mu \mapsto \mu|_{v_0 + L}.
\end{equation}

Recall that the logarithm sheaf $\sLog_L$ on $T_L$ introduced in (\cite{bkl}, \S 3.4) is defined as the local system associated to the $R[L]$-algebra $R[\![L]\!]\cong \cD_{\pol}(L, R)$. The $\Ga_0$-action on $L$ provides $\sLog_L$ with the structure of a $\Ga_0$-equivariant sheaf. The stalks of $\sLog_L$ admit a similar description as those of $\sD_L$ (see 
Remark \ref{remark:stalksdlayer}), namely for $x\in T_L$ we have $\sLog_{L, x}\cong \cD_{\pol, L}(q_L^{-1}(x), R)$. Since $t_0$ is stabilised by $\Ga_0$ and disjoint from $C$ there exists a canonical map (see \cite{bkl}, Def.\ 3.32)
\begin{eqnarray}
\label{eisenbkl}
&&\Eis_L(t_0): (R[C]^0)^{\Ga_0}\lra H^{n-1}(\Ga_0, (\sLog_L)_{t_0}(\ep)) \cong H^{n-1}(\Ga_0, \cD_{\pol}(v_0 + L , R)(\ep)),\\
&&\hspace{3cm} \alpha \longmapsto \Eis_{L, \alpha}(t_0).\nonumber
\end{eqnarray}
Its definition is similar to that of the map \eqref{eisenmap} above. We are going to review the main steps of the construction in the beginning of the proof of the Prop.\ \ref{prop:eisencompare} below. 

\begin{remark} \rm If the order of $t_0$ is invertible in $R$ then there exists a canonical isomorphism of $\Ga_0$-modules $\cD_{\pol}(v_0 + L, R)=(\sLog_L)_{t_0}\cong \cD_{\pol}(L, R)\, =\, R[\![L]\!]$ (see \cite{bkl}, 3.24) so that $\Eis_{L, \alpha}(t_0)\in H^{n-1}(\Ga_0, R[\![L]\!](\ep))$ in this case. 
\end{remark}

Composing \eqref{eisenclasscomp} with the canonical map $\cD_{\lpol}(v_0 + L, R)\to \cD_{\pol}(v_0 + L , R)$ yields a $\Ga_0$-equivariant homomorphism $\cD_{\lpol}(v+ \cM, R)\, \lra\, \cD_{\pol}(v_0 + L, R)$. Together with the inclusion $\Ga_0\hra \Ga$ it induces a homomorphism 
\begin{equation}
\label{eisenclass2}
H^{n-1}(\Ga, \cD_{\lpol}(v+ \cM, R)(\ep)) \, \lra \, H^{n-1}(\Ga_0, \cD_{\pol}(v_0 + L, R)(\ep)).
\end{equation}

\begin{prop}
\label{prop:eisencompare} 
For every $\alpha\in (R[C]^0)^{\Ga}$ the image of the adelic Eisenstein class $\bEis_{\alpha}(t)$ under the homomorphism \eqref{eisenclass2} is the Eisenstein class $\Eis_{L, \alpha}(t_0)$ of Beilinson, Kings and Levin. 
\end{prop}

\begin{proof} Let 
\begin{equation}
\label{eisenbkl2}
H^{\bu}(T_L\setminus C_L, \Ga_0, \sLog_L(\ep)) \, \lra \, H^{\bu}(\{t_0\}, \Ga_0, \io_{t_0}^*(\sLog_L(\ep))) \, \cong \, H^{\bu}(\Ga_0, (\sLog_L)_{t_0}(\ep))
\end{equation}
be the canonical homomorphism of equivariant cohomology groups induced by the $\Ga_0$-equivariant embedding $\io_{t_0}:\{t_0\}\hra T_L\setminus C_L$. Using similar arguments as in the proof of \eqref{cohsupppolmeas4} above one shows that 
\begin{equation}
\label{eisenbkl3}
H^{n-1}(T_L\setminus C_L, \Ga_0, \sLog_L(\ep)) \, \cong \, H^{n-1}(T_L\setminus C_L, \sLog_L(\ep))^{\Ga_0}\, \cong \,\ker((\bigoplus_{c\in C_L} \sLog_c)^{\Ga_0} \to R).
\end{equation}
Since all points in $C_L$ are torsion points of $T_L$ whose orders are invertible in $R$ there are canonical isomorphisms of $\Ga_0$-modules $\bigoplus_{c\in C_L} \sLog_c \cong \Maps(C, \cD_{\pol}(L , R))$ (see \cite{bkl}, 3.24). Therefore the source of the map \eqref{eisenbkl2} in degree $n-1$ can be identified with 
$\ker\left(\aug: \Maps_{\Ga_0}(C, \cD_{\pol}(L , R)) \to R\right)$. The map \eqref{eisenbkl} is the restriction of \eqref{eisenbkl2} (in degree $n-1$) to the subgroup $(R[C]^0)^{\Ga_0}$.

Note that in the construction of \eqref{eisenbkl} we could have replaced the sheaf $\sLog_L$ with $\sD_L$ and the coefficients $\cD_{\pol}(v_0 + L , R)(\ep)$ in the target of \eqref{eisenbkl} with $\cD_{\lpol}(v_0 + L , R)(\ep)$. Indeed, it follows again from Prop.\ \ref{prop:cohomld} (c) (applied to $T= T^{\sL(L)}$; compare Remark \ref{remark:slfinal}) that (\cite{bkl}, Prop.\ 3.24 and Cor.\ 3.28) hold as well if we replace $\sLog_L$ with the sheaf $\sD_L$. 

Consider the diagram 
\begin{equation*}
\label{eisenbkl4}
\scriptsize{
\begin{CD}
\ker\left(\aug: \Maps_{\Ga_0}(C, \cD_{\pol}(L , R)) \to R\right) @>\eqref{eisenbkl3}>\cong > H^{n-1}(T_L\setminus C_L, \Ga_0, \sLog_L(\ep))@> 1 >> H^{n-1}(B_L, \Ga_0, \io_{t, L}^*\sLog_L(\ep))\\
@A 2 AA@A 3 AA@A 4 AA\\
\ker\left(\aug: \Maps_{\Ga_0}(C, \cD_{\lpol}(L, R)) \to R\right) @> 5 >\cong > H^{n-1}(T_L\setminus C_L, \Ga_0, \sD_L(\ep))@> 6 >> H^{n-1}(B_L, \Ga_0, \io_{t, L}^*\sD_L(\ep))\\
@V 7 VV@V 8 VV@V 9 V\cong V\\
\ker\left(\aug: \Maps_{\Ga_0}(C, \cD_{\lpol, b}(\cM, R)) \to R\right) @>\eqref{pullbacktorus2}> \cong >H^{n-1}(T\setminus \brC, \Ga_0, \sD(\ep))@>\eqref{pullbacktorus}>> H^{n-1}(B, \Ga_0, \io_t^*\sD(\ep))\\
@A\incl AA@A 10 AA@A 11 AA\\
\ker\left(\aug:\Maps_{\Ga}(C, \cD_{\lpol, b}(\cM, R)) \to R\right)  @>\eqref{pullbacktorus2}>\cong > H^{n-1}(T\setminus \brC, \Ga, \sD(\ep))@>\eqref{pullbacktorus}>> H^{n-1}(B, \Ga, \io_t^*\sD(\ep))\\
\end{CD}
}
\end{equation*}
Here the maps $2$, $3$ and $4$ are all induced by the canonical map $\cD_{\lpol}(L, R)\to \cD_{\pol}(L , R)$ respectively the associated morphism of sheaves $\sD_L\to \sLog_L$. The isomorphism $5$ is defined similarly to  \eqref{eisenbkl3} and the 
maps $1$ and $6$ similarly to \eqref{pullbackcoh4}. 
$7$ is induced by the canonical map $\cD_{\lpol}(L, R)\to \cD_{\lpol, b}(\cM, R), \mu \mapsto \mu_!$ (see \eqref{bdistext0}). The map $8$ is defined by first identifying source and target with $H^{n-1}(T_L\setminus C_L, \sD_L(\ep))^{\Ga_0}$ and $H^{n-1}(T\setminus \brC, \sD(\ep))^{\Ga_0}$ respectively. It is then obtained by passing in the second map of \eqref{supp5} to $\Ga_0$-invariants. The map $9$ is the isomorphism \eqref{compequivcoh} of Cor.\ \ref{coro:classicequiv} (for $X=B$ and $\sF=\sD$). 
$10$ and $11$ are the obvious maps. 

Let $c\in H^{n-1}(B_L, \Ga_0, \io_{t,L}^*\sD_L(\ep))$ denote the image of $\alpha\in  (R[C]^0)^{\Ga}\subseteq (R[C]^0)^{\Ga_0}$ under the composition of $5$ and $6$. By \eqref{cohbtdl2tor} we can identify the group in the lower right corner of the diagram with $H^{n-1}(\Ga, \cD_{\lpol}(v +\cM, R)(\ep))$. Hence $\alpha$ is mapped under the composition of the lower two horizontal maps to $\bEis_{\alpha}(t)$.
Therefore the commutativity of the diagram shows that $\bEis_{\alpha}(t)$ is mapped to $c$ under the composition of $11$ with the inverse of $9$. On the other hand by the definition of \eqref{eisenbkl} the image of $c$ under the composition of $4$ with the map 
\[
H^{n-1}(B_L, \Ga_0, \io_{t, L}^*\sLog_L(\ep))\,\lra \, H^{n-1}(\{t_1\}, \Ga_0, \io^*(\io_{t, L}^*\sLog_L)(\ep))\, \cong\,  H^{n-1}(\Ga_0, (\sLog_L)_{t_0}(\ep))
\]
induced by the inclusion $\io: \{t_1\}\hra B_L$ is easily seen to be equal to $\Eis_{\alpha}(t_0)$. The proves the assertion.
\end{proof}

\section{Eisenstein classes and special values of partial zeta functions} 
\label{section:partzeta}

The aim of this section is to relate the adelic Eisenstein classes defined in the previous section to special values of partial zeta functions and to Stickelberger elements. Throughout this section $F$ denotes a totally real number field of degree $n\ge 2$ over $\bQ$. We choose an ordering $\xi_1, \ldots, \xi_n: F\to \bR$ of the set of field embeddings $\Hom(F, \bR)$. Note that this choice provides $F_{\infty} = F\otimes \bR$ with an orientation.
We recall the notion of a partial zeta function associated to a ray class of $F$. For that we fix an ideal $\fm\subseteq \cO_F$, $\fm\ne (0)$. For a ray class $\fA\in \cI^{\fm}/\cP^{\fm}$ the partial zeta function $\zeta(\fm, \fA, s)$ is defined as 
\begin{equation*}
\label{partzeta1}
\zeta(\fm, \fA, s)\,\, =\,\, \sum_{\fa\in \fA, \fa \subseteq \cO_F} \Norm(\fa)^{-s}
\end{equation*}
for $\Re(s) >1$. It admits an analytic continuation to the whole complex plane except for a single simple pole at $s=1$. 

Recall \eqref{partzeta2} that there are also partial zeta function $\zeta_S(\sigma,s)$ associated to an elements $\sigma$ of the Galois group $G$ of an abelian extension of $K/F$. It is given by $\zeta_S(\sigma, s)=\sum_{(\fa,S)=1, \sigma_{\fa}=\sigma}\, \Norm(\fa)^{-s}$ if $\Real(s) >1$. 
Here $S$ is a finite set of nonarchimedean places of $F$ containing all places that are ramified in $K$ and the sum is taken over all ideals $\fa\subseteq \cO_F$ that are relatively prime to the elements in $S$ and such that their image $\sigma_{\fa}\in G$ under the Artin map is equal to $\sigma$. If $\fm$ is the nonarchimedean part of a cycle of declaration of $K/F$ (cf.\ \cite{neukirch}, p.\ 103) so that $K\subseteq F^{\fm}$ and if $S$ consists of all prime divisors of $\fm$ then we have 
\begin{equation}
\label{partzeta3}
\zeta_S(\sigma,s)\,\, =\,\,\sum_{\fA} \zeta(\fm, \fA, s)
\end{equation} 
Here the sum is taken over the finite number of ray classes $\fA\in \cI^{\fm}/\cP^{\fm}$ that are mapped to $\sigma$ under the Artin map. 

Recall as well the $T$-smoothed Stickelberger element $\Theta_{S, T}(K/F, s)$ defined in \eqref{stick2} in the introduction (where $T$ is an additional finite set of nonarchimedean places of $F$ disjoint from $S$). We will consider in this section only the case when $T$ consist of a single place $\fq$ where we simply write $\Theta_{S, \fq}(K/F, s)$ instead of $\Theta_{S, \{\fq\}}(K/F, s)$. Note that we have
\begin{equation}
\label{stickelberger2}
\Theta_{S, \fq}(K/F, s)\, =\, \sum_{\sigma\in G} \left(\zeta_S(\sigma,s) - \Norm(\fq)^{1-s} \zeta_S(\sigma \sigma_{\fq}^{-1}, s)\right) [\sigma^{-1}].
\end{equation} 

\paragraph{The Eisenstein class $\bEis_{\fb, \fc}$} We fix two non-trivial coprime ideals $\fb, \fc\subsetneq \cO_F$ with $2\not\in \fc$. Let $\Si$ be the set prime factor of $\fb\cdot \fc$ and put 
\begin{equation*}
\label{gammaf}
\Ga:= \cOscu :=\, \{\gamma\in \cOsu\mid \gamma \equiv 1 \mod \fc \cOs\}.
\end{equation*}
We consider the Eisenstein classes \eqref{eisenclass} defined in the previous section in the case $V=F$, $\cM = \fc \cOs$ and $\Ga$. Thus we have $\sL= \fc \cdot \cI^{\Si}$ (i.e.\ $\sL$ is the set of fractional ideal of the form $\fc \cdot \fa$ with $\fa$ coprime to $\Si$) and $\hT= \bA^{\Si}/\fc \cOs$. Put 
\begin{equation*}
\label{torcp}
C = T[\fb]:= \fb^{-1}\fc \cOs/\fc \cOs \qquad \mbox{and} \qquad t:= 1 + \fc \cOs\in F/\fc \cOs.
\end{equation*}
Note that $C$ is a finite $\Ga$-stable subset of $F/\fc\cOs$ and and that $t\in F/\fc\cOs$ is a $\Ga$ fixed-point with $t\not\in T[\fb]$.
We choose a coefficient ring $R\subseteq \bC$ such that $\Norm(\fb)$ is invertible in $R$. In this set-up the map
\eqref{eisenmap} will be denoted by
\begin{equation*}
\label{eisentotreal}
\bEis_{\fc} = \bEis_{\fc}(t): (R[T[\fb]]^0)^{\Ga}\lra H^{n-1}(\Ga, \cD_{\lpol}(1+\fc\cOs, R)(\ep)), \quad \alpha \mapsto \bEis_{\alpha, \fc}.
\end{equation*}
Here $\ep: F^* \to \{\pm 1\}$ is the character given by $\ep(x) = \sgn(\Norm_{F/\bQ}(x))$ for $x\in F^*$. We consider the following special choice for $\alpha$ (following \cite{bkl}) 
\begin{equation*}
\label{torsum}
\alpha[\fb]: = \, \Norm(\fb) \cdot 0 - \sum_{c\in T[\fb]} c
\end{equation*}
and define
\begin{equation}
\label{eisentotreal2}
\bEis_{\fb, \fc} :=\bEis_{\alpha[\fb], \fc}\in H^{n-1}(\Ga, \cD_{\lpol}(1+\fc\cOs, R)(\ep)).
\end{equation}
In order to establish the relation of the class \eqref{eisentotreal2} to partial zeta values in Theorem \ref{theorem:eisparzeta} below 
we will rephrase Prop.\ \ref{prop:eisencompare} in a special cases. Namely, if $\fm, \fa \subseteq \cO_F$ are coprime ideals that are also relatively prime to $\Si$ then we apply \ref{prop:eisencompare} to the element $t_0:=1 + \fc \fm \fa^{-1}$ of $T_{\fc\fm\fa^{-1}}=F_{\infty}/\fc\fm\fa^{-1}$ (it lies in the image of $\io_{t, \fc \fm \fa^{-1}}$ by Remark \ref{remark:imiotl}). 
We note that $E_{\fc}$ is the subgroup of $\Ga=\cOscu$ that stabilizes the lattice $\fc\fm\fa^{-1}\subseteq \cOs$ and that $E_{\fm \fc}$ is the stabilizer of $t_0\in T_{\fc\fm\fa^{-1}}$ in $E_{\fc}$. We denote the class \eqref{eisenbkl} for $L:= \fc \fm \fa^{-1}$, the group $\Gamma_0=E_{\fm \fc, +}$  and $\alpha = \alpha[\fb]$ by 
\begin{equation*}
\label{eisenbkltotreal}
\Eis_{\fc \fm \fa^{-1}, \fb}(t_0)\in H^{n-1}(E_{\fc\fm, +},\cD_{\pol}(1+\fc\fm\fa^{-1}, R)).
\end{equation*}
Consider the pair of essentially dual maps 
\begin{eqnarray}
\label{extzeroideal}
&& j: \Int(1 + \fc\fm\fa^{-1}, \bZ)\to\Int_{\loc, b}(1 + \fc \cOs, \bZ), \qquad  f\mapsto f_!   \\
&& j^{\vee} : \cD_{\lpol}^{\Sigma}(1+ \fc \cOs, R) \to \cD_{\pol}(1+\fc\fm\fa^{-1}, R), \qquad \mu \mapsto \mu|_{1 + \fc\fm\fa^{-1}}.
\label{distrestr}
\end{eqnarray} 
Together with the inclusion $E_{\fc\fm, +}\hra \Ga$ they induce homomorphisms
\begin{eqnarray}
\label{cores}
&& \cor := \cor_{E_{\fc\fm, +}}^{\Ga} \circ j_*: H_{n-1}(E_{\fc\fm, +}, \Int(1 + \fc\fm\fa^{-1}, \bZ)) \to  H_{n-1}(\Ga, \Int_{\loc, b}(1 + \fc\cOs, \bZ)(\ep)), \\
\label{resdual} 
&& \res:= (j^{\vee})_* \circ \res_{E_{\fc\fm, +}}^{\Gamma}: H^{n-1}(\Ga, \cD_{\lpol}(1+\fc\cOs, R)(\ep)) \to H^{n-1}(E_{\fc\fm, +},\cD_{\pol}(1+\fc\fm\fa^{-1}, R)).
\end{eqnarray}
By Prop.\ \ref{prop:eisencompare} we have 
\begin{equation}
\label{eisenres}
\res(\bEis_{\fb, \fc}) \, =\, \Eis_{\fc \fm \fa^{-1}, \fb}(t_0).
\end{equation}

\paragraph{Homology classes associated to ray classes} Let $A$ be another ring (later $A=\bZ$ or $A=\bZ[G]$ is the group ring of  an abelian Galois group). The pairing \eqref{eval} induces a cap-product pairing 
\begin{equation*}
\label{cap}
\cap: H^{n-1}(\Ga, \cD_{\lpol}(1+\fc\cOs, R)(\ep))\times H_{n-1}(\Ga, \Int_{\loc, b}(1 + \fc\cOs, A)(\ep))\, \lra \, A\otimes R.
\end{equation*}
Our aim is to show that the cap-product of the Eisenstein class \eqref{eisentotreal2} with certain homology classes associated canonically to ray classes yields the values of partial zeta functions at non-positive integers. To define these homology classes we follow mostly (\cite{dassp}, \S 3 and \S 5.3). 

Firstly, we present a general set-up to produce non-trivial homology classes in degree $n-1$. To fix ideas, let $\cU\subseteq U^{\Si}$ be a closed subgroup containing $U_{\infty}$ and let $M$ an abelian group. We consider the $(\bA^{\Si})^*$-module $C((\bA^{\Sigma})^*/\cU, M)$ of locally constant maps $\varphi: (\bA^{\Sigma})^*/\cU\to M$ and its submodule $C_c((\bA^{\Sigma})^*/\cU, M)$ consisting of those $\varphi$ that have compact support. The $(\bA^{\Si})^*$-action is defined by $(x \varphi)(a\cU) :=\varphi(x^{-1}a\cU)$ for $x\in (\bA^{\Si})^*$, $a\cU\in (\bA^{\Sigma})^*/\cU$ and $\varphi \in C((\bA^{\Sigma})^*/\cU,  M)$. 

Now assume that $M$ is an $\Ga$-module. We equip $C((\bA^{\Sigma})^*/\cU,  M)$ with an $\Ga$-action defined by $(\ga \varphi)(a\cU) := \ga \cdot \varphi(\ga^{-1}a\cU)$ for $\ga\in \Ga$, $\varphi \in C((\bA^{\Sigma})^*/\cU,  M)$ and $a\cU\in (\bA^{\Sigma})^*/\cU$. If $\cU_1 \subseteq \cU_2$ are closed subgroups of $U^{\Si}$ containing $U_{\infty}$ and if $\pr: (\bA^{\Sigma})^*/\cU_1\to (\bA^{\Sigma})^*/\cU_2$ denotes the projection then we can (and will) identify $C((\bA^{\Sigma})^*/\cU_2, M)$ with its image under the monomorphism $C((\bA^{\Sigma})^*/\cU_2, M)\to C((\bA^{\Sigma})^*/\cU_1, M),\varphi \mapsto \varphi\circ \pr$. With this convention we have $C_c((\bA^{\Sigma})^*/\cU_2, M)= C_c((\bA^{\Sigma})^*/\cU_1, M)\cap C((\bA^{\Sigma})^*/\cU_2, M)$. 

Suppose that $\cU_1, \cU_2, \cU_3\subseteq U^{\Sigma}$ are closed subgroups containing $U_{\infty}$ with $\cU_3\subseteq \cU_1\cap \cU_2$. Consider the $\Ga$-equivariant pairing 
\begin{equation}
\label{deltapair2}
C((\bA^{\Sigma})^*/\cU_1,  M) \times C((\bA^{\Sigma})^*/\cU_2, \bZ) \lra 
C((\bA^{\Sigma})^*/\cU_3, M), \qquad (\varphi, \psi) \mapsto \varphi\odot \psi
\end{equation}
defined by $(\varphi\odot \psi)(a\cU_3)= \psi(a\cU_1)\cdot \varphi(a\cU_2)$ for every $a\in (\bA^{\Sigma})^*$. Note that if $\psi$ has compact support then $\varphi\odot \psi$ has compact support as well, i.e.\ \eqref{deltapair2} restricts to a pairing
\begin{equation}
\label{deltapair2a}
C((\bA^{\Sigma})^*/\cU_1,  M) \times C_c((\bA^{\Sigma})^*/\cU_2, \bZ) \lra 
C_c((\bA^{\Sigma})^*/\cU_3, M).
\end{equation}
The latter induces a cap-product pairing 
\begin{equation}
\label{cap1}
\cap: H^i(\Ga, C((\bA^{\Sigma})^*/\cU_1, M)) \times H_j(\Ga, C_c((\bA^{\Sigma})^*/\cU_2, \bZ)) \, \lra  \, H_{j-i}(\Ga, C_c((\bA^{\Sigma})^*/\cU_3, M))
\end{equation} 
for every $i,j\in \bZ$. 

For an ideal $\fm\subseteq \cO_F$ coprime to $\Sigma$ there exists a canonical homology class (cf.\ \cite{dassp}, \S 3.1 \footnote{In loc.\ cit.\ we considered only the case $\fm = \cO_F$}) 
\begin{equation}
\label{theta}
\vartheta_{\fm}\in H_{n-1}(\Ga, C_c((\bA^{\Sigma})^*/U_{\fm}^{\Sigma}, \bZ)).
\end{equation} 
We recall its definition. By Dirichlet's unit theorem the homology group $H_{n-1}(E_{\fm \fc, +}, \bZ)$ is a free $\bZ$-module of rank one. Due to the chosen ordering of the embeddings $F\hookrightarrow \bR$ there is a canonical choice of a generator $\eta_{\fc\fm}\in H_{n-1}(E_{\fm \fc, +}, \bZ)$. Let $\cF\subseteq (\bA^{\Sigma})^*/U_{\fm}^{\Sigma}$ be a fundamental domain for the action of $\Ga/E_{\fm \fc, +}$ on $(\bA^{\Sigma})^*/U_{\fm}^{\Sigma}$. Since 
$C_c((\bA^{\Sigma})^*/U_{\fm}^{\Sigma}, \bZ)\, \cong \, \Ind^{\Ga}_{E_{\fm \fc, +}} C(\cF,\bZ)$ as $\Ga$-modules, by Shapiro's Lemma we have
\begin{equation*}
\label{shapiro}
H_{n-1}(\Ga, C_c((\bA^{\Sigma})^*/U_{\fm}^{\Sigma}, \bZ))\, \cong \, H_{n-1}(E_{\fm \fc, +}, C(\cF,\bZ)) \, \cong \, C(\cF,\bZ)\otimes H_{n-1}(E_{\fm \fc, +}, \bZ).
\end{equation*}
The homology class $\vartheta_{\fm}$ is the class that is mapped to $1_{\cF}\otimes\eta_{\fc\fm}$ under this isomorphism. 
If $\fm=\cO_F$ then we write $\vartheta$ instead of $\vartheta_{\fm}$. 

For $\cU_2= U^{\Si}$ and $\cU_1=\cU_3=\cU\subseteq U^{\Si}$ taking the cap-product with $\vartheta$ in \eqref{cap1} yields the map
\begin{equation*}
\label{cap2}
H^0(\Ga, C((\bA^{\Sigma})^*/\cU,  M))\, \lra  \, H_{n-1}(\Ga, C_c((\bA^{\Sigma})^*/\cU, M)), \quad \rho\mapsto \rho\cap \vartheta.
\end{equation*} 
Now consider the special case $A=\bZ=M$ (with trivial $\Ga$-action) and $\cU=U_{\fm}^{\Si}$ where $\fm\subseteq \cO_F$ is an ideal that is coprime to $\Si$. For $a\in (\bA^{\Sigma})^*/U_{\fm}^{\Si}$ let $\io_a: \bZ \to C_c((\bA^{\Sigma})^*/U_{\fm}^{\Si},  \bZ)$ be the homomorphism that maps $1$ to the characteristic function of the subset $\{a\}\subseteq (\bA^{\Sigma})^*/U_{\fm}^{\Si}$. Together with the inclusion $E_{\fc\fm, +}\hra \Ga$ it induces a homomorphism 
\begin{equation*}
\label{shapiro2}
\sh_a:= \cor_{E_{\fc\fm, +}}^{\Ga} \circ \, (\io_a)_*: H_{n-1}(E_{\fc\fm, +}, \bZ)\, \cong \, H_{n-1}(\Ga, C_c(\cA, \bZ))\stackrel{(j_!)*}{\lra} H_{n-1}(\Ga, C_c((\bA^{\Sigma})^*/U^{\Sigma}, \bZ)).
\end{equation*}
where  $\cA: = \Ga \cdot a\subseteq (\bA^{\Sigma})^*/U_{\fm}^{\Si}$ is the $\Ga$-orbit of $a$ and $j:\cA \hra (\bA^{\Sigma})^*/U^{\Sigma}$ is the inclusion. Note that if $a_1, \ldots, a_h$ is a system of representatives of the $\Ga=\cOscu$-orbits in $(\bA^{\Sigma})^*/U_{\fm}^{\Sigma}$ then we have 
\begin{equation}
\label{fundclass}
\vartheta_{\fm} \, = \, \sum_{i=1}^h \sh_{a_i}(\eta_{\fc\fm}).
\end{equation}

\begin{lemma}
\label{lemma:thetacap1}
Let $1_{\cA}\in H^0(\Ga, C((\bA^{\Sigma})^*/U_{\fm}^{\Si},  \bZ))$ be the characteristic function of $\cA$. We have $1_{\cA} \cap \vartheta = \sh_a(\eta_{\fc\fm})$.
\end{lemma}

\begin{proof} Let $\io: C_c((\bA^{\Sigma})^*/U^{\Sigma}, \bZ)\, \to C_c((\bA^{\Sigma})^*/U_{\fm}^{\Sigma}, \bZ)$ be the inclusion. It can be easily seen (by using $\cor_{E_{\fc\fm, +}}^{E_{\fc, +}}(\eta_{\fc}) = \eta_{\fc\fm}$) that the induced homomorphism 
\begin{equation*}
\label{prstar}
\io_*: H_{n-1}(\Ga, C_c((\bA^{\Sigma})^*/U^{\Sigma}, \bZ))\, \lra \, H_{n-1}(\Ga, C_c((\bA^{\Sigma})^*/U_{\fm}^{\Sigma}, \bZ))
\end{equation*}
maps $\vartheta$ to $\vartheta_{\fm}$. If we denote the cap-product \eqref{cap1} for $\cU_1= \cU_2= \cU_3= U_{\fm}^{\Si}$ by $\cap'$ then together with \eqref{fundclass} we obtain
\[
1_{\cA} \cap \vartheta  \,=\, 1_{\cA} \cap' \io_*(\vartheta) \,= \,1_{\cA} \cap' \vartheta_{\fm}\,=\, \sum_{i=1}^h 1_{\cA} \cap' \sh_{a_i}(\eta_{\fc\fm}).
\]
The assertion now follows from $1_{\cA} \cap' \sh_{a_i}(\eta_{\fc\fm})= \sh_{a_i}(\eta_{\fc\fm})=\sh_{a}(\eta_{\fc\fm})$ if $a_i\in \cA$ and $1_{\cA} \cap'\sh_{a_i}(\eta_{\fc\fm})=0$ if $a_i\not\in \cA$.
\end{proof}

Let $S$ be a finite set of nonarchimedean places of $F$ disjoint from $\Si$. In order to define homology classes in $H_{n-1}(\Ga, \Int_{\loc, b}(1 + \fc\cOs, \bZ))$ that are related to partial zeta values we recall the definition of the $(\bA_f^{\Si})^*$-equivariant homomorphism
\begin{equation}
\label{dasspmap1}
\Delta_{S, f}^{\Sigma}: C_c( (\bA_f^{\Sigma})^*/ U_f^{S, \Sigma}, \bZ)\,=\, C_c( F_S^*\times (\bA_f^{S, \Sigma})^*/ U_f^{S, \Sigma}, \bZ) \,\lra \, C_c(\bA_f^{\Sigma}, \bZ)
\end{equation}
introduced in (\cite{dassp}, \S 5.3). 
Since $(\bA_f^{S, \Sigma})^*/U_f^{S, \Sigma}$ is canonically isomorphic to the group of fractional ideals $\cI^{S, \Sigma}:= \cI^{S\cup \Sigma}$ that are coprime to $S\cup \Si$, we can identify the space $(\bA_f^{\Sigma})^*/ U_f^{S, \Sigma}$ with $F_S^*\times \cI^{S, \Sigma}$. There exists a canonical isomorphism (see \cite{dassp}, \S 2 or \cite{dashon}, Prop.\ 5.3)
\begin{equation*}
\label{adelicfunctions2}
C_c(F_S^*, \bZ) \otimes \bZ[\cI^{S, \Sigma}]\, \cong\, C_c(F_S^*, \bZ) \otimes C_c(\cI^{S, \Sigma}, \bZ)\,\stackrel{\cong}{\lra} \, C_c( F_S^*\times \cI^{S, \Sigma}, \bZ)\,= \,C_c( (\bA_f^{\Sigma})^*/ U_f^{S, \Sigma}, \bZ).
\end{equation*}
Hence we can (and will) identify the source of  \eqref{dasspmap1} with the module $C_c(F_S^*, \bZ) \otimes \bZ[\cI^{S, \Sigma}]$. Define
\begin{eqnarray}
\label{dasspmap2}
&& \delta_S: C_c(F_S^*, \bZ) \,\lra \,C_c(F_S, \bZ) , \qquad f \mapsto f_!\\
&& \delta_f^{S, \Sigma}: \bZ[\cI^{S, \Sigma}]\,\lra \, C_c(\bA_f^{S, \Sigma}, \bZ), \quad \sum_{\fa \in \cI^{S, \Sigma}} m_{\fa} \, [\fa]  \mapsto \sum_{\fa \in \cI^{S, \Sigma}} m_{\fa} 1_{\widehat{\fa}^{S, \Sigma}}
\label{dasspmap2a}
\end{eqnarray} 
where $1_{\widehat{\fa}^{S, \Sigma}}$ is the characteristic function of $\widehat{\fa}^{S, \Sigma} := \widehat{\fa}^{S\cup \Sigma} \subseteq \bA_f^{S, \Sigma}$. Furthermore let
\begin{equation}
\label{adelicfunctions3}
C_c(F_S, \bZ) \otimes C_c(\bA_f^{S, \Sigma}, \bZ)\,\lra \, C_c(\bA_f^{\Sigma}, \bZ), \quad \varphi_S \otimes \varphi^S\mapsto (\varphi_S \circ \pr_1) \cdot (\varphi^S\circ \pr_2)
\end{equation}
where $\pr_1: (\bA_f^{\Sigma})\to F_S$, $\pr_2: \bA_f^{\Sigma}\to \bA_f^{S, \Sigma}$ denote the projections. The map \eqref{dasspmap1} is defined
as the composite of $\delta_S \otimes  \delta_f^{S, \Sigma}$ with \eqref{adelicfunctions3}. 

Note that the restriction of a locally constant map $\varphi: \bA_f^{\Si} \to \bZ$ with compact support to the subset $1 + \fc\cOs\subseteq \bA_f^{\Sigma}$ is contained in $\Int_{\loc, b}(1+ \fc\cOs, \bZ)$ (in fact if we denote the restriction by $f$ then there exists a fractional ideal $\fa\in \cI^{\Si}$ such that $f|_{x+ \fc\fa}$ is constant for every $x\in 1 + \fc\cOs$ and non-zero for only finitely many cosets $x+\fc\fa$). Therefore 
\begin{equation}
\label{dasspmap2b}
C_c(\bA_f^{\Sigma}, \bZ)\, \lra \, \Int_{\loc, b}(1 + \fc\cOs, \bZ), \qquad \varphi \mapsto \varphi|_{1 + \fc\cOs}
\end{equation}
is a well-defined $\Ga$-equivariant homomorphism. By abuse of notation we denote the composite of \eqref{dasspmap1} and \eqref{dasspmap2b} by
\begin{equation}
\label{dasspmap1a}
\Delta_{S, f}^{\Sigma}: C_c( (\bA_f^{\Sigma})^*/ U_f^{S, \Sigma}, \bZ) \,\lra \, \Int_{\loc, b}(1 + \fc\cOs, \bZ)
\end{equation}
as well. Finally, we incorporate the archimedean places as well, i.e.\ we extend \eqref{dasspmap1a} to an $\Ga$-equivariant homomorphism
\begin{equation}
\label{dasspmap3}
\Delta_S^{\Sigma}: C_c( (\bA^{\Sigma})^*/ U^{S, \Sigma} , \bZ) \,\lra \, \Int_{\loc, b}(1 + \fc\cOs, \bZ)(\ep)
\end{equation}
as follows. For a function $\varphi: (\bA^{\Sigma})^*/ U^{S, \Sigma} = (\bA_f^{\Sigma})^*/ U_f^{S, \Sigma} \times F_{\infty}^*/U_{\infty} \to \bZ$ with compact support and $x_{\infty}\in F_{\infty}^*/U_{\infty}$ we let $\varphi(\wcdot, x_{\infty})$ denote the map $(\bA_f^{\Sigma})^*/ U_f^{S, \Sigma}\to \bZ, x \mapsto \varphi(x, x_{\infty})$. We define \eqref{dasspmap3} by 
\begin{equation*}
\label{dasspmap3a}
\Delta_S^{\Sigma}(\varphi)\, =\, \sum_{x_{\infty}U_{\infty} \in F_{\infty}^*/U_{\infty}} \ep(x_{\infty}) \cdot \Delta_{S, f}^{\Sigma}(\varphi(\wcdot, x_{\infty}U_{\infty})).
\end{equation*}
More generally, if $\cU\subseteq U^{\Si}$ is a closed subgroup containing $U^{S, \Sigma}$ and $M$ is an abelian group then we define
\begin{equation}
\label{dasspmap4}
\Delta_S^{\Sigma}: C_c( (\bA^{\Sigma})^*/ \cU , M) \,\lra \, \Int_{\loc, b}(1 + \fc\cOs, M)(\ep)
\end{equation}
as follows. Note that we can identify the source with the module $C_c((\bA^{\Sigma})^*/\cU, \bZ)\otimes M$ and that $C_c( (\bA^{\Sigma})^*/ \cU, \bZ)$ can be viewed naturally as a submodule of $C_c( (\bA^{\Sigma})^*/ U^{S, \Sigma} , \bZ)$. Therefore we can denote elements of the source of  \eqref{dasspmap4} as finite sums $\sum_{i=1}^r \psi_i \otimes m_i$ with $\psi_1, \ldots, \psi_r\in C_c((\bA^{\Sigma})^*/\cU, \bZ)$ and $m_1, \ldots, m_r\in M$ and define 
\begin{equation}
\label{dasspmap4a}
\Delta_S^{\Sigma}\left(\sum_{i=1}^r \varphi_i \otimes f_i\right):=\, \sum_{i=1}^r\Delta_S^{\Sigma}(\varphi_i) \cdot m_i.
\end{equation}
Note that if $M$ is an $\Ga$-module then \eqref{dasspmap4} is $\Ga$-equivariant.

For $M = \Int_{\loc, b}(1 + \fc\cOs, A)$ composing \eqref{dasspmap4} with the map
\begin{eqnarray}
\label{dasspmap4b}
&& \Int_{\loc, b}(1 + \fc\cOs, \Int_{\loc, b}(1 + \fc\cOs, A)) = \Int_{\loc, b}(1 + \fc\cOs, \bZ)\otimes \Int_{\loc, b}(1 + \fc\cOs, A)\\
&&\hspace{2cm}\lra  \Int_{\loc, b}(1 + \fc\cOs, A), \quad f_1 \otimes f_2 \mapsto f_1 \cdot f_2 \nonumber
\end{eqnarray}
induces an $\Ga$-equivariant homomorphism
\begin{equation}
\label{dasspmap5}
\wDelta_S^{\Sigma}: C_c((\bA^{\Sigma})^*/\cU,  \Int_{\loc, b}(1 + \fc\cOs, A)) \, \lra \, \Int_{\loc, b}(1 + \fc\cOs, A)(\ep).
\end{equation}

\begin{remark}
\label{remark:changeS}
\rm We recall the dependence of the map \eqref{dasspmap4a} on the set $S$ (compare \cite{dassp}, Remark 5.5). Let $v\in S$, put $S'= S\setminus \{v\}$ and assume that $\cU$ contains $U^{S', \Sigma}$. Then we have 
\begin{equation*}
\label{dschangeS}
\Delta_S^{\Sigma}(\varphi) \, =\, \Delta_{S'}^{\Sigma}(\varphi - [\varpi_v] \cdot \varphi)
\end{equation*}
for every $\varphi\in C_c( (\bA^{\Sigma})^*/ \cU , M)$. Here $\varpi_v\in F_v^*$ is a uniformizer and $[\varpi_v]\in(\bA^{\Sigma})^*$ denotes the adele whose component at $v$ is a uniformizer $\varpi_v\in F_v^*$ and whose other components are $=1$.
Similarly, for the map \eqref{dasspmap5} we obtain
\begin{equation}
\label{wdschangeS}
\wDelta_S^{\Sigma}(\varphi) \, =\, \wDelta_{S'}^{\Sigma}(\varphi - [\varpi_v] \cdot \varphi)
\end{equation}
for every $\varphi\in C_c( (\bA^{\Sigma})^*/ \cU , \Int_{\loc, b}(1 + \fc\cOs, A))$.
\end{remark}

For open subgroups $\cU_1, \cU_2\subseteq U^{\Si}$ with $U^{S, \Sigma}\subseteq \cU_1, \cU_2$ we consider the $\Ga$-equivariant pairing 
\begin{equation}
\label{deltapair}
\langle\wcdot , \wcdot \rangle_S: C((\bA^{\Sigma})^*/\cU_1,  \Int_{\loc, b}(1 + \fc\cOs, A)) \times C_c((\bA^{\Sigma})^*/\cU_2, \bZ) \lra \Int_{\loc, b}(1 + \fc\cOs, A)(\ep)
\end{equation}
defined as the composition of \eqref{deltapair2a} (for $\cU_3 := \cU_1 \cap \cU_2$ and $M= \Int_{\loc, b}(1 + \fc\cOs, A)$) with 
the homomorphism $\wDelta_S^{\Sigma}$ (for $\cU= \cU_3$), i.e.\ we have 
\begin{equation*}
\label{deltapair1}
\langle \varphi, \psi\rangle_S \, =\, \wDelta_S(\varphi\odot \psi)
\end{equation*}
for $\varphi \in C((\bA^{\Sigma})^*/\cU_1,  \Int_{\loc, b}(1 + \fc\cOs, A))$ and $\psi\in C_c((\bA^{\Sigma})^*/\cU_2, \bZ)$. 
The pairing \eqref{deltapair} induces a cap-product pairing 
\begin{eqnarray}
\label{cap3}
&& \cap: H^i(\Ga, C((\bA^{\Sigma})^*/\cU_1,  \Int_{\loc, b}(1 + \fc\cOs, A))) \times H_j(\Ga, C_c((\bA^{\Sigma})^*/\cU_2, \bZ))\hspace{2cm}\\
&& \hspace{4cm} \, \lra  \quad H_{j-i}(\Ga, \Int_{\loc, b}(1 + \fc\cOs, A)(\ep))
\nonumber
\end{eqnarray} 
for every $i,j\in \bZ$. 

In particular for $\cU_1=\cU_2= U_{\fm}^{\Si} \subseteq U^{\Si}$ and $i=0$ taking the cap-product with the homology class $\sh_a(\eta_{\fc\fm})$ associated to a point $a\in (\bA^{\Sigma})^*/U_{\fm}^{\Si}$ yields a map 
\begin{equation*}
\label{cap4}
H^0(\Ga, C((\bA^{\Sigma})^*/U_{\fm}^{\Si},  \Int_{\loc, b}(1 + \fc\cOs, A)))\, \lra  \, H_{n-1}(\Ga, \Int_{\loc, b}(1 + \fc\cOs, A)(\ep)), \quad \varphi\mapsto \varphi\cap \,\sh_a(\eta_{\fc\fm}).
\end{equation*} 

\begin{lemma}
\label{lemma:rescor}
For $\varphi\in H^0(\Ga, C((\bA^{\Sigma})^*/U_{\fm}^{\Si},  \Int_{\loc, b}(1 + \fc\cOs, A)))$ we have 
\[
\varphi \cap \,\sh_a(\eta_{\fc\fm}) \, = \, \cor_{E_{\fc\fm, +}}^{\Ga} (\varphi(a) \cap \eta_{\fc\fm}).
\]
\end{lemma}

\begin{proof} Let $\ev_a : C((\bA^{\Sigma})^*/U_{\fm}^{\Si}, \Int_{\loc, b}(1 + \fc\cOs, \bZ)) \to  \Int_{\loc, b}(1 + \fc\cOs, \bZ),\varphi\mapsto \varphi(a)$ denote the evaluation map at $a$. We have 
\[
\langle \varphi , \io_a(m) \rangle \, =\, m \cdot \ev_a(\varphi)
\]
for every $\varphi\in C((\bA^{\Sigma})^*/U_{\fm}^{\Si}, \Int_{\loc, b}(1 + \fc\cOs, \bZ))$ and $m\in \bZ$. Standard functorial properties of the cap-product with respect to restrictions and corestrictions therefore imply 
\[
\cor_{E_{\fc\fm, +}}^{\Ga}\left(\left((\ev_a)_* \circ \res_{E_{\fc\fm, +}}^{\Ga}( \kappa )\right) \cap \eta_{\fc\fm}\right) \, =\, \kappa\cap \sh_a(\eta_{\fc\fm})
\]
for every $\kappa \in H^i(\Ga, C((\bA^{\Sigma})^*/U_{\fm}^{\Si},  \Int_{\loc, b}(1 + \fc\cOs, A)))$.
\end{proof}

\begin{remarks}
\label{remarks:pairchangeS}
\rm (a) For the applications in section \ref{section:lvalues} we need to address the dependence of the pairing \eqref{deltapair} on the set $S$. Let $v\in S$, put $S'= S\setminus \{v\}$ and assume that $U^{S', \Sigma}\subseteq \cU_3= \cU_1 \cap \cU_2$. By \eqref{wdschangeS} we have 
\begin{equation}
\label{deltapair3}
\langle \varphi, \psi\rangle_S \, =\,\langle  [\varpi_v]\cdot \varphi, \psi - [\varpi_v] \cdot \psi \rangle_{S'} + \langle \varphi - [\varpi_v] \cdot \varphi,  \psi \rangle_{S'}
\end{equation}
for every $\varphi \in C((\bA^{\Sigma})^*/\cU_1,  \Int_{\loc, b}(1 + \fc\cOs, A))$ and $\psi\in C_c((\bA^{\Sigma})^*/\cU_2, \bZ)$.
\medskip

\noi (b) We also need a more concrete description of the pairing \eqref{deltapair} in the case $\cU : =\cU_1\subseteq \cU_2= U^{\Si}$. For that we fix an open subgroup $V\subseteq U_S$ such that $V \times U^{S, \Si}\subseteq \cU$. For that we write elements of $(\bA^{\Sigma})^*$ as triples $a=(a_1, a_2, a_{\infty})$ with $a_1=(a_v)_{v\in S}\in F_S^*$, $a_2= (a_v)_{v\not\in S\cup \Si, v\nmid \infty} (\bA_f^{S, \Sigma})^*$ and $a_{\infty}\in F_{\infty}^*$. Let $V\subseteq U_S$ be an open subgroup with $V \times U^{S, \Si}\subseteq \cU$. We have 
\begin{equation*}
\label{fxas}
\Delta_S^{\Si}(1_{a (V \times U^{S, \Si})}) \, =\, 1_{\fX_{a, V, S}}
\end{equation*}
where $\fX_{a, V, S}$ denotes the set
\begin{equation*}
\label{fxas2}
\fX_{a, V, S} = \left\{x\in 1+ \fc \cO_{\Si}\mid x\in a_1 V\,  \text{and} \,\ord_v(x) \ge \ord_v(a_v)\, \forall \, v\not \in S\cup \Si \cup S_{\infty} \right\}
\end{equation*}
Note that if $\fa\in \cI^{\Si}$ denotes the fractional ideal associated to the idele $(a_1, a_2)\in (\bA_f^{\Sigma})^*$ then we have  $\fX_{a, V, S}\subseteq \fa\cap 1+ \fc \cO_{\Si}$. Since 
\[
\varphi\odot \psi = \sum_{a (V \times U^{S, \Si}) \in (\bA^{\Sigma})^*/(V \times U^{S, \Si})} \psi(a U^{\Si}) 1_{a (V \times U^{S, \Si})} \otimes \varphi(a\cU)
\] 
we get 
\begin{equation}
\label{deltapair2b}
\langle \varphi, \psi\rangle_S \, =\,\sum_{a (V \times U^{S, \Si}) \in (\bA^{\Sigma})^*/(V \times U^{S, \Si})}  \psi(a U^{\Si})\, 1_{\fX_{a, V, S}} \cdot \varphi(a \cU)
\end{equation} 
\enddemo
\end{remarks}

Now we consider \eqref{cap3} for $\cU_1=\cU\subseteq \cU_2=U^{\Sigma}$. Taking the cap-product with $\vartheta$ yields a map
\begin{equation}
\label{cap5}
H^0(\Ga, C((\bA^{\Sigma})^*/\cU,  \Int_{\loc, b}(1 + \fc\cOs, A)))\, \lra  \, H_{n-1}(\Ga, \Int_{\loc, b}(1 + \fc\cOs, A)(\ep)), \quad \rho\mapsto \rho\cap \vartheta.
\end{equation} 
We define certain canonical elements in the source of \eqref{cap5}. For that put $N:= \Norm_{F/\bQ}: F\to \bQ$ and let $\bN=\bN^{\Sigma}: (\bA^{\Sigma})^*\to \bQ^*$ be the {\it idele norm character} given by $\bN((a_v)_v)= \prod_{v\not\in \Sigma} \bN_v(a_v)$ where 
\begin{equation*}
\bN_v(a) \,=\,  \left\{\begin{array}{ll} |a|_v^{-1} &  \mbox{if $v$ is nonarchimedean,}\\
\sgn(a) &  \mbox{if $v$ is archimedean.}
\end{array}\right.
\end{equation*}
For an idele $a=(a_v)_{v\not\in\Si}\in (\bA^{\Si})^*$ let $\fa= \left\{x\in \cO_{\Si}\mid \ord_v(x) \le \ord_v(a_v) \text{for every $v\not \in \Si$, $v\nmid \infty$}\, \right\}$ be the associated fractional ideal. Note that $\bN(a) = \pm \Norm(\fa)$ and $\bN(a) = \Norm(\fa)$ if $a_v>0$ for all $v\nmid \infty$. Moreover if $a=\ga \in \cOsu\subseteq (\bA^{\Si})^*$ then we have $\bN(a) = N(\ga)$. Consider the map 
\begin{equation}
\label{normfunc}
\cN: (\bA^{\Sigma})^*/U^{\Sigma}\lra  \Int_{\loc, b}(1 + \fc\cOs, \bZ), \quad a U^{\Si} \mapsto \bN(a)^{-1} 1_{\fa \cap 1+\fc\cOs} \cdot N.
\end{equation}
If we identify the group $(\bA^{\Sigma})^*/U^{\Sigma}$ with $\cI^{\Si}\times F_{\infty}^*/U_{\infty}$ and correspondingly write its elements as pairs $(\fa, a_{\infty} U_{\infty})$ then we have 
\begin{equation}
\label{normfunc1}
\cN(\fa, a_{\infty} U_{\infty}) \, =\, \ep(a_{\infty})\,\Norm(\fa)^{-1} 1_{\fa \cap 1+\fc\cOs} \cdot N
\end{equation}
where $\ep(a_{\infty})= \prod_{v\mid \infty} \sgn(a_v)$ for $a_{\infty}= (a_v )_{v\mid \infty}\in \prod_{v\mid\infty} F_v^*$.

Since $\Norm(\fa)^{-1} N(x) \in \bZ$ for every $x\in \fa$ we see that the right hand side of \eqref{normfunc1} lies indeed in $\Int_{\loc, b}(1 + \fc\cOs, \bZ)$. The map \eqref{normfunc} is $\Ga$-equivariant because of
\begin{equation}
\label{normfunc2}
\cN(\ga \cdot a U^{\Si}) =  \bN(\ga a)^{-1} \cdot 1_{\ga\fa \cap 1+ \fc\cOs} \cdot N =  \bN(a)^{-1} \ga \cdot (1_{\fa \cap 1+ \fc\cOs} \cdot N)=  \ga \cN(a U^{\Si})
\end{equation}
for every $\ga\in \Ga$ and $aU^{\Si}\in (\bA^{\Sigma})^*/U^{\Sigma}$.
It follows that we can view $\cN^k$ for $k\in \bZ_{\ge 0}$ as an element of $H^0(\Ga, C((\bA^{\Sigma})^*/U^{\Sigma}, \Int_{\loc, b}(1 + \fc\cOs, \bZ)))$. We remark that for $k=0$ the map $\cN^0:  (\bA^{\Sigma})^*/U^{\Sigma}\to  \Int_{\loc, b}(1 + \fc\cOs, \bZ)$ is given by $\cN^0( a U^{\Si}) = 1_{\fa \cap 1+\fc\cOs}$. 

Let $\fm\subseteq \cO_F$ be again an ideal coprime to $\Si$. For a ray class $\fA\in \cI^{\fc\fm}/\cP^{\fc\fm}$ and $k\in \bZ_{\ge 0}$ we introduce a certain homology class 
\begin{equation}
\label{homclassray}
\varrho_{\fA}^k \, =\, \varrho_{\fA, S, \fc}^k \in H_{n-1}(\Ga, \Int_{\loc, b}(1 + \fc\cOs, \bZ)(\ep))
\end{equation}
as follows. We identify the ray class group $\cI^{\fc\fm}/\cP^{\fc\fm}$ with the idele class group $(\bA^{\Sigma})^*/U_{\fm}^{\Sigma}\Ga$, so that  we can (and will) view $\fA$ as an $\Ga$-orbit in $(\bA^{\Sigma})^*/U_{\fm}^{\Sigma}$. The class \eqref{homclassray} is given by 
\begin{equation*}
\label{homclassray2}
\varrho_{\fA}^k \,:= \, (\cN^k\cup 1_{\fA}) \cap \vartheta \in H_{n-1}(\Ga, \Int_{\loc, b}(1 + \fc\cOs, \bZ)(\ep))
\end{equation*}
where $1_{\fA}\in H^0(\Ga, C((\bA^{\Sigma})^*/U_{\fm}^{\Sigma}, \bZ))$ denotes the characteristic function of $\fA$ and where we view
$\cN^k$ for $k\in \bZ_{\ge 0}$ as an element of $H^0(\Ga, C((\bA^{\Sigma})^*/U^{\Sigma}, \Int_{\loc, b}(1 + \fc\cOs, \bZ)))$.\footnote{Note that for $k=0$ the map $\cN^0:  (\bA^{\Sigma})^*/U^{\Sigma}\to  \Int_{\loc, b}(1 + \fc\cOs, \bZ)$ is given by $\cN^0( a U^{\Si}) = 1_{\fa \cap 1+\fc\cOs}$ where $\fa$ denotes again the fractional ideal associated to $aU^{\Sigma}$.} Moreover the cup-product $\cN^k\cup 1_{\fA}$ is induced by the pairing \eqref{deltapair2} for $M= \Int_{\loc, b}(1 + \fc\cOs, \bZ)$. 

\begin{theorem}
\label{theorem:eisparzeta}
We have 
\begin{equation}
\label{eisparzeta}
\bEis_{\fb, \fc} \cap \varrho_{\fA}^k  \, =\, (-1)^{n-1} \Norm(\fb) \, \zeta(\fc\fm, \fA, -k) + (-1)^n \Norm(\fb)^{-k} \, \zeta(\fc\fm, \fA\fB, -k)
\end{equation}
where $\fB=[\fb]\in \cI^{\fc\fm}/\cP^{\fc\fm}$ denotes the ray class of $\fb$. 
\end{theorem}

For the proof we choose a fractional ideal $\fa\subseteq \cO_F$ coprime to $S\cup \Si$ such that $\fa^{-1}$ is a representative of $\fA$. If we identify 
$(\bA^{\Sigma})^*/U_{\fm}^{\Sigma}$ with the product $F_S^*/U_{\fm, S} \times \cI^{S, \Sigma} \times F_{\infty}^*/U_{\infty}$ then $\fA$ is the $\Ga/E_{\fc\fm, +}$-orbit of the element 
\begin{equation*}
\label{nicecf}
a:=(1,  \fa^{-1}, 1)\in F_S^*/U_{\fm, S} \times \cI^{S, \Sigma} \times F_{\infty}^*/U_{\infty}\, =\, (\bA^{\Sigma})^*/U_{\fm}^{\Sigma}.
\end{equation*}
Let $\cN_0\in \Int(1 + \fc\fm\fa^{-1}, \bZ)$ denote the polynomial function $\cN_0(x) := \Norm(\fa)\Norm(x)$ for $x\in 1 + \fc \fm \fa^{-1}$.  
Since $N(\vep)=1$ for every $\vep\in E_{\fc\fm, +}$ we have $\cN_0\in H^0(E_{\fc\fm, +}, \Int(1 + \fc\fm\fa, \bZ))$. 

\begin{lemma}
\label{lemma:shapirohom}
We have $\varrho_{\fA}^k = \cor(\cN_0^k\cap \eta_{\fc\fm})$.\footnote{Recall that $\cor$ has been defined in \eqref{cores}.}
\end{lemma}

\begin{proof} Note that $\cN^k(a) = j(\cN^k_0)$. Together with Lemmas \ref{lemma:thetacap1} and \ref{lemma:rescor} we obtain 
\begin{eqnarray*}
\varrho_{\fA}^k & = & (\cN^k\cup 1_{\fA}) \cap \vartheta \, =\, \cN^k \cap (1_{\fA} \cap \vartheta) \, =\, \cN^k \cap\, \sh_a(\eta_{\fc\fm})\, =\, 
\cor_{E_{\fc\fm, +}}^{\Ga} (\cN(a)^k \cap \eta_{\fc\fm})\\
& = & \cor_{E_{\fc\fm, +}}^{\Ga} (j_*(\cN_0^k) \cap \eta_{\fc\fm})\, =\, \cor(\cN_0^k\cap \eta_{\fc\fm}).
\end{eqnarray*}
\end{proof}

\begin{proof}[Proof of Theorem \ref{theorem:eisparzeta}] By \eqref{eisenres} and Lemma \ref{lemma:shapirohom} we have 
\begin{eqnarray*}
\bEis_{\fb, \fc} \cap \varrho_{\fA}^k & = & \bEis_{\fb, \fc} \cap \cor(\cN_0^k\cap \eta_{\fc\fm})\, =\, \res(\bEis_{\fb, \fc}) \cap (\cN_0^k\cap \eta_{\fc\fm})\\
& = & \Eis_{\fc \fm \fa^{-1}, \fb}(t_0)\cap (\cN_0^k\cap \eta_{\fc\fm})\, =\, \cN_0^k\cap (\Eis_{\fc \fm \fa^{-1}, \fb}(t_0)\cap \eta_{\fc\fm}).
\end{eqnarray*}
To finish the proof we recall the relation between $\Eis_{\fc \fm \fa^{-1}, \fb}(t_0)\cap \eta_{\fc\fm}\in H_0(E_{\fc\fm, +},\cD_{\pol}(1+\fc\fm\fa^{-1}, R))$ and values partial zeta functions from (\cite{bkl}, \S 5). For that we enlarge $R$, i.e.\ we pass to the complex numbers $R=\bC$. 

To rephrase the result of Beilinson, Kings and Levin suitable for our framework we fix an ideal $\fn\subseteq \cO_F$ coprime to $\fc\fm$ (in the formula \eqref{bkliso} below $\fn$ will be either $\fa$ or $\fa\fb$). Put $L = \fc \fm \fn^{-1}$, $H= 1+ L$, $\Ga_0=E_{\fc\fm, +}$ and $t_0= t+L\in F_{\infty}/L$. As in (\cite{bkl}, \S 5) we consider the following composition of isomorphisms 
\begin{equation}
\label{bkliso}
\Phi_{\fn}: \bC[\![w]\!]\, \stackrel{\cong}{\lra}\, \cD_{\pol}(L, \bC)_{\Ga_0}  \,  \stackrel{\cong}{\lra}\, \cD_{\pol, L}(1+ L, \bC)_{\Ga_0}.
\end{equation}
The first isomorphism is given as follows. We view the embeddings $\xi_1, \ldots, \xi_n: F\to \bR$ as polynomial functions on $L$, i.e.\ we have $(\xi_i)|_L\in\Int(L, \bC)$. In fact the powers $\uxi^{\um}$, $\um\in (\bZ_{\ge 0})^n$ form a $\bC$-basis of $\Int(L, \bC)$. As in Example \ref{example:distchar0} for $i=1, \ldots, n$ let $z_i \in \cD_{\pol}(L, \bC)$ be given by $z_i(\uxi^{\um}) = 1$ if $\um = e_i$ and $z_i(\uxi^{\um}) = 0$ if $\um\ne e_i$. If we put $w:= z_1\cdot \ldots \cdot z_n$ then we have $\cD_{\pol}(L, \bC)^{\Ga_0} = \bC[\![w]\!]$ (see \cite{bkl}, Lemma 5.4). Now the isomorphism $\bC[\![w]\!]\cong \cD_{\pol}(L, \bC)_{\Ga_0}$ is induced by the projection $\cD_{\pol}(L, \bC)\to \cD_{\pol}(L, \bC)_{\Ga_0}$. For the second isomorphism in \eqref{bkliso} note that a polynomial function $f: L \to \bC$ extends uniquely to a polynomial function $\tf: \bZ+ L\to \bC$. We obtain an isomorphism $\Int(L, \bC) \to \Int_L(1+ L, \bC), f\mapsto \tf|_{1 + L}$, hence dually an isomorphism 
$\cD_{\pol, L}(1+ L, \bC)\to \cD_{\pol}(L, \bC)$. Passing to $\Ga_0$-coinvariants yields the isomorphism $\cD_{\pol}(L, \bC)_{\Ga_0}  \cong \cD_{\pol, L}(1+ L, \bC)_{\Ga_0}$. 

By (\cite{bkl}, 3.44, 5.6 and (5.10)) we have 
\begin{eqnarray}
\label{bkliso2}
\Eis_{\fc \fm \fa^{-1}, \fb}(t_0)\cap \eta_{\fc\fm}& =& (-1)^{n-1}\Norm(\fb)\,\Phi_{\fa} \left(\sum_{j=0}^{\infty} \Norm(\fa)^{-j} \zeta(\fc\fm, \fA, -j) \frac{w^j}{(j!)^n}\right) \\
& & + (-1)^n \,\Phi_{\fa\fb} \left(\sum_{j=0}^{\infty} \Norm(\fa\fb)^{-j} \zeta(\fc\fm, \fA\fB, -k) \frac{w^j}{(j!)^n}\right).
\nonumber
\end{eqnarray}
Note that the second summand is initially an element of $\cD_{\pol}(1+ \fc\fm(\fa\fb)^{-1}, \bC)_{\Ga_0}$. However since we work with complex coefficients we can identify it canonically with the group $\cD_{\pol}(1+ \fc\fm\fa^{-1}, \bC)_{\Ga_0}$. 

To finish the proof note that by \eqref{normdist} we have 
\[
\cN_0^k\cap \Phi_{\fa}(w^j) \, =\, \cN_0^k\cap \Phi_{\fa\fb}(w^j)\, =\, \left\{ \begin{array}{cc} \Norm(\fa)^k (k!)^n & \mbox{if $k = j$,}\\
                                                      0 & \mbox{otherwise.}
                                                      \end{array}\right.
\]
We conclude
\begin{eqnarray*}
\bEis_{\fb, \fc} \cap \varrho_{\fA}^k & = & \cN_0^k\cap (\Eis_{\fc \fm \fa^{-1}, \fb}(t_0)\cap \eta_{\fc\fm})\\
& = & (-1)^{n-1} \Norm(\fb) \zeta(\fc\fm, \fA, -k) + (-1)^n \Norm(\fb)^{-k} \zeta(\fc\fm, \fA\fB, -k).
\end{eqnarray*}
\end{proof}

Now we will give formulas for the values of the partial zeta functions $\zeta_S(\sigma, s)$ and the Stickelberger elements $\Theta_{S, \fq}(K/F, s)$ (see \eqref{partzeta2} and \eqref{stickelberger2}) at non-positive integers in terms of a cap-product of an Eisenstein class with certain homology classes. For that we make specific choices for the coprime ideals $\fb$ and $\fc$ of $\cO_F$.  

To fix ideas let $K/F$ be a finite abelian extension with Galois group $G$. Let $\fp$ be a fixed nonarchimedean place of $F$ and let $S$ be a finite and set of nonarchimedean places of $F$ with $\fp\not\in S$ and such that $S'= S\cup \{\fp\}$ contains all places that are ramified in $K/F$. Let $\ff$ be the nonarchimedean part of the conductor of $L/K$. For $\fb$ we choose any non-trivial ideal of $\cO_F$ coprime to $S'$ and for $\fc$ we choose a sufficiently high power of $\fp$. Namely if $\fp^f$ is the exact power of $\fp$ that divides $\ff$ (the case $f=0$ is allowed) then we put $\fc :=\fp^m$ where $m$ is any integer $\ge \max(f, n+1)$. 

For $\sigma\in G$ and $k\in \bZ_{\ge 0}$ we define
\begin{equation*}
\label{homclassgal}
\varrho_{\sigma}^k \, =\, \varrho_{\sigma, S, \fc}^k\,:= \, (\cN^k\cup 1_{\rec^{-1}(\sigma)}) \cap \vartheta \in H_{n-1}(\Ga, \Int_{\loc, b}(1 + \fp^m\cOs, \bZ)(\ep))
\end{equation*}
where $\Si$ denotes again the set of all prime factors of $\fb\fc$ and where 
$1_{\rec^{-1}(\sigma)}\in H^0(\Ga, C((\bA^{\Sigma})^*/ U^{S, \Sigma}, \bZ))$ is the characteristic function of the preimage of $\sigma$ under the reciprocity map $\rec : (\bA^{\Sigma})^*/U^{S, \Sigma}\to G$. 

\begin{coro}
\label{coro:eisparzeta2}
We have 
\begin{equation}
\label{eisparzeta2}
\bEis_{\fb, \fc} \cap \varrho_{\sigma}^k  \, =\, (-1)^{n-1} \Norm(\fb) \, \zeta_{S'}(\sigma, -k) + (-1)^n \Norm(\fb)^{-k} \, \zeta_{S'}(\sigma\sigma_{\fb}, -k).
\end{equation}
\end{coro}

Recall that $\sigma_{\fb}$ denotes the image of $\fb$ under the Artin map $\cI^S\to G$. 

\begin{proof} We choose an ideal $\fm\subseteq \cO_F$ whose set of prime factors is equal to $S$ and so that $\fc\fm$ is a multiple of $\ff$ (hence $K\subseteq F^{\fc\fm}$). If $\rec: \cI^{\fc\fm}/\cP^{\fc\fm} \to \Gal(F^{\fc\fm}/F)$ denotes the reciprocity isomorphism then we have 
\begin{equation*}
\label{homclassgal2}
\varrho_{\sigma}^k \,= \,\sum_{\fA} \varrho_{\fA}^k 
\end{equation*}
where the sum is taken over the ray classes $\fA\in \cI^{\fc\fm}/\cP^{\fc\fm}$ with $\rec(\fA)|_K = \sigma$. Thus the assertion follows from \eqref{partzeta3} and \eqref{eisparzeta}.
\end{proof}

Finally, we give a formula similar to \eqref{eisparzeta2} for the Stickelberger elements $\Theta_{S', \fq}(K/F, -k)$. Recall that $\fq$ denotes an additional place not contained in $S$. Now we choose $\fb:=\fq$ so that $\Si=\{\fp, \fq\}$. Moreover if $q$ denotes the characteristic of the residue field of $\fq$ then we can choose the coefficient ring $R$ to be $\bZ[1/q]$. As in (\cite{dassp}, \S 5.4) by composing the reciprocity map $\rec: (\bA^{\Si})^*/U^{S, \Si}\to G$ with the inclusion $G\hra \bZ[G]$ we view it as an homomorphism 
\begin{equation}
\label{homrec}
\rec: (\bA^{\Si})^*/U^{S, \Si}\lra \bZ[G].
\end{equation}
The fact that $\Ga$ lies in its kernel implies $\rec\in H^0(\Ga, C((\bA^{\Si})^*/U^{S, \Si}, \bZ[G]))$. For $k\in \bZ_{\ge 0}$ we define
\begin{equation}
\label{homclassrec}
\varrho_{K/F}^k\, =\, \varrho_{K/F, S, \fc, \fq}^k\,:= \, (\cN^k\cup (\inv\circ \rec)) \cap \vartheta \in H_{n-1}(\Ga, \Int_{\loc, b}(1 + \fp^m\cOs, \bZ[G])(\ep))
\end{equation}
where $\inv$ denotes the involution $\bZ[G]\to  \bZ[G], \, \sum_{\sigma\in G} n_{\sigma} [\sigma]\mapsto \sum_{\sigma\in G} n_{\sigma} [\sigma^{-1}]$.

\begin{coro}
\label{coro:eisstick}
With $\bEis : = \bEis_{\fq, \fp^m}$ we have
\begin{equation}
\label{eisstick}
\bEis \cap \varrho_{K/F}^k  \, =\, (-1)^n \Norm(\fq)^{-k} \, [\sigma_{\fq}] \cdot \Theta_{S', \fq}(-k).
\end{equation}
\end{coro}

Note that \eqref{eisstick} is an equality in the group ring $\bZ[1/q][G]$.

\begin{proof} Under the canonical isomorphism 
\begin{eqnarray*}
H_{n-1}(\Ga, \Int_{\loc, b}(1 + \fp^m\cOs, \bZ[G])(\ep)) & \cong & H_{n-1}(\Ga, \Int_{\loc, b}(1 + \fp^m\cOs, \bZ)(\ep))\otimes \bZ[G] \\
& = & \bigoplus_{\sigma\in G} H_{n-1}(\Ga, \Int_{\loc, b}(1 + \fp^m\cOs, \bZ)(\ep)) \otimes [\sigma]
\end{eqnarray*}
the class \eqref{homclassrec} decomposes as $\varrho_{K/F}^k = \sum_{\sigma\in G} \varrho_{\sigma^{-1}}^k\otimes [\sigma]$.
Together with \eqref{eisparzeta2} we deduce
\begin{eqnarray*}
\bEis \cap \varrho_{K/F}^k & = & (-1)^n \sum_{\sigma\in G}\left( \Norm(\fq)^{-k}  \zeta_{S'}(\sigma^{-1}\sigma_{\fq}, -k) - \Norm(\fq) \zeta_{S'}(\sigma^{-1}, -k) \right) [\sigma]\\
& = & (-1)^n \Norm(\fq)^{-k}\, [\sigma_{\fq}] \cdot \Theta_{S', \fq}(-k).
\end{eqnarray*}
\end{proof}

\section{Divisibility properties of Stickelberger elements} 
\label{section:lvalues}

As at the end of last section we let $S$ be a finite set of nonarchimedean places of $F$ of cardinality $r$ and let $\Si=\{\fp, \fq\}$ where $\fp$, $\fq$ are two distinct fixed nonarchimedean places of $F$ not contained in $S$. As before we put $\Ga := \cOscu$ and consider $\bZ[1/q]$-coefficients where $q$ is the residue characteristic of $\fq$. We also put $\fc = \fp^m$ (for $m$ sufficiently large). 
The aim of this section is to prove Theorem \ref{theorem:highervanishing}. In order to obtain this refinement of Corollary \ref{coro:eisstick} we work in a more general framework than in the last section, namely we take the cap-product of our adelic Eisenstein class $\Eis:= \Eis_{\fq, \fp^m}$ with certain hyperhomology groups and we work with more general characters than the reciprocity map \eqref{homrec}. 

To begin with let $R$ be an arbitrary $\bZ[1/q]$-algebra and let $\cU$ be an open subgroup of $(\bA^{\Si})^*$ that contains $U^{S, \Si}$. For this data we consider the trilinear map
\begin{eqnarray}
\label{triple}
&& \beta_S: \cD_{\lpol}(1+\fp^m\cOs, \bZ[1/q])(\ep) \times C((\bA^{\Sigma})^*/\cU,  \Int_{\loc, b}(1 + \fp^m\cOs, R))\times C_c((\bA^{\Sigma})^*/U^{\Sigma},  \bZ) \to R\\
&&\hspace{2cm}  (\mu, \Psi, \varphi) \mapsto \beta_S(\mu, \Psi, \varphi) := \int_{1+\fp^m\cOs} \langle \Psi, \varphi \rangle_S(x)\, d\mu(x).\nonumber
\end{eqnarray}
Here $\langle \Psi, \varphi\rangle_S\in \Int_{\loc, b}(1 + \fp^m\cOs, R)$ denotes the image of the pair $(\Psi, \varphi)$ under \eqref{deltapair}. Note that \eqref{triple} is $\Ga$-equivariant, i.e.\ we have 
\begin{equation}
\label{triple2}
\beta_S(\ga \mu, \ga \Psi, \ga \varphi) \, =\, \beta_S(\mu, \Psi, \varphi)
\end{equation}
for every $\ga\in \Ga$. 

Now fix $k\in \bZ_{\ge 0}$ and let
\begin{equation*}
\label{homchi}
\chi: (\bA^{\Si})^*/\cU\lra R^*
\end{equation*}
be a homomorphism. We define  
\begin{equation*}
\label{chinormk}
\chi \cN^k: (\bA^{\Si})^*/\cU \, \lra\,  \Int_{\loc, b}(1 + \fp^m\cOs, R),\quad a \cU \mapsto (\cN^k)(a U^{\Si}) \otimes \chi(a \cU)\end{equation*}
and consider the map \eqref{triple} for $\Psi= \chi\cN^k$ fixed, i.e.\ we consider the pairing
\begin{equation}
\label{kchipair}
\langle \wcdot, \wcdot \rangle_{\chi, k, S}:= \beta_S(\wcdot, \chi \cN^k, \wcdot): \cD_{\lpol}(1+\fp^m\cOs, \bZ[1/q])(\ep) \times C_c((\bA^{\Sigma})^*/U^{\Sigma},  \bZ) \lra R.
\end{equation}
Thus for $\mu\in \cD_{\lpol}(1+\fp^m\cOs, \bZ[1/q])(\ep)$ and $\varphi\in C_c((\bA^{\Sigma})^*/U^{\Sigma},  \bZ)$ we have 
\begin{equation*}
\label{kchipair9}
\langle \mu, \varphi \rangle_{\chi, k, S} :=  \beta_S(\mu, \chi\cN^k, \varphi)\,=\,\int_{1+\fp^m\cOs} \langle \chi \cN^k, \varphi\rangle_S(x) \, d\mu(x).
\end{equation*}
Note that \eqref{triple2} and \eqref{normfunc2} implies that $\langle \ga \mu, \ga \varphi \rangle_{\chi, k, S}=\chi(\ga) \cdot \langle \varphi, \mu \rangle_{\chi, k, S}$ for every $\ga \in \Ga$. Therefore \eqref{kchipair} induces cap-product pairings
\begin{equation*}
\label{kchipair3}
\cap_{\chi, k, S}: H^i(\Ga, \cD_{\lpol}(1+\fp^m\cOs, \bZ[1/q])(\ep))\times H_j(\Ga, C_c((\bA^{\Sigma})^*/U^{\Sigma}, \bZ)) \, \lra  \, H_{j-i}(\Ga, R(\chi))
\end{equation*} 
for $i,j\in \bZ$. In particular for $i=j=n-1$ we can consider the homology class 
\begin{equation*}
\label{eisstick2}
\bEis \cap_{\chi, k, S} \,\vartheta  \in H_0(\Ga, R(\chi))
\end{equation*}

For a place $v$ of $F$ with $v\ne \fp, \fq$ we let $\chi_v: F_v^* \to R^*$ be the local component of $\chi$ at $v$. Our main technical result is

\begin{theorem}
\label{theorem:vanishing}
Let  $\{\fa_v\}_{v\in S\cup S_{\infty}}$ be a collection of ideals of $R$ that satisfies following properties
\medskip

\noi (i) For $v\in S$ the local component $\chi_v$ is unramified modulo $\fa_v$\footnote{By this we mean $\chi_v(u) \equiv 1 \mod \fa_v$ for every $u\in U_v$} and $\chi_v(\varpi_v) \equiv \Norm(v)^k\!\! \mod \fa_v$.
\medskip

\noi (ii) For $v\in S_{\infty}$ we have $\chi_v(-1) \equiv (-1)^k \!\!\mod \fa_v$.
\medskip

\noi Put $\barR = R/\prod_{v\in S\cup S_{\infty}}\fa_v$ and let $\pi: R\to \barR$ be the projection. Then we have 
\begin{equation}
\label{highvanishing}
\pi_*\left(\bEis \cap_{\chi, k, S}\, \vartheta\right) \,= \,0 \quad \text{in} \quad H_0(\Ga, \barR(\chi)).
\end{equation}
\end{theorem}

The proof of Theorem \ref{theorem:vanishing} is given in several steps. It uses ideas introduced in (\cite{spiess1}, \S 3), (\cite{dassp}, \S 3) and  \cite{hirose}. In order to deal with the vanishing of the left hand side of \eqref{highvanishing} modulo the product over the infinite places $\prod_{v\in S_{\infty}}\fa_v$ we need the following\footnote{See also \cite{dassp}, Prop.\ 3.8 for a related argument.}

\begin{lemma}
\label{lemma:archvanish}
The image of the pairing \eqref{kchipair} is contained in $\prod_{v\in S_{\infty}} \fa_v$. 
\end{lemma}

\begin{proof} By shrinking $\cU$ if necessary we may assume that it is of the form $\cU = \cU_f \times U_{\infty}$ where $\cU_f$ is an open subgroup of $(\bA_f^{\Si})^*$ with $U_f^{S, \Si}\subseteq \cU_f \subseteq U_f^{\Si}$. The first key ingredient in the proof is the formula 
\begin{equation}
\label{projfin} 
\langle \Psi, \varphi\rangle_S\, =\, \langle P_{f, \ep}(\Psi), P_f(\varphi)\rangle_{S, f} 
\end{equation} 
for $\Psi \in C((\bA^{\Sigma})^*/\cU,  \Int_{\loc, b}(1 + \fc\cOs, R))$ and $\varphi\in C_c((\bA^{\Sigma})^*/U^{\Sigma},  \bZ)$. We explain the terms on the right hand side of \eqref{projfin}. Firstly, we define the map $P_{f, \ep}$ as 
\begin{eqnarray}
\label{projfep}
&& P_{f, \ep}: C((\bA^{\Sigma})^*/\cU,  \Int_{\loc, b}(1 + \fc\cOs, R))\lra C((\bA_f^{\Sigma})^*/\cU_f,  \Int_{\loc, b}(1 + \fc\cOs, R))(\ep)\\
&& \hspace{2cm}\Psi \mapsto P_{f, \ep}(\Psi)= \sum_{x_{\infty}U_{\infty} \in F_{\infty}^*/U_{\infty}} \ep(x_{\infty}) \cdot \Psi(\wcdot, x_{\infty}U_{\infty})\nonumber
\end{eqnarray}
and the map $P_f$ by
\begin{equation*}
\label{projf2} 
P_f: C_c((\bA^{\Sigma})^*/U^{\Si}, \bZ)\lra C_c((\bA_f^{\Sigma})^*/U_f^{\Si}, \bZ), \quad \varphi \mapsto P_f(\varphi)= \sum_{x_{\infty}U_{\infty} \in F_{\infty}^*/U_{\infty}} \varphi(\wcdot, x_{\infty}U_{\infty}).
\end{equation*}
The pairing 
\begin{equation*}
\label{deltapairfin}
\langle\wcdot , \wcdot \rangle_{S, f}: C((\bA_f^{\Sigma})^*/\cU_f,  \Int_{\loc, b}(1 + \fc\cOs, R))(\ep) \times C_c((\bA_f^{\Sigma})^*/U_f^{\Si}, \bZ) \lra \Int_{\loc, b}(1 + \fp^m\cOs, R)(\ep)
\end{equation*}
is defined similarly to \eqref{deltapair} by using the map $\Delta_{S, f}^{\Sigma}$ instead of $\Delta_S^{\Si}$ in its definition. More concretely we define
\begin{equation}
\label{dasspmap5a}
\wDelta_{S, f}^{\Sigma}: C_c((\bA_f^{\Sigma})^*/\cU_f,  \Int_{\loc, b}(1 + \fp^m \cOs, R))(\ep) \, \lra \, \Int_{\loc, b}(1 + \fp^m \cOs, R)(\ep).
\end{equation}
similar as \eqref{dasspmap5} by replacing in its definition the map \eqref{dasspmap3} with the map \eqref{dasspmap1a}. We then define \eqref{deltapairfin}
as the composite of the obvious pairing 
\[
C((\bA_f^{\Sigma})^*/\cU_f,  \Int_{\loc, b}(1 + \fc\cOs, R)) \times C_c((\bA_f^{\Sigma})^*/U_f^{\Si}, \bZ) \lra 
C_c((\bA_f^{\Sigma})^*/\cU_f,  \Int_{\loc, b}(1 + \fc\cOs, R))
\]
with \eqref{dasspmap5a}.

Secondly, we show that $P_{f,\ep}(\Psi)$ has values in $\Int_{\loc, b}(1 + \fc\cOs, \prod_{v\in S_{\infty}}\fa_v)$. Let $a=(a_v)_{v\not\in\Si\cup S_{\infty}}\in (\bA_f^{\Si})^*$ be an 
idele and let $\fa\in \cI^{\Si}$ be the associated fractional ideal. We have 
\begin{equation*}
\label{projfep2}
P_{f,\ep}(\chi \cN^k)(a \cU_f) \, =\, \Norm(\fa)^{-k} 1_{\fa \cap 1+\fp^m\cOs} \cdot N^k\otimes \left(\chi_f(a\cU_f)\cdot  \sum_{a_{\infty}U_{\infty} \in F_{\infty}^*/U_{\infty}} \chi_{\infty}(a_{\infty})\ep(a_{\infty})^{k+1}\right)
\end{equation*}
where $\chi_f$ (resp.\ $\chi_{\infty}$) denotes the nonarchimedean (resp.\ the archimedean) component of $\chi$. Since 
\[
\sum_{a_{\infty}U_{\infty} \in F_{\infty}^*/U_{\infty}} \chi_{\infty}(a_{\infty})\ep(a_{\infty})^{k+1}=\prod_{v\in S_{\infty}} (1 + \chi_v(-1)(-1)^{k+1}) \in \prod_{v\in S_{\infty}} \fa_v
\]
we get
\begin{equation}
\label{projfep3}
P_{f,\ep}(\chi \cN^k)(a \cU_f) \in \Int_{\loc, b}(1 + \fc\cOs, \bZ) \otimes \prod_{v\in S_{\infty}}\fa_v.
\end{equation}
Now using \eqref{projfin} and \eqref{projfep3} we conclude 
\begin{equation*}
\label{kchipair2}
\langle \mu, \varphi \rangle_{\chi, k, S} \,=\,\int_{1+\fp^m\cOs} \langle P_{f,\ep}(\chi \cN^k), P_f(\varphi)\rangle_{S, f}(x) \, d\mu(x) \in \prod_{v\in S_{\infty}}\fa_v.
\end{equation*}
\end{proof}

For the proof of Theorem \ref{theorem:vanishing} we need further preparation. To begin with we recall the definition of the homology class 
\begin{equation}
\label{thetaS}
\vartheta^S\in H_{n+r-1}(\Ga, C_c((\bA^{S, \Sigma})^*/U^{S, \Sigma},  \bZ)).
\end{equation} 
introduced in (\cite{dassp}, \S 3). Its definition is similar to that of the class \eqref{theta}. One only has to replace the group $E_+$ with the group $E_{S,+}$ of totally positive $S$-units of $F$ (which is free-abelian of rank $n+r-1$) and use the fact that 
\begin{equation*}
\label{thetaS2}
H_{n+r-1}(\Ga, C_c((\bA^{S, \Sigma})^*/U^{S, \Sigma},  \bZ))\,\cong \, C(\cF, \bZ)\otimes H_{n+r-1}(E_{S,+}, \bZ)
\end{equation*}
where now $\cF$ denotes a fundamental domain for the action of $\Ga/E_{S, +}$ on $(\bA^{S, \Sigma})^*/U^{S, \Sigma}$. The class $\vartheta^S$ corresponds to $1_{\cF}\otimes \eta_S$ under this isomorphism where $\eta_S$ denotes a generator of the group $H_{n+r-1}(E_{S,+}, \bZ)$. In order to have a canonical choice for $\eta_S$ we have to fix an ordering $v_1, \ldots, v_r$ of the primes in $S$. 

We will reinterpret the class \eqref{thetaS} as a $\Ga$-hyperhomology class in degree $n-1$. For that we introduce the following $(\bA^{\Si})^*$-modules. Given a subset $S_1\subseteq S$ define 
\begin{equation*}
\label{thetaS3}
\cC_c(S_1, \bZ) := C_c(F_{S_1} \times (\bA^{S_1, \Sigma})^*/U^{S_1, \Sigma},  \bZ)^{U_{S_1}} = C_c(F_{S_1} \times (\bA_f^{S_1, \Sigma})^*\times F_{\infty}^*/U_{\infty}, \bZ)^{U_{S_1} \times U_f^{S_1, \Sigma}}
\end{equation*}
i.e.\ $\cC_c(S_1, \bZ)$ consists of locally constant functions $\varphi: F_{S_1} \times (\bA_f^{S_1, \Sigma})^*\times F_{\infty}^*/U_{\infty}\to \bZ$ with compact support such that $\varphi(u_1 x_1, u_2 x_2, x_3U_{\infty}) = \varphi(x_1, x_2, x_3U_{\infty})$ for every $(x_1, x_2, x_3)\in F_{S_1}\times (\bA^{S_1, \Sigma})^*\times F_{\infty}^*$ and $(u_1, u_2)\in U_{S_1} \times U_f^{S_1, \Sigma}$. Note that if $S_2\subseteq S_1$ then $F_{S_2} \times (\bA_f^{S_2, \Sigma})^*\times F_{\infty}^*/U_{\infty}$ is an open subset of $F_{S_1} \times (\bA_f^{S_1, \Sigma})^*\times F_{\infty}^*/U_{\infty}$. Hence we can view $\cC_c(S_2, \bZ)$ as a submodule of $\cC_c(S_1, \bZ)$, namely $\cC_c(S_2, \bZ)$ consists of those $\varphi\in \cC_c(S_1, \bZ)$ with $\supp(\varphi) \subseteq F_{S_2} \times (\bA_f^{S_2, \Sigma})^*\times F_{\infty}^*/U_{\infty}$. 

\begin{remark}
\label{remark:anncs}
\rm Let $S_1\subseteq S$ be a subset, let $v\in S_1$ and let $\varphi\in \cC_c(S_1, \bZ)$. Then the function $\varphi- [\varpi_v] \varphi$ vanishes at every element $x\in F_{S_1} \times (\bA_f^{S_1, \Sigma})^*\times F_{\infty}^*/U_{\infty}$ whose $v$-component is $=0$, i.e.\ we have $\varphi - [\varpi_v] \varphi \in \cC_c(S_2, \bZ)$ with $S_2 : = S_1\setminus \{v\}$. \footnote{Recall that $\varpi_v$ is a prime element of $\cO_v$.}
Thus if we extend $(\bA^{\Si})^*$-action on $\cC_c(S_1, \bZ)$ to an action of the group ring $\bZ[(\bA^{\Si})^*]$ then we have 
$(1- [\varpi_v]) \cdot \cC_c(S_1, \bZ)\subseteq \cC_c(S_2, \bZ)$. Hence we get 
\begin{equation}
\label{anncs1}
\prod_{v\in S_1} (1 - [\varpi_v]) \cdot \cC_c(S_1, \bZ)\subseteq  \cC_c(\emptyset, \bZ)= C_c((\bA^{\Sigma})^*/U^{\Sigma},  \bZ).
\end{equation}
\enddemo
\end{remark}

We consider the bounded complex of $\Ga$-modules 
\begin{equation*}
\label{thetaS4}
\begin{CD}
\cC_{\bu}: \quad 0 @>>> \cC_0 @>\partial_0 >> \cC_{-1} @>\partial_{-1} >> \ldots @>>> \cC_{1-r} @>\partial_{1-r}>>\cC_{-r}@>>> 0 
\end{CD}
\end{equation*}
defined by 
\begin{equation*}
\label{thetaS5}
\cC_{-i} \, =\, \bigoplus_{\tiny{\begin{array}{c} S_1 \subseteq S\\ \# S_1 = i\end{array}}} \, \cC_c(S_1, \bZ)
\end{equation*}
for $i\in \{0, 1, \ldots, r\}$ and $\cC_i=0$ if $i\not\in \{-r, \ldots, -1, 0\}$. The boundary map $\partial_{-i}: \cC_{-i}\to \cC_{-i-1}$ for $i\in \{0, 1, \ldots, r-1\}$ is defined as follows: if $S_2$ and $S_1$ are subsets of $S$ of cardinality $i$ and $i +1$ respectively then the $(S_2, S_1)$-component of $\partial_{-i}$ is $=0$ if $S_2\not\subset S_1$ and is $(-1)^{\nu} \cdot \incl: \cC_c(S_2, \bZ)\hra \cC_c(S_1, \bZ)$ if $S_2\subseteq S_1$ and $S_1 = \{v_{j_0}, \ldots,  v_{j_i}\}$, $S_2 =  S_1\setminus \{v_{j_{\nu}} \}$ with $1\le j_0< \ldots < j_i \le r$. 

For the homology of $\cC_{\bu}$ and hyperhomology of $\Ga$ with coefficients in $\cC_{\bu}$ we have

\begin{lemma}
\label{homthetaS5}
For every $i\in \bZ$ we have
\medskip

\noi (a) $\sH_i(\cC_{\bu}) \, =\, \left\{ \begin{array}{cc} C_c((\bA^{S, \Sigma})^*/U^{S, \Sigma}, \bZ) & \mbox{if $i = -r$,}\\
                                                      0 & \mbox{otherwise.}
                                                      \end{array}\right.
$
\medskip

\noi (b) $H_i(\Ga, \cC_{\bu}) = H_{i-r}(\Ga, C_c((\bA^{S, \Sigma})^*/U^{S, \Sigma}, \bZ))$.
\end{lemma}

\begin{proof} (a) For a nonarchimedean place $v$ we consider the short exact sequence 
\begin{equation*}
\begin{CD}
 0 @>>>  C_c(F_v^*, \bZ) @> \varphi\mapsto \varphi_! >>  C_c(F_v, \bZ) @> \varphi \mapsto \varphi(0) >>\bZ @>>> 0
\end{CD}
\end{equation*}
(compare \cite{dassp}, \S 3.1). Note that it remains exact if we pass to $U_v$-invariants. Thus if let $\cC_{\bu}^v$ be the complex $[C_c(F_v^*, \bZ)^{U_v} \hra  C_c(F_v, \bZ)^{U_v}]$ concentrated in degree $0$ and $-1$ then its homology vanishes 
except in degree $-1$ where it is $=\bZ$. Now the assertion follows from the fact that $\cC_{\bu}$ is isomorphic to the complex $\left(\bigotimes_{i=1}^r \cC_{\bu}^{v_i}\right) \otimes C_c((\bA^{S, \Sigma})^*/U^{S, \Sigma}, \bZ)$.

(b) There is a homological spectral sequences (see \cite{weibel}, 5.7.8)
\[
E_{ij}^2 = H_i(\Ga, \sH_j(\cC_{\bu})) \, \Longrightarrow\, E_{i+j}= H_{i+j}(\Ga, \cC_{\bu}).
\]
By (a) it degenerates and we have $E_i \cong E_{i-r,r}^2$. 
\end{proof}

By the above Lemma we can view \eqref{thetaS} as a hyperhomology class 
\begin{equation*}
\label{thetaS6}
\vartheta^S\in H_{n-1}(\Ga, \cC_{\bu}) = H_{n+r-1}(\Ga, C_c((\bA^{S, \Sigma})^*/U^{S, \Sigma},  \bZ)).
\end{equation*} 

\begin{remark}
\label{remark:thetaStheta}
\rm Note that $\cC_0 = C_c((\bA^{\Sigma})^*/U^{\Sigma},  \bZ)$. Thus if we view $C_c((\bA^{\Sigma})^*/U^{\Sigma},  \bZ)$ as a complex concentrated in degree $0$ and if $\io: \cC_{\bu} \to C_c((\bA^{\Sigma})^*/U^{\Sigma},  \bZ)$ denotes the forgetful chain map (i.e.\ it is the identity in degree $0$) then it induces a homomorphism 
\begin{equation}
\label{thetaS7}
H_{n-1}(\Ga, \cC_{\bu}) \lra H_{n-1}(\Ga, C_c((\bA^{\Sigma})^*/U^{\Sigma},  \bZ)).
\end{equation} 
Lemma 3.1 of \cite{dassp} can be rephrased by stating that the class $\vartheta^S$ is mapped to $\vartheta$ under the homomorphism \eqref{thetaS7}.
\end{remark}

The following Lemma is the key technical result we use to deal with the places in $S$ in \eqref{highvanishing}.

\begin{lemma}
\label{lemma:nkchiextend}
Let $v\in S$ such that we have
\begin{equation}
\label{kchipair4}
\chi_v(\varpi_v) = |x|_v^{-k} \qquad \forall\, x\in \cO_v, x\ne 0.
\end{equation}
Then we have
\begin{equation*}
\label{kchipair5}
\langle \mu, \varphi \rangle_{\chi, k, S} \, =\, \beta_{S_0}(\mu, [\varpi_v]\cdot (\chi\cN^k), (1- [\varpi_v]) \varphi) 
\end{equation*}
for every $\mu\in \cD_{\lpol}(1+\fp^m\cOs, \bZ[1/q])$ and $\varphi\in C_c((\bA^{\Sigma})^*/U^{\Sigma},  \bZ)$.\footnote{Recall that $\varpi_v$ is a prime element of $\cO_v$.}
\end{lemma}

Note that property \eqref{kchipair4} is equivalent to require that $\chi_v$ is unramified and that $\chi_v(\varpi_v) = \Norm(v)^k$ holds.

\begin{proof} We may assume that $\cU = V \times U^{S_0, \Si}$ for some open subset $V\subseteq U_{S_0}$. By \eqref{deltapair3} we have
\begin{equation}
\label{kchipair6}
\langle \chi\cN^k, \varphi\rangle_S \, =\, \langle  [\varpi_v]\cdot (\chi \cN^k), (1 - [\varpi_v]) \varphi \rangle_{S_0} + \langle \chi\cN^k - [\varpi_v] \cdot \left(\chi\cN^k\right), \varphi \rangle_{S_0}.
\end{equation}
Thus to prove \eqref{kchipair5} it suffices to verify
\begin{equation}
\label{kchipair7}
 \langle \chi\cN^k - [\varpi_v] \cdot \left(\chi\cN^k\right),  \varphi \rangle_{S_0}\, =\, 0
\end{equation}
for every $\varphi\in C_c((\bA^{\Sigma})^*/U^{\Sigma},  \bZ)$. In fact by Remark \ref{remarks:pairchangeS} (b) (and in particular \eqref{deltapair2b}) it is enough to show that 
\begin{equation}
\label{kchipair8}
1_{\fX_{a, V, S_0}} \cdot \left(\chi \cN^k- [\varpi_v] \cdot (\chi \cN^k)\right)(a) \, =\, 0
\end{equation}
for every $a\in (\bA^{\Sigma})^*$. 

Let $\fr$ denote the prime ideal of $\cO_F$ corresponding to the place $v$ and put $\varrho:= \chi_v(\varpi_v)$ so that $\varrho =\Norm(\fr)^k\in \bZ \cdot 1_R \subseteq R$. 
Fix $a = (a_1, a_2) \in (\bA^{\Sigma})^*= F_v^* \times (\bA^{v, \Sigma})^*$ and put $\nu=\ord_v(a_1)$ so that $\chi_v(a_1) = \varrho^\nu$. Let $\chi_2: (\bA^{v, \Sigma})^*\to R^*$ be the composition 
\[
(\bA^{v, \Sigma})^* \cong 1\times (\bA^{v, \Sigma})^*\hra F_v^* \times (\bA^{v, \Sigma})^* =(\bA^{\Sigma})^*\stackrel{\pr}{\lra} (\bA^{\Sigma})^*/\cU \lra R^*.
\] 
We have $\chi(a \cU)= \varrho^\nu\cdot  \chi_2(a_2)$ and by \eqref{normfunc1} 
\begin{equation*}
\label{chinormk2}
\cN(a U^{\Si}) \,=\, \ep(a_{\infty}) \Norm(\fr^{\nu} \fa_2)^{-1} 1_{\fr^{\nu} \fa_2 \cap 1+\fp^m\cOs} \cdot N.
\end{equation*}
Here $\fa_2\in \cI^{\{v\}\cup \Si}$ denotes the fractional ideal associated to $a_2$ and $a_{\infty} \in F_{\infty}^*$ its archimedean component. Because of  $\varrho =\Norm(\fr)^k$ we get
\begin{eqnarray}
\label{chinormk3}
(\chi\cN^k)(a \cU) & = & \ep(a_{\infty})^k \cdot \left( \Norm(\fr^{\nu}\fa_2)^{-k} 1_{\fr^{\nu} \fa_2 \cap 1+\fp^m\cOs} \cdot N^k\right) \otimes (\varrho^{\nu}\cdot  \chi_2(a_2))\\
& = & \ep(a_{\infty})^k \cdot \left( \Norm(\fr^{\nu-1}\fa_2)^{-k} 1_{\fr^{\nu} \fa_2 \cap 1+\fp^m\cOs} \cdot N^k\right) \otimes (\varrho^{\nu-1}\cdot  \chi_2(a_2)).
\nonumber
\end{eqnarray}
Applying \eqref{chinormk3} to the adele $a[\varpi_v]^{-1}$ yields
\begin{equation}
\label{chinormk4}
\left([\varpi_v] \cdot (\chi\cN^k)\right)(a) \, = \, \ep(a_{\infty})^k \cdot \left( \Norm(\fr^{\nu-1}\fa_2)^{-k} 1_{\fr^{\um-1} \fa_2 \cap 1+\fp^m\cOs} \cdot N^k\right) \otimes (\varrho^{\nu-1}\cdot  \chi_2(a_2)).
\end{equation}
Therefore by combining \eqref{chinormk3} and \eqref{chinormk4} we obtain 
\begin{equation*}
\label{chinormk5}
\left(\chi\cN^k- [\varpi_v] \cdot (\chi\cN^k)\right)(a)\, =\, -\ep(a_{\infty})^k \cdot \left( \Norm(\fr^{\nu-1}\fa_2)^{-k} 1_{\fX} \cdot N^k\right) \otimes (\varrho^{\nu-1}\cdot  \chi_2(a_2))
\end{equation*}
where $\fX = (\fr^{m-1}\fa_2 \setminus \fr^{\nu}\fa_2) \cap (1+\fp^m\cOs)$. Now \eqref{kchipair8} follows from
\begin{equation*}
1_{\fX_{a, V, S_0}} \cdot \left(\chi \cN^k- [\varpi_v] (\chi \cN^k)\right)(a) = \pm \left( \Norm(\fr^{m-1}\fa_2)^{-k} 1_{\fX_{a, V, S_0} \cap \fX} \cdot N^k\right) \otimes (\varrho^{\nu-1}  \chi_2(a_2)) = 0
\end{equation*}
where we have used $\fX_{a, V, S_0} \cap \fX\subseteq \fr^{\nu}\fa_2 \cap (\fr^{m-1}\fa_2 \setminus \fr^{\nu}\fa_2) = \emptyset$.
\end{proof}

Now the key observation is that if $S_1\subseteq S$ is a subset such that \eqref{kchipair4} holds for every $v\in S_1$ then the pairing \eqref{kchipair}
extends canonically to a pairing $\cD_{\lpol}(1+\fp^m\cOs, \bZ[1/q])(\ep) \times \cC_c(S_1, \bZ) \to R$. Namely, we define a pairing 
\begin{equation}
\label{kchipext}
\langle \wcdot, \wcdot \rangle_{\chi, k, S_1, S}: \cD_{\lpol}(1+\fp^m\cOs, \bZ[1/q])(\ep) \times \cC_c(S_1,  \bZ) \lra R
\end{equation}
by 
\begin{equation*}
\label{kchipext2}
\langle \mu, \varphi \rangle_{\chi, k, S_1, S} := \beta_{S\setminus S_1}(\mu, \prod_{v\in S_1}  [\varpi_v]\cdot (\chi\cN^k), \prod_{v\in S_1} (1 - [\varpi_v]) \varphi) 
\end{equation*}
for $\mu\in \cD_{\lpol}(1+\fp^m\cOs, \bZ[1/q])(\ep)$ and $\varphi\in \cC_c(S_1,  \bZ)$. Note that it is well-defined by \eqref{anncs1}. The pairing \eqref{kchipext} has the following properties

\begin{lemma}
\label{lemma:kchiextc}
Let $S_1\subseteq S$ be a subset such that \eqref{kchipair4} holds for every $v\in S_1$. 
\medskip

\noi (a) If $S_2\subseteq S_1$ then we have 
\begin{equation}
\label{kchipext3}
\langle \mu, \varphi \rangle_{\chi, k, S_1, S}= \langle \mu, \varphi \rangle_{\chi, k, S_2, S}
\end{equation}
for every $\mu\in \cD_{\lpol}(1+\fp^m\cOs, \bZ[1/q])(\ep)$ and $\varphi\in \cC_c(S_2,  \bZ)$. 
\medskip

\noi (b) The pairing \eqref{kchipext} is $\Ga$-equivariant, i.e.\ we have 
\[
\langle \ga \mu, \ga \varphi \rangle_{\chi, k, S_1, S}\, =\, \chi(\ga) \cdot \langle \varphi, \mu \rangle_{\chi, k, S_1, S}
\] 
for every $\ga \in \Ga$, $\mu\in \cD_{\lpol}(1+\fp^m\cOs, \bZ[1/q])(\ep)$ and $\varphi\in \cC_c(S_1,  \bZ)$.
\medskip

\noi (c) Let $\{\fa_v\}_{v\in S_{\infty}}$ be a collection of ideals of $R$ such that condition (ii) of Thm\ \ref{theorem:vanishing} holds.
Then the image of the pairing \eqref{kchipext} is contained in $\prod_{v\in S_{\infty}} \fa_v$. 
\end{lemma}

Note that (a) for $S_2=\emptyset$ implies in particular that the pairing \eqref{kchipext} extends \eqref{kchipair}, i.e.\ we have 
\[
\langle \mu, \varphi \rangle_{\chi, k, S_1, S}= \langle \mu, \varphi \rangle_{\chi, k, S}
\]
for every $\mu\in \cD_{\lpol}(1+\fp^m\cOs, \bZ[1/q])(\ep)$ and $\varphi\in C_c((\bA^{\Sigma})^*/U^{\Sigma},  \bZ)$. 

\begin{proof} For (a) it suffices to consider the case $S_1 = S_2 \cup \{w\}$, $w\not\in S_2$. Note that if $S_2=\emptyset$, $S_1=\{w\}$ then 
this is just the statement of Lemma \ref{lemma:nkchiextend} above. The proof there can be easily adapted to the more general case so we only sketch the argument. Put $S_3=S\setminus S_1$, $S_4=S\setminus S_2= S_3 \cup\{w\}$, $\Psi := \prod_{v\in S_2}  [\varpi_v]\cdot (\chi\cN^k)$ and $\varphi': =  \prod_{v\in S_2} (1 - [\varpi_v]) \varphi$. Similar to \eqref{kchipair6}, Remark \ref{remarks:pairchangeS} (b) yields
\begin{equation*}
\label{kchipext4}
\langle \Psi, \varphi'\rangle_{S_4} \, =\, \langle  [\varpi_w] \Psi, (1 - [\varpi_w]) \varphi'\rangle_{S_3} +  \langle (1- [\varpi_w]) \Psi , \varphi' \rangle_{S_3}
\end{equation*}
and \eqref{kchipext3} follows once we have established that $\langle (1- [\varpi_w]) \Psi , \varphi' \rangle_{S_3}=0$. This can be seen by an almost identical argument as \eqref{kchipair7}.

The second assertion follows again from \eqref{triple2} and \eqref{normfunc2} and the third can be proved by the same arguments as Lemma \ref{lemma:archvanish}. Note that one only has to verify that \eqref{projfep3} holds if we replace the function $\chi\cN^k$ with $\prod_{v\in S_1}  [\varpi_v]\cdot (\chi\cN^k)$. 
\end{proof}

Now assume that $\{\fa_v\}_{v\in S\cup S_{\infty}}$ be a collection of ideals of $R$ such that conditions (i) and (ii) of Thm.\ \ref{theorem:vanishing} are satisfied. Recall that the second condition is equivalent to require that we have $\chi_v(x) \equiv |x|_v^{-k}\!\! \mod \fa_v$ for every $x\in \cO_v$, $x\ne 0$ and $v\in S$. For $v\in S$ we put $R(v) := R/\fa_v$ and let $\chi(v): (\bA^{\Si})^*/\cU\to R(v)^*$ be the reduction of $\chi$ modulo $\fa_v$. More generally, for a non-empty subset $S_1\subseteq S$ we put $R(S_1) := R/\fa_{S_1}$ where $\fa_{S_1} = \sum_{v\in S_1} \fa_v$ and let $\chi^{(S_1)}: = \pr\circ \chi: (\bA^{\Si})^*/\cU\stackrel{\chi}{\lra} R^*\stackrel{\pr}{\lra} R(S_1)^*$ be the reduction of $\chi$ modulo $\fa_{S_1}$. For $S_1=\emptyset$ we put $R(\emptyset) = R$ and $\chi^{(\emptyset)}=\chi$. 

We consider the bounded complex 
\begin{equation}
\label{rmoda}
\begin{CD}
\cR_{\bu}: \quad 0 @>>> \cR_0 @>\partial_0 >> \cR_{-1} @>\partial_{-1} >> \ldots @>>> \cR_{1-r} @>\partial_{1-r}>>\cR_{-r}@>>> 0 
\end{CD}
\end{equation}
defined by 
\begin{equation*}
\label{rmoda2}
\cR_{-i} \, =\, \bigoplus_{\tiny{\begin{array}{c} S_1 \subseteq S\\ \# S_1 = i\end{array}}} \, R(S_1)
\end{equation*}
for $i\in \{0, 1, \ldots, r\}$ and $\cR_i=0$ if $i\not\in \{-r, \ldots, -1, 0\}$. The boundary map $\partial_{-i}: \cR_{-i}\to \cR_{-i-1}$ for $i\in \{0, 1, \ldots, r-1\}$ is defined as follows. Let $S_2$ and $S_1$ are subsets of $S$ of cardinality $i$ and $i +1$ respectively. The $(S_2, S_1)$-component of $\partial_{-i}$ is $=0$ if $S_2\not\subset S_1$. It is $(-1)^{\nu} \cdot \pr$, where $\pr: R(S_2)\to R(S_1)$ is the natural projection, if $S_2\subseteq S_1$ and $S_1 = \{v_{j_0}, \ldots,  v_{j_i}\}$, $S_2 =  S_1\setminus \{v_{j_{\nu}} \}$ with $1\le j_0< \ldots < j_i \le r$. Alternatively, the complex \eqref{rmoda} may be defined as follows. For $v\in S$ let $\cR_{\bu}^v$ be the complex $[R\stackrel{\pr}{\lra} R/\fa_v]$ concentrated in degree $0$ and $-1$. Then we have 
\begin{equation*}
\label{rmoda3}
\cR_{\bu} \, =\, \cR_{\bu}^{v_1}\otimes_R \ldots \otimes_R \cR_{\bu}^{v_1}.
\end{equation*}

We can extend the pairing \eqref{kchipair} to a pairing $\cD_{\lpol}(1+\fp^m\cOs, \bZ[1/q])\times \cC_{\bu} \to \cR_{\bu}$, i.e.\ to a chain map 
\begin{equation}
\label{kchipext5}
\Xi_{\bu}: \cD_{\lpol}(1+\fp^m\cOs, \bZ[1/q])(\ep) \otimes  \cC_{\bu} \lra \cR_{\bu} 
\end{equation}
by defining $\Xi_{-i}$ for $i\in \{0, 1, \ldots, r\}$ as the direct sum (the sum over all $S_1\subseteq S$ with $\# S_1=i$) of the maps 
\begin{equation*}
\label{kchipext6}
\cD_{\lpol}(1+\fp^m\cOs, \bZ[1/q])(\ep) \otimes \cC_c(S_1,  \bZ) \lra R(S_1), \quad \mu \otimes \varphi \mapsto \langle\mu, \varphi \rangle_{\chi^{(S_1)}, k, S_1, S}
\end{equation*}
In fact by Lemma \ref{lemma:kchiextc} (c) we see that the image of \eqref{kchipext5} is contained in $\cR_{\bu} \otimes_R (\prod_{v\in S_{\infty}} \fa_v)$. That \eqref{kchipext5} is indeed a chain map follows from Lemma \ref{lemma:kchiextc} (a). If we take into account the $\Ga$-actions then Lemma \ref{lemma:kchiextc} (b) implies that 
\[
\Xi_{\bu}: \cD_{\lpol}(1+\fp^m\cOs, \bZ[1/q])(\ep) \otimes  \cC_{\bu} \to \cR_{\bu}\otimes_R \left(\prod_{v\in S_{\infty}}\fa_v\right)(\chi)
\] 
is a chain map of $\Ga$-modules. It therefore induces cap-product pairings
\begin{equation*}
\label{kchipext7}
\cap_{\chi, k, S}: H^i(\Ga, \cD_{\lpol}(1+\fp^m\cOs, \bZ[1/q])(\ep))\times H_j(\Ga, \cC_{\bu}) \, \lra  \, H_{j-i}(\Ga, \cR_{\bu}\otimes_R \left(\prod_{v\in S_{\infty}}\fa_v\right)(\chi))
\end{equation*} 
for $i,j\in \bZ$. 

Now we make specific choices for $R$, $\chi$ and the ideals $\fa_v$, $v\in S\cup S_{\infty}$ namely we choose them to be ``universal''. Firstly, choose an ideal $\fm\subseteq \cO_F$ whose set of prime factors is contained in $S$ and let $\chi$ be the universal homomorphism from $(\bA^{\Si})^*/U_{\fm}^{\Si}$ into the group of units of an $\bZ[1/q]$-algebra. More precisely, we consider the $\bZ[1/q]$-group algebra $R=R^{\univ}:= \bZ[1/q][(\bA^{\Si})^*/U_{\fm}^{\Si}]$ and let $\chi$ be the obvious homomorphism 
\begin{equation*}
\label{chiuniv}
\chi= \chi^{\univ}:(\bA^{\Si})^*/U_{\fm}^{\Si}\lra R^*, \quad xU_{\fm}^{\Si}\mapsto [xU_{\fm}^{\Si}].
\end{equation*}
For $v\in S\cup S_{\infty}$ we choose for $\fa_v$ to be the smallest ideal of $R$ such that condition (i) resp.\ (ii) of Thm.\ \ref{theorem:vanishing} holds. Specifically, if $m_v$ denotes the exponent of $v\in S$ occurring the prime decomposition of the ideal $\fm$ then we set
\begin{equation*}
\label{rvuniv}
\wfa_v\, :=\, \left\{\begin{array}{cc} \ker(\bZ[1/q][F_v^*/U_v^{(m_v)}]\stackrel{\pr}{\lra} \bZ[1/q][F_v^*/U_v]/([\varpi_v]- \Norm(v)^k)) & \mbox{if $v\in S$,}\\
\ker(\bZ[1/q][F_v^*/U_v] \stackrel{\pr}{\lra} \bZ[1/q][F_v^*/U_v]/\left([(-1) U_v] + (-1)^{k+1} [U_v]\right))& \mbox{if $v\in S_{\infty}$}
\end{array}\right.
\end{equation*}
and put
\begin{equation*}
\label{rvuniv2}
\fa_v=\fa_v^{\univ} \,:= \,\wfa_v \otimes \bZ[1/q][(\bA^{\Si, v})^*/U_{\fm}^{\Si, v}].
\end{equation*}
Note that $\fa_v$ is indeed contained in $R$ since $\bZ[1/q][(\bA^{\Si, v})^*/U_{\fm}^{\Si, v}]$ is a free $\bZ[1/q]$-algebra.
With these choices for $R$ and $\fa_v$ we obtain 

\begin{lemma}
\label{homologycr}
For every $i\in \bZ$ we have
\medskip

\noi (a)  $\sH_i\left(\cR_{\bu}\otimes_R \left(\prod_{v\in S_{\infty}}\fa_v\right)\right) \, =\, \left\{ \begin{array}{cc} \prod_{v\in S\cup S_{\infty}}\fa_v & \mbox{if $i = 0$,}\\
                                                      0 & \mbox{otherwise.}
                                                      \end{array}\right.$
                                                      \medskip
                                                      
\noi (b) $H_i(\Ga, \cR_{\bu}\otimes_R \left(\prod_{v\in S_{\infty}}\fa_v\right)) = H_i(\Ga, (\prod_{v\in S\cup S_{\infty}}\fa_v)(\chi))$.
\end{lemma}

\begin{proof} (a) For $v\in S_{\infty}$ the quotient $(\bZ[1/q][F_v^*/U_v])/\wfa_v$ is isomorphic to 
\[
\bZ[1/q][F_v^*/U_v]/\left([(-1) U_v] + (-1)^{k+1} [U_v]\right)\cong \bZ[1/q],
\]
hence it is in particular a flat $\bZ[1/q]$-algebra. It follows that we have
\[
\prod_{v\in S_{\infty}} \fa_v \,=\,  \bZ[1/q][(\bA^{\Si, \infty})^*/U_{\fm}^{\Si, \infty}] \otimes_{\bZ[1/q]} \bigotimes_{v\in S_{\infty}}\wfa_v
\]
For $v\in S$ we put $R_v:= \bZ[1/q][F_v^*/U_v^{(m_v)}]$, $\barR_v:= \bZ[1/q][F_v^*/U_v]/([\varpi_vU_v]-\Norm(v)^k)$ and let $\pi_v: R_v \to \barR_v$ be the canonical projection so that $\wfa_v$ is its kernel. If we denote by $\wcR_{\bu}^v$ the complex $[R_v \stackrel{\pi_v}{\lra}\barR_v]$ concentrated in degree $0$ and $-1$ and put $M^S:=\bZ[1/q][(\bA_f^{S, \Si})^*/U_f^S]\otimes_{\bZ[1/q]} (\bigotimes_{v\in S_{\infty}}\wfa_v)$ then we have 
\begin{equation*}
\label{rmoda4}
\cR_{\bu}\otimes_R \left(\prod_{v\in S_{\infty}}\fa_v\right)\ \, =\, \cR_{\bu}^{v_1}\otimes_{\bZ[1/q]} \ldots \otimes_{\bZ[1/q]} \cR_{\bu}^{v_1}\otimes_{\bZ[1/q]} M^S.
\end{equation*}
Note that $R_v$, $\barR_v$ and $M^S$ are flat as $\bZ[1/q]$-modules. Indeed, $R_v$ and $R^S$ are free as $\bZ[1/q]$-modules and we have $\barR_v\cong \bZ[1/q][T^{\pm 1}]/(T-1)\cong \bZ[1/q]$ if $k=0$ and $\barR_v\cong \bZ[1/q][T^{\pm 1}]/(T-\Norm(v)^k)\cong \bZ[1/q, 1/\Norm(v)]$ if $k>0$. It follows 
\begin{eqnarray*} 
\sH_i\left(\cR_{\bu}\otimes_R \left(\prod_{v\in S_{\infty}}\fa_v\right)\right) & = & \left\{ \begin{array}{cc} \wfa_{v_1}\otimes_{\bZ[1/q]} \ldots \otimes_{\bZ[1/q]} \wfa_{v_r}\otimes_{\bZ[1/q]} M^S& \mbox{if $i = 0$,}\\
                                                      0 & \mbox{otherwise.}
                                                      \end{array}\right.\\
&=& \left\{ \begin{array}{cc}  \prod_{v\in S\cup S_{\infty}}\fa_v & \mbox{if $i = 0$,}\\
                                                      0 & \mbox{otherwise.}
                                                      \end{array}\right.
                                                      \end{eqnarray*}
(b) follows immediately from (a) using the homological spectral sequences (\cite{weibel}, 5.7.8).
\end{proof}

\begin{proof}[Proof of Theorem \ref{theorem:vanishing}] We first prove the assertion for $R=R^{\univ}$, $\chi=\chi^{\univ}$, $\{\fa_v\}_{v\in S\cup S_{\infty}}=\{\fa_v^{\univ}\}_{v\in S\cup S_{\infty}}$ and $\pi=\pi^{\univ}: R^{\univ} \to \barR^{\univ} = R^{\univ}/\prod_{v\in S\cup S_{\infty}}\fa_v^{\univ}$. Consider the commutative diagram
\begin{equation*}
\label{hypercap}
\begin{CD}
H_{n-1}(\Ga, \cC_{\bu}) @>\eqref{thetaS7} >> H_{n-1}(\Ga, C_c((\bA^{\Sigma})^*/U^{\Sigma}, \bZ)) @. \\
@VV \xi \mapsto  \bEis \cap_{\chi, k, S} \xi V @VV \zeta \mapsto \bEis \cap_{\chi, k, S} \zeta V \\
H_0(\Ga,  \cR_{\bu}\otimes_R \left(\prod_{v\in S_{\infty}} \fa_v\right)(\chi)) @>>> H_0(\Ga, R(\chi)) @> \pi_* >> H_0(\Ga, \barR(\chi))
\end{CD}
\end{equation*}
The first lower horizontal map is induced by the composite of the forgetful map 
\[
\cR_{\bu}\otimes_R \left(\prod_{v\in S_{\infty}}\fa_v\right)\lra \cR_0 \otimes_R \left(\prod_{v\in S_{\infty}}\fa_v\right)= \prod_{v\in S_{\infty}}\fa_v
\]
with the inclusion $\prod_{v\in S_{\infty}}\fa_v\hra R$. Therefore by Lemma \ref{homologycr} (b) the lower sequence of the diagram can be identified with the exact sequence 
\[
\begin{CD}
H_0(\Ga, (\prod_{v\in S\cup S_{\infty}}\fa_v)(\chi))@> \io_* >> H_0(\Ga, R(\chi)) @> \pi_* >>H_0(\Ga, \barR(\chi))
\end{CD}
\]
where $\io: \prod_{v\in S\cup S_{\infty}}\fa_v \hra R$ is the inclusion. Using Remark \ref{remark:thetaStheta} we obtain
\[
\pi_*\left(\bEis \cap_{\chi, k, S}\, \vartheta\right)\, =\, \pi_* \circ \io_*\left(\bEis \cap_{\chi, k, S}\, \vartheta^S \right)=(\pi\circ \io)_*\left(\bEis \cap_{\chi, k, S}\, \vartheta^S \right)\,=\, 0.
\]
Now we consider the case of an arbitrary ring $R$, character $\chi: (\bA^{\Si})^*/\cU\to R^*$ and collection of ideals $\{\fa_v\}_{v\in S\cup S_{\infty}}$ of $R$ satisfying the conditions (i) and (ii). We choose an ideal $\fm\subseteq \cO_F$ whose set of prime factors is contained in $S$ and such that $U_{\fm}^{\Si}\subseteq \cU$. Then $\chi$ induces a $\bZ[1/q]$-algebra homomorphism $X: R^{\univ}\to R$ such that $\chi= X\circ \chi^{\univ}$ and $X(\fa_v^{\univ}) \subseteq \fa_v$ for all $v\in S\cup S_{\infty}$. Therefore $X$ induces a homomorphism $\barX: \barR^{\univ}\to \barR$ and standard properties of the cap-product imply 
\[
\pi_*\left(\bEis \cap_{\chi, k, S}\, \vartheta\right) \,= \, \barX_* \left(\pi^{\univ}_*\left(\bEis \cap_{\chi^{\univ}, k, S}\, \vartheta\right)\right) \,= \,0 .
\]
\end{proof}

Now we turn to the proof of Thm.\ \ref{theorem:highervanishing}. Therefore in the following $K/F$ denotes a finite abelian extension with Galois group $G$ and $S$ a finite set of nonarchimedean places of $F$ that contains all places that are ramified in $K/F$. We fix a place $\fp\in S$ and an additional nonarchimedean place $\fq$ of $F$ not contained in $S$.

\begin{prop}
\label{prop:highervanishing2}
For every $k\in \bZ_{\ge 0}$ we have 
\begin{equation}
\label{highvanishing3}
\Theta_{S, \fq}(K/F, -k)\in \left(\prod_{v\in S \cup S_{\infty}, v\ne \fp} \cI_v^{(k)}\right)\otimes \bZ[1/q].
\end{equation}
Here $q$ denotes again the residue characteristic of $\fq$. 
\end{prop}

\begin{proof} We apply Thm.\ \ref{theorem:vanishing} to the ring $R= \bZ[1/q][G]$, the character $\chi= \inv\circ \rec : (\bA^{\fp, \fq})^*/\cU\lra \bZ[1/q][G]$ (where $\cU$ is a sufficiently small open subgroup of $U^{\fp, \fq}$ that contains $U^{S, \fq}$) and the collection of ideals $\fa_v\otimes \bZ[1/q]:= \cI_v^{(k)}\otimes \bZ[1/q]$ for $v\in S\setminus \{\fp\}$.
Thus by \eqref{highvanishing} and Cor.\ \ref{coro:eisstick} we have \footnote{Note that we have changed the notation; what is called $S$ and $S'$ in loc.\ cit.\ is called $S\setminus \{\fp\}$ and $S$ here.}
\[
\pi_*\left((-1)^n \Norm(\fq)^{-k} \, [\sigma_{\fq}] \cdot \Theta_{S, \fq}(-k)\right) \, =\, \pi_*\left(\bEis \cap_{\chi, k, S\setminus \{\fp\}}\, \vartheta\right)\, =\, 0\quad \text{in}\quad H_0(\Ga, \barR(\chi))
\]
where $\pi: \bZ[1/q][G] \to \barR:=\left(\bZ[G]/\prod_{v\in S \cup S_{\infty}, v\ne \fp} \cI_v^{(k)}\right)\otimes \bZ[1/q]$ is the projection. Note that because $\chi$ is trivial on $\Ga$, we have $\pi_* = \pi: R= H_0(\Ga, R(\chi)) \to H_0(\Ga, \barR(\chi)) = \barR$, hence 
\[
(-1)^n \Norm(\fq)^{-k} \, [\sigma_{\fq}] \cdot \Theta_{S, \fq}(-k)\,\in \ker(\pi)= \prod_{v\in S \cup S_{\infty}, v\ne \fp} \cI_v^{(k)}\otimes \bZ[1/q].
\]
Since $\Norm(\fq)$ and $[\sigma_{\fq}]$ are units in $\bZ[1/q][G]$ we conclude $\Theta_{S, \fq}(K/F, -k)\in \prod_{v\in S \cup S_{\infty}, v\ne \fp} \cI_v^{(k)}$.
\end{proof}

For $s\in \bC$ put $\delta_T(s):= \prod_{\fq\in T} (1-\Norm(\fq)^{1-s} [\sigma_{\fq}^{-1}])\in \bC[G]$ so that $\Theta_{S, T}(K/F, s)\, =\,\delta_T(s) \Theta_S(K/F, s)$. The following Lemma generalizes Lemmas 62 and 63 in \cite{hirose}.\footnote{In \cite{hirose} only the case $k=0$ is considered.}  

\begin{lemma}
\label{lemma:annh0}
For $k\in \bZ_{\ge 0}$ we have
\medskip

\noi (a) $\delta_T(-k) \in \Ann_{\bZ[G]}\left(\bigoplus_{v\in T_K} H^0(k(v), (\bQ/\bZ)'(k+1))\right)$.
\medskip

\noi (b) Let $m:= \#(H^0(K, \bQ/\bZ(k+1)))$, let $X\subseteq \Spec \cO_F[1/m]$ be an nonempty open subset disjoint from $S$ and let $\cJ(X)\subseteq \bZ[G]$ be the ideal generated by the set $\{1-\Norm(\fq)^{k+1} [\sigma_{\fq}^{-1}] \mid \fq\in X\}$. Then, 
\[
\Ann_{\bZ[G]}(H^0(K, \bQ/\bZ(k+1))) \, = \, \cJ(X).
\]
\end{lemma}

\begin{proof} (a) Let $\fq\in T$ and let $v$ be a place of $K$ above $\fq$. Since the Frobenius $\sigma_{\fq}$ acts by multiplication with $\Norm(\fq)^{k+1}$ on $H^0(k(v), (\bQ/\bZ)'(k+1))$ we have $(1-\Norm(\fq)^{1 + k} [\sigma_{\fq}^{-1}])\cdot H^0(k(v), (\bQ/\bZ)'(k+1))= 0$.
Because $\delta_T(-k)$ is a multiple of $1-\Norm(\fq)^{1 + k} [\sigma_{\fq}^{-1}]$ we conclude that
\[
\delta_T(-k) \cdot \left(\bigoplus_{v\in T_K} H^0(k(v), (\bQ/\bZ)'(k+1))\right)\, =\, 0.
\]
(b) The action of $G$ on $H^0(K, \bQ/\bZ(k+1))$ is given by a character $\psi: G \to (\bZ/m\bZ)^*$. By replacing $K$ with the fixed field of $\ker(\psi)$ if necessary we may assume that $\psi: G \to (\bZ/m\bZ)^*$ is injective.
The character $\psi$ induces a ring homomorphism $\Psi: \bZ[G]\to \bZ/m\bZ$. If $\fq\in X$ then $\psi(\sigma_{\fq}) = \Norm(\fq)^{k+1}\!\!\pmod m$ hence $1-\Norm(\fq)^{k+1} [\sigma_{\fq}^{-1}]\in \ker(\Psi)$. In particular we have 
\begin{equation}
\label{decomplaw1}
\Norm(\fq)^{k+1} \equiv 1 \mod m \quad \Longleftrightarrow \quad \psi(\sigma_{\fq}) = 1  \quad \Longleftrightarrow \quad \sigma_{\fq} = 1.
\end{equation}
Note that the ring homomorphism $\bZ\to \bZ[G]/\cJ(X)$ is surjective with non-trivial kernel. Indeed, for any $\sigma\in G$ there exists $\fq\in X$ with $\sigma_{\fq}=\sigma$, so we have $[\sigma]-\Norm(\fq)^{k+1} = [\sigma_{\fq}] (1-\Norm(\fq)^{k+1} [\sigma_{\fq}^{-1}]) \in J(X)$. Let $m'\in \bZ_{\ge 1}$ with $m' \bZ=\bZ \cap \cJ(X)$. We will show $m'=m$ which implies 
\[
\Ann_{\bZ[G]}(H^0(K, \bQ/\bZ(k+1)))\,= \,\ker(\Psi)\,= \,\cJ(X).
\]
Firstly, since $m'\in \cJ(X) \subseteq \ker(\Psi)$ we have $m\mid m'$. Note that for $\fq\in X$ we have 
\begin{equation}
\label{decomplaw2}
\Norm(\fq)^{k+1}\equiv 1 \mod m \quad \Longrightarrow \quad \Norm(\fq)^{k+1}\equiv 1 \mod m'.
\end{equation}
Indeed, if $1-\Norm(\fq)^{k+1}\in m\bZ$ then $\sigma_{\fq} =1$ by \eqref{decomplaw1} and therefore $1- \Norm(\fq)^{k+1} = (1-\Norm(\fq)^{k+1} [\sigma_{\fq}^{-1}])\in J(X) \cap \bZ=m'\bZ$. To prove $m= m'$ we consider the action of $\cG_F$ on the group $\mu_{m'}(\barQ)^{\otimes^{k+1}}$. It is given by a character $\psi': \cG_F \to (\bZ/m'\bZ)^*$. We denote by $K'$ the fixed field of $\ker(\psi')$, by $G'$ the Galois group of $K'/F$ and by $\Psi': \bZ[G'] \to \bZ/m'\bZ$ the ring homomorphism induced by $\psi'$. As before the ideal $J(X)'\subseteq \bZ[G']$ generated by the set $\{1-\Norm(\fq)^{k+1} [(\sigma'_{\fq})^{-1}] \mid \fq\in X\}$ is contained in $\ker(\Psi')$ (here $\sigma'_{\fq}$ denotes the Frobenius at $\fq$ in $G'$). Now by \eqref{decomplaw1} and \eqref{decomplaw2} a prime $\fq\in X$ that splits completely in $K/F$ splits completely in $K'/F$ as well. Since $K\subseteq K'$ this implies by \v{C}ebotarev's density theorem that $K=K'$, hence $m'=m$.
\end{proof}

\begin{proof}[Proof of Thm.\ \ref{theorem:highervanishing}]  
It suffices to show 
\begin{equation}
\label{anntheta}
\cJ:= \Ann_{\bZ[G]}(H^0(K, \bQ/\bZ(k+1)))  \subseteq \cI : = \{ x\in \bZ[G]\mid x\cdot \Theta_S(K/F, -k)\in \prod_{v\in S \cup S_{\infty}, v\ne \fp} \cI_v^{(k)}\}.
\end{equation}
Indeed, the injectivity of \eqref{galoisres} together with Lemma \ref{lemma:annh0} (a) implies that $\delta_T(-k)\in \cJ$. Hence \eqref{anntheta} yields 
\[
\Theta_{S, T}(K/F, -k)\, =\, \delta_T(-k) \cdot \Theta_S(K/F, -k) \in \prod_{v\in S \cup S_{\infty}, v\ne \fp} \cI_v^{(k)}.
\]
To prove \eqref{anntheta} note that by \eqref{highvanishing3}, for every $\fq\in X$ there exists 
$n_{\fq}\in \bZ_{\ge 0}$ with 
\begin{equation}
\label{anntheta2}
\Norm(\fq)^{n_{\fq}}(1-\Norm(\fq)^{k+1} [\sigma_{\fq}^{-1}]) \in \cI. 
\end{equation}
By Lemma \ref{lemma:annh0} (b) for any nonempty open subset $X\subseteq \Spec \cO_F[1/m]$ disjoint from $S$ we have 
$\cJ(X) = \cJ$. In fact since $\bZ[G]$ is noetherian the Lemma implies that we can choose two disjoint finite sets $X_1, X_2\subseteq \Spec \cO_F[1/m]\setminus S$ such that $\cJ(X_1) = \cJ = \cJ(X_2)$ and so that $M_1:= \prod_{\fq\in X_1} \Norm(\fq)^{n_{\fq}}$ and $M_2:= \prod_{\fq\in X_2} \Norm(\fq)^{n_{\fq}}$ are coprime. Now \eqref{anntheta2} implies $M_i\cdot \cJ = M_i \cdot \cJ(X_i) \subseteq \cI$ for $i=1,2$ hence $\cJ = \gcd(M_1, M_2) \cdot \cJ \subseteq \cI$. 
\end{proof}

\addappendix \paragraph{Topologies and Sheaves} Let $X$ be a set. We consider a collection $\fOpen(X)$ of subsets of $X$ with the property that $\emptyset, X\in \fOpen(X)$ and that for any pair $U, V\in \fOpen(X)$ we have $U\cap V, U\cup V\in \fOpen(X)$. The pair $(X, \fOpen(X))$ will be called a {\it pre-site}. 
Elements of $\fOpen(X)$ will be called open subsets of $X$. For an arbitrary subset $X'\subseteq X$ we let $\fOpen(X')$ be the set of subsets of $X'$ of the form $U\cap X'$ with $U\subseteq X$ open. Clearly $(X', \fOpen(X'))$ is again a pre-site. Let $(X, \fOpen(X))$ and $(Y, \fOpen(Y))$ be pre-sites. A pre-continuous morphism $f: (X, \fOpen(X))\to (Y, \fOpen(Y))$ is map $f:X\to Y$ such we have $f^{-1}(V) \in \fOpen(X)$ for every $V\in \fOpen(Y)$. 

Let $R$ be a ring. A contravariant functor $\cF:\fOpen(X)\to \Mod_R$ with $\cF(\emptyset) =0$ will be called an {\it $R$-presheaf on $(X, \fOpen(X))$}. Here we view $\fOpen(X)$ as a category; its set of objects are the open subsets of $X$ and the only morphisms are inclusion maps. A morphism of $R$-presheaves is a morphism of functors. As usual if $V\subseteq U$ are open subsets of $X$ then the image of the inclusion $V\hra U$ under $\cF$ will be denoted by $\cF(U) \to \cF(V), s\mapsto s|_V$. The category of $R$-presheaves on $(X, \fOpen(X))$ will be denoted by $\PSh(X, \fOpen(X); R)$ (or by $\PSh(X, R)$ for short, if it is clear which collection of subsets $\fOpen(X)$ of $X$ we consider). It is an abelian category with enough injectives. If $(Y, \fOpen(Y))$ is another pre-site and if $f: X\to Y$ is pre-continuous then for $\cF\in \PSh(X, R)$ we define $f_*(\cF)\in \PSh(Y, R)$ as usual 
by $f_*(\cF)(V):=\cF(f^{-1}(V))$ for every open $V\subseteq Y$.

We collect certain facts and notions from (\cite{johnstone}, Ch.\ C 2, \S 1), (\cite{ks1}, \S\S 16-18) and (\cite{stacks}, 7.47) adapted to the specific framework that is relevant to us. Firstly, we recall the notion of a coverage and of a site (as defined in \cite{johnstone}, Ch.\ C 2, Def.\ 2.1.1).

\begin{df}
\label{df:coverage} 
(a) Let $(X, \fOpen(X))$ be a pre-site. A coverage $\fC$ on $\fOpen(X)$ is function assigning to each open subset $U$ of $X$ a collection $\fC(U)$ of families $\{U_i\}_{i\in I}$ of elements in $\fOpen(U)$ called ($\fC$-)coverings of $U$. It has the following property
\medskip

\noi (C) For open subsets $U$, $V$ of $X$ and a covering $\{U_i\}_{i\in I}$ of $U$, the family $\{U_i\cap V\}_{i\in I}$ is a covering of $U\cap V$.
\medskip

\noi (b) The triple $(X, \fOpen(X), \fC)$ consisting of a pre-site $(X, \fOpen(X))$ and a coverage $\fC$ on $\fOpen(X)$ will be called a site. \footnote{Note that a site as defined above is usually not a site in the sense (\cite{ks1}, 17.2.1) and (\cite{stacks}, Tag 00VH).}
\medskip

\noi (c) A continuous morphism $f:(X, \fOpen(X), \fC_X)\to (Y, \fOpen(Y), \fC_Y)$ between sites is a map $f:X\to Y$ that is pre-continuous (i.e.\ we have $f^{-1}(V) \in \fOpen(X)$ for every open subset $V\subseteq Y$) and is cover-preserving (i.e.\ for every $V \in \fOpen(Y)$ and a covering $\cV=\{V_i\}_{i\in I}$ of $V$ the collection of pre-images $f^{-1}(\cV):= \{f^{-1}(V_i)\}_{i\in I}$ is a covering of $f^{-1}(V)$).
 \end{df}

Note that a topological space $X$ defines in a canonical way a site in the above sense which -- by abuse of notation -- will be denoted by $X$ as well. Moreover a map between topological spaces $f: X\to Y$ is continuous if and only if it is continuous as a morphism between sites. 

Let $X=(X, \fOpen(X), \fC)$ be a site. An $R$-sheaf on $X$ is an $R$-presheaf that satisfies the sheaf property for $\fC$-coverings. More precisely we have

\begin{df}
\label{df:sheafsite}
An $R$-sheaf $\cF$ on $X=(X, \fOpen(X), \fC)$ is an $R$-presheaf on $(X, \fOpen(X))$ such that for every open subset $U\subseteq X$ and every $\fC$-covering $\{U_i\}_{i\in I}$ of $U$ the sequence of $R$-modules 
\[
\begin{CD}
0 @>>> \cF(U) @> s\mapsto (s|_{U_i})_{i\in I} >> \prod_{i\in I} \cF(U_i) @> (s_i)_{i\in I}\mapsto ({s_i}|_{U_i\cap U_j}- {s_j}|_{U_i\cap U_j})_{(i,j)\in I^2} >> \prod_{(i,j)\in I^2} \cF(U_i\cap U_j) 
\end{CD}
\]
is exact. 
\end{df}

A homomorphism of $R$-sheaves on $X$ is a morphism of presheaves. The category of $R$-sheaves on $X$ will be denoted  by $\Sh(X, R)=\Sh(X, \fOpen(X), \fC; R)$. 

\begin{prop}
\label{prop:sheafnice}
(a) Let $X$ be a site. The category $\Sh(X, R)$ is $R$-linear, abelian and has enough injectives. Moreover it satisfies the axiom (AB5). 
\medskip

\noi (b) Let $f:X\to Y$ be a continuous morphism between sites. The functor 
\begin{equation}
\label{fimage}
f_*: \Sh(X, R)\lra \Sh(Y, R), \quad \cF \mapsto f_*(\cF)
\end{equation}
is well-defined and admits an exact right adjoint
\begin{equation}
\label{finv}
f^*: \Sh(Y, R)\lra \Sh(X, R),\quad \cG \mapsto f^*(\cG).
\end{equation}
\end{prop}

For the proof we need some preparation. Let $U\subseteq X$ be open. Recall that a subset $\cS$ of $\fOpen(U)$ is called a {\it sieve over $U$} if for every pair of open subsets $U$, $V$ of $X$ with $U\in \cS$ and $V\subseteq U$ we have $V\in \cS$. A coverage $\fC$ on $\fOpen(X)$ is called {\it sifted} if every covering is a sieve. 
Given an arbitrary coverage $\fC$ on $\fOpen(X)$ there is an associated sifted coverage $\barfC$ which defines the same category of $R$-sheaves. We recall its definition. Given an open subset $U \subseteq X$ and a collection $\cU=\{U_i\}_{i\in I}$ of open subsets of $U$ the set of open subsets of $U$
\begin{equation*}
\label{gensieve} 
\cS(\cU) :=\, \{V\in \fOpen(X)\mid \exists i\in I, V \subseteq U_i\}
\end{equation*} 
is a sieve. It is called {\it the sieve generated by $\cU$}. If $\fC$ be a coverage on $\fOpen(X)$ then $\barfC$ is the coverage defined by $\barfC(U) :=\, \{\cS(\cU)\mid \cU\in \fC(U)\}$
for every open $U\subseteq X$. By (\cite{johnstone}, Ch.\ C 2, Lemma 2.1.3) we have $\Sh(X, \fOpen(X), \fC; R)= \Sh(X, \fOpen(X), \barfC; R)$. We recall

\begin{df}
\label{df:gtop}
(a) A sifted coverage $\fC$ on $\fOpen(X)$ is called {\it Grothendieck topology} if the following holds 
\medskip

\noi (GT1) For every open $U\subseteq X$ the maximal sieve over $U$ (i.e.\ the sieve $\fOpen(U)$) is a $\fC$-covering.
\medskip

\noi (GT2) If $\cU$ is a $\fC$-covering of $U$ and if $\cS$ is another sieve over $U$ such that $\cS_{V}:=\{W\in \fOpen(V)\mid W\in \cS\}$ is a $\fC$-covering of $V$ for every $V\in \cU$ then $\cS$ is an $\fC$-covering of $U$. 
\medskip

\noi A triple $(X, \fOpen(X), \fC)$ consisting of a pre-site $(X, \fOpen(X))$ and a Grothendieck coverage $\fC$ on $\fOpen(X)$ will be called a Grothendieck site. \footnote{Note that a Grothendieck site is a site in the sense of (\cite{ks1}, 17.2.1) as well as in the sense of (\cite{stacks}, Tag 00VH).}
\medskip

\noi (b) Let $X=(X, \fOpen(X), \fC_X)$ and $Y=(Y, \fOpen(Y), \fC_Y)$ be Grothendieck sites. 
A Grothendieck continuous morphism $f:X\to Y$ is a map of the underlying sets that is pre-continuous and so that for every open subset $V \subseteq Y$ and $\fC_Y$-covering $\cV$ of $V$ the sieve $\cS(f^{-1}(\cV))$ generated by $f^{-1}(\cV)= \{f^{-1}(W)\mid W\in \cV\}$ is a $\fC_X$-covering of $f^{-1}(V)$. 
\end{df}

By (\cite{johnstone}, Ch.\ C 2, Prop.\ 2.1.9) given a coverage $\fC$ on $\fOpen(X)$ there is a smallest Grothendieck coverage $\wfC$ on $\fOpen(X)$ containing $\fC$ (i.e.\ we have $\fC(U) \subseteq \wfC(U)$ for every open $U\subseteq X$) and for this coverage we have $\Sh(X, R)=\Sh(X, \fOpen(X), \fC; R)= \Sh(X, \fOpen(X), \wfC; R)$. By (\cite{ks1}, 9.6.2 and 18.1.6) we conclude in particular that the category $\Sh(X, R)$ is abelian with enough injectives (see \cite{ks1}, 9.6.2 and 18.1.6). This proves the first part of Prop.\ \ref{prop:sheafnice}. For (b) we need the following 

\begin{lemma}
\label{lemma:contGcont}
Let $f: (X, \fOpen(X), \fC_X)\to (Y, \fOpen(Y), \fC_Y)$ be a continuous morphism between sites. Then $f$ is a Grothendieck continuous morphism, when viewed as a morphism between the Grothendieck sites $(X, \fOpen(X), \wfC_X)\to (Y, \fOpen(Y), \wfC_Y)$.
\end{lemma}

\begin{proof} We define a coverage $\fD$ on $\fOpen(Y)$ by letting $\fD(V)$ for $V\in \fOpen(Y)$ consists of all sieves $\cS'$ over $V$ such that the sieve $\cS(f^{-1}(\cS'))$ over $f^{-1}(V)$ generated by $f^{-1}(\cS')=\{f^{-1}(W)\mid W\in \cS'\}$ is a $\wfC_X$-covering of $f^{-1}(V)$. One easily verifies that $\fD$ is a Grothendieck coverage that contains $\barfC_Y$ hence also $\wfC_Y$.
\end{proof}

\begin{proof}[Proof of Prop.\ \ref{prop:sheafnice} (b)] By the Lemma and (\cite{ks1}, 17.5.1) we have $f_*(\cF) \in \Sh(Y, R)$ 
for $\cF\in \Sh(X, R)$. This proves that \eqref{fimage} is well-defined and that it admits a right adjoint \eqref{finv}. Since the map $f^{-1}: \fOpen(Y)\to \fOpen(X), V\mapsto f^{-1}(V)$ -- when viewed as a functor -- is clearly left exact (compare \cite{ks1}, Def.\ 3.3.1) the functor \eqref{finv} is exact by (\cite{ks1}, 17.5.2 (iv)).
\end{proof}

\begin{remark}
\label{remark:generators} 
\rm The key step in the proof that the category $\Sh(X, R)$ has enough injectives given in \cite{ks1} is to show that $\Sh(X, R)$ has a system of generators. In fact in (\cite{ks1}, 18.1.6) it is shown that for every open subset $U\subseteq X$ there exists a sheaf $R_U\in \Sh(X, R)$ with the property 
\begin{equation*}
\label{generator}
\Hom_{\Sh(X, R)}(R_U, \cF)\,\cong\, \cF(U)
\end{equation*}
for every $\cF\in \Sh(X, R)$, i.e.\ the sheaf $R_U$ represents the functor $\Sh(X, R)\to \Mod_R, \cF\mapsto \cF(U)$. The family of sheaves $\{R_U\}_{U\in \fOpen(X)}$ is then a system of generators of $\Sh(X, R)$.
\end{remark}

\begin{coro}
\label{coro:constantsheaf}
The functor of global sections 
\begin{equation}
\label{globalsection}
H^0(X, \wcdot): \Sh(X, R)\lra \Mod_R, \quad \cF \mapsto H^0(X, \cF):= \cF(X)
\end{equation}
has an exact left adjoint
\begin{equation*}
\label{constantsheaf}
\Mod_R\lra \Sh(X, R), \quad M\mapsto \uM_X.
\end{equation*}
We call $\uM_X$ the sheaf on $X$ associated to the $R$-module $M$. 
\end{coro}

\begin{proof} The category of sites has a final object $\pt$, the one-pointed topological space. Let $f:X\to \pt$ be the unique morphism. For $\pt$ we have $\Sh(\pt, R) = \Mod_R$ and $f_*$ corresponds to 
the functor \eqref{globalsection} under this identification. Therefore the assertion is special case of Prop.\ \ref{prop:sheafnice} (b), namely the case of he morphism $f:X\to \pt$.
\end{proof}

Recall that the sheaf cohomology groups $H^{\bu}(X, \cF)$ are defined as the right derived functors of \eqref{globalsection}. More precisely for $i\in \bZ_{\ge 0}$ the $i$-th right derived functor of \eqref{globalsection} will be denoted by
\begin{equation*}
\label{cohomology}
H^i(X, \wcdot): \Sh(X, R)\lra \Mod_R, \quad \cF \mapsto H^i(X, \cF).
\end{equation*}
We note that if $f:X\to Y$ is a continuous morphism between sites then by Prop.\ \ref{prop:sheafnice} (b) the functor $f_*$ is left exact and preserves injectives. As in the case of usual topological spaces this implies that there exists a Leray spectral sequence 
\begin{equation}
\label{lerayss}
E_2^{rs} = H^r(Y, R^s f_* \cF) \, \Longrightarrow\, E^{r+s}=H^{r+s}(X, \cF)
\end{equation}
for every $R$-sheaf $\cF$ on $X$.

\paragraph{$G$-equivariant sheaf} Let $G$ be a group and let $X=(X, \fOpen(X), \fC_X)$ be a site equipped with a continuous $G$-action $G\times X\to X, (g, x)\mapsto g\cdot x$. By that we mean that the underlying set $X$ is equipped with a $G$-action $G\times X\to X, (g, x)\mapsto g\cdot x$ and that $g\cdot : X\to X, x\mapsto g\cdot x$ is a continuous morphism of sites for every $g\in G$. 

\begin{df}
\label{df:equisheaf}
Let $R$ be a ring and let $X$ be a site equipped with a continuous $G$-action.
\medskip

\noi (a) A $G$-equivariant $R$-sheaf $\cF$ on $X$ is an $R$-sheaf together with a collection of isomorphisms $\rho_{g, \cF}: \cF\to g_*\cF$ such that $(g_1)_*(\rho_{g_2, \cF}) \circ \rho_{g_1, \cF}=\rho_{g_2\cdot g_1, \cF}$ for all $g_1, g_2\in G$. 
\medskip

\noi (b) A morphism of $G$-equivariant $R$-sheaves $\alpha:\cF\to \cG$ is a morphism of sheaves that commutes -- in the obvious sense -- with the isomorphisms $\rho_{g, \cF}, \rho_{g, \cG}$ for every $g\in G$. 
\medskip

\noi (c) The category of $G$-equivariant $R$-sheaves on $X$ will be denoted by $\Sh(X, G, R)$. 
\end{df}

The collection of isomorphisms in (a) will be called the {\it $G$-action} on $\cF$. Note that it induces a $G$-action on the $R$-module of global sections $H^0(X, \cF)= \cF(X)$. 

\begin{prop}
\label{prop:injequiv}
(a) The category $\Sh(X, G, R)$ is an $R$-linear abelian category which satisfies axiom (AB5). 
\medskip

\noi (b) The forgetful functor 
\begin{equation}
\label{forget}
\Sh(X, G, R)\lra \Sh(X, R)
\end{equation} 
is exact and has an exact left adjoint. In particular it preserves injectives.
\medskip

\noi (c) The category $\Sh(X, G, R)$ has enough injectives. 
\medskip

\noi (d) The functor of global sections $H^0(X, \wcdot): \Sh(X, G, R)\to \Mod_{R[G]}$ has an exact left adjoint.
\end{prop}

\begin{proof} (a) and the exactness of \eqref{forget} follow easily from Prop.\ \ref{prop:sheafnice} (a) and from the fact that the kernel and cokernel of a morphism of $G$-equivariant $R$-sheaves $\alpha:\cF\to \cG$ in the category $\Sh(X, R)$ carry again a canonical $G$-action. The left adjoint 
\begin{equation}
\label{indgamma}
\Ind^G: \Sh(X, R)\lra \Sh(X, G, R)
\end{equation}
of \eqref{forget} can be defined as follows. For $\cF\in \Sh(X, R)$ the underlying sheaf of $\Ind^{G}(\cF)$ is given by 
\begin{equation*}
\label{indgamma2}
\Ind^{G}(\cF) = \bigoplus_{g\in G} g_*(\cF).
\end{equation*}
For $g\in G$ we define the isomorphism 
\begin{equation*}
\label{indgamma3}
\rho_{g}:\Ind^{G}(\cF)= \bigoplus_{h\in G} h_*(\cF) \lra g_*(\Ind^{G}(\cF))=\bigoplus_{h\in G} (gh)_*(\cF)\end{equation*}
by just permuting the individual summands $h_*(\cF)$ of $\Ind^{G}(\cF)$ indexed by $h\in G$ via left multiplication with $g^{-1}$, i.e.\ the summand $h_*(\cF)$ is mapped identically to the summand $g_*((g^{-1}h)_*(\cF))$ of $g_*(\Ind^{G}(\cF))$. Clearly the functor \eqref{indgamma} is exact.

(c) It suffices to see that $\Sh(X, G, R)$ admits generators. For an open subset $U\subseteq X$ let $R_U$ be the sheaf representing the functor $\Sh(X, R)\to \Mod_R, \cF\mapsto \cF(U)$ (see Remark \ref{remark:generators}). Then $\Ind^G R_U$ represents the functor $\Sh(X, G, R)\to \Mod_R, \cF\mapsto \cF(U)$. This implies that the family $\{\Ind^G R_U\}_{U\subseteq \hX\open}$ is a system of generators of $\Sh(X, G, R)$.

(d) follows immediately from Cor.\ \ref{coro:constantsheaf} using the fact that a $G$-action on a $R$-module $M$ induces a $G$-action on the sheaf $\uM_X$ associated to $M$. 
\end{proof}

\begin{remark}
\label{remark:twists}
\rm In section \ref{section:latspaces} we have considered equivariant sheaves twisted by a character $\chi:G \to R^*$. To explain this notion let $(\cF, (\rho_{g, \cF})_{g\in G})$ be a $G$-equivariant $R$-sheaf. We define the {\it twisted $G$-equivariant $R$-sheaf} $\cF(\chi)$ by $\cF(\chi)= \cF$ and $\rho_{g, \cF(\chi)}:= \chi(g) \cdot \rho_{g, \cF}$ for $g\in G$.
\end{remark}

We recall the definition of equivariant sheaf cohomology and the spectral sequence linking it to ordinary sheaf cohomology. For $\cF\in \Sh(X, G, R)$ taking the $G$-invariant elements of $H^0(X, \cF)$ define a left exact functor 
\begin{equation}
\label{sections}
\Sh(X, G, R) \lra \Mod_R, \quad \cF \mapsto H^0(X, \cF)^G.
\end{equation}
Its right derived functors will be denoted by 
\begin{equation*}
\label{rsections}
H^i(X, G, \wcdot): \Sh(X, G, R) \lra \Mod_R, \quad \cF \mapsto H^i(X, G, \cF).
\end{equation*}
The cohomology groups $H^{\bu}(X, G, \cF)$ are called the equivariant cohomology of $X$ with coefficient in the $G$-equivariant $R$-sheaf $\cF$. 

\begin{prop}
\label{prop:covss}
For $\cF \in \Sh(X, G, R)$ there exists a spectral sequence
\begin{equation}
\label{hsscomp}
E_2^{rs} = \Ext_{R[G]}^r(R, H^s(X, \cF)) \, \Longrightarrow\, E^{r+s}=H^{r+s}(X, G, \cF).
\end{equation}
\end{prop}

\begin{proof} The functor \eqref{sections} factors as 
\begin{equation}
\label{sections2}
\begin{CD}
\Sh(X, G, R) @> \cF \mapsto H^0(X, \cF) >> \Mod_{R[G]} @> M \mapsto \Hom_{R[G]}(R, M) >> \Mod_R.
\end{CD}
\end{equation}
By Prop.\ \ref{prop:injequiv} (d) the first functor has an exact left adjoint, hence it preserves injectives. Therefore there exists a Grothendieck spectral sequence corresponding to the decomposition \eqref{sections2} of the functor \eqref{sections}. Finally by Prop.\ \ref{prop:injequiv} (b) for the $E_2$-terms of said spectral sequence we have $E_2^{rs} = \Ext_{R[G]}^r(R, H^s(X, \cF))$. 
\end{proof}

Let $f: X\to Y$ be a continuous morphism between sites. For $\cF\in \Sh(Y, R)$ passing in the adjunction map $\cF\to f_*(f^*(\cF))$ to global sections yields a homomorphism
\begin{equation}
\label{pullbacksections}
H^0(Y, \cF)\lra H^0(X, f^*(\cF)).
\end{equation} 
As in (\cite{iversen}, Ch.\ II.5) the exactness of $f^*$ (see Prop.\ \ref{prop:sheafnice} (b)) implies that \eqref{pullbacksections} extends to a morphism of $\delta$-functors
\begin{equation}
\label{apppullbackcoh}
H^i(Y, \cF)\lra H^i(X, f^*(\cF))\qquad i\ge 0, \,\cF\in  \Sh(Y, R).
\end{equation}
Similarly, if $X$ and $Y$ are equipped with a continuous $G$-action and if $f$ is $G$-equivariant and if $\cF$ is a $G$-equivariant sheaf on $Y$ then \eqref{pullbacksections} is a homomorphism of $R[G]$-modules hence passing in \eqref{pullbacksections} to $G$-invariants
leads to homomorphism $H^0(Y, G, \cF)\to H^0(X, G, f^*(\cF))$ that  extends to a morphism of $\delta$-functors
\begin{equation}
\label{apppullbackcoh2}
H^i(Y, G, \cF)\lra H^i(X, G, f^*(\cF)) \qquad i\ge 0, \,\cF\in  \Sh(Y, G, R).
\end{equation}
as well.

\paragraph{Some homological algebra}
Let $F:\cA\to \cB$ and $G:\cB\to\cC$ be left exact additive functors between abelian categories and assume that $\cA$ and $\cB$ have enough injectives. There exists obvious morphisms
\begin{eqnarray}
\label{edge1}
&& e_A^{(1), i}: R^i G(F(A))\lra R^i(G\circ F)(A)\qquad \forall \,A\in \cA, i\ge 0,\\
&& e_A^{(2), i}: R^i(G\circ F)(A) \lra G(R^i F(A))\qquad \forall \,A\in \cA, i\ge 0.
\label{edge2}
\end{eqnarray}
If there exists a Grothendieck spectral sequence associated to the composition of functors $G\circ F: \cA\to \cC$ then \eqref{edge1} are just edge morphisms. Also the following Lemma would be a simple consequence of said spectral sequence.
However we will apply it in situations where the existence of the Grothendieck spectral sequence is unclear.

\begin{lemma}
\label{lemma:nogrss}
Let $n\in \bZ_{\ge 1}$ such that $R^i F(A)=0$ for all $i=1, \ldots, n$. Then the morphism \eqref{edge1} is an epimorphism for $i=0, 1, \ldots, n$ and 
\begin{equation}
\label{edge2a}
e_A^{(2), n+1}: R^{n+1}(G\circ F)(A) \lra G(R^{n+1} F(A))
\end{equation}
is an isomorphism. In particular, if additionally we have $F(A)=0$ then $R^i(G\circ F)(A)=0$ for all  $i=0, 1, \ldots, n$.
\end{lemma}

\begin{proof} We prove the assertion by induction on $n$. So assume that $n\ge 1$ and that the assertion holds for $n-1$. Let $0\lra A\lra I\lra A'\lra 0$ be a short exact sequence in $\cA$ with $I$ being injective. Then by assumption the sequence $0\lra F(A)\lra F(I)\lra F(A')\lra 0$ is exact as well and we have $R^i F(A') \cong R^{i+1} F(A) =0$ for $i=1, \ldots, n-1$. The assertion regarding \eqref{edge1} follows by a diagram chase in  
\begin{equation*}
\label{edge3}
\begin{CD}
\ldots R^i G(F(A)) @>>> R^i G(F(I))@>>> R^i G(F(A'))@>>> R^{i+1} G(F(A))\ldots\\
@VV e_A^{(1), i} V @VV e_I^{(1), i} V@VV e_{A'}^{(1), i}V@VV e_A^{(1), i+1} V\\
\ldots R^i(G\circ F)(A) @>>> R^i(G\circ F)(I)@>>> R^i(G\circ F)(A')@>>> R^{i+1}(G\circ F)(A)\ldots
\end{CD}
\end{equation*}
using the fact that $R^i(G\circ F)(I)=0$ for $i\ge 1$ and that $e_{A'}^i$ is an epimorphism for $i=1, \ldots, n-1$. 
Moreover that \eqref{edge2a} is an isomorphism follows from the commutativity of the diagram
\begin{equation*}
\label{edge4}
\begin{CD}
R^n(G\circ F)(A') @> \cong >> R^{n+1}(G\circ F)(A)\\
@VV e_A^{(2), n} V @VV e_A^{(2), n+1} V\\
G(R^n F(A')) @>\cong >> G(R^{n+1} F(A)).
\end{CD}
\end{equation*}
\end{proof}

\end{document}